\documentclass[12pt]{article}

\usepackage{amsmath,amssymb}
\usepackage{graphicx}
\usepackage{times}
\usepackage[margin=0.78in]{geometry}
\usepackage{color}
\usepackage{epsf}
\usepackage{rotating}
\usepackage{supertabular}
\usepackage{comment}
\usepackage{lastpage}
\usepackage{fancyhdr}

\begin{document}
\title{\boldmath\bf The  Complex Gradient Operator {\LARGE \\ and the {\LARGE $\mathbb{CR}$}-Calculus}}%
\author{{\it June 25, 2009} \\ \ \\ \textsf{Ken Kreutz--Delgado} \smallskip \\
 Electrical and Computer Engineering \\
Jacobs School of Engineering \\ University of California, San
Diego \\ \ \\
{\small VERSION UCSD-ECE275CG-S2009v1.0} \\ \medskip \emph{\small
Copyright $\copyright$ 2003-2009, All Rights Reserved}}
\date{}
\maketitle



\section{Introduction}

Often signals and system parameters are most conveniently
represented as complex-valued vectors. This occurs, for example, in
array processing \cite{VT02}, as well as in communication systems
\cite{GS01} when processing narrowband signals using the
\emph{equivalent complex baseband} representation \cite{He95}.
Furthermore, in many important applications one attempts to optimize
a scalar \emph{real}-valued measure of performance over the complex
parameters defining the signal or system of interest. This is the
case, for example, in LMS adaptive filtering where complex filter
coefficients are adapted on line. To effect this adaption one
attempts to optimize the performance measure by adjustments of the
coefficients along its stochastic gradient direction \cite{Ha96,AS03}.

\medskip
However, an often confusing aspect of complex LMS adaptive
filtering, and other similar gradient-based optimization procedures,
is that the partial derivative or gradient used in the adaptation of
complex parameters is \emph{not} based on the standard complex
derivative taught in the standard mathematics and engineering
complex variables courses \cite{Ch96}-\cite{Fl72}, which exists if
and only if a function of a complex variable $z$ is \emph{analytic}
in $z$.\footnote{I.e., \emph{complex}-analytic.} This is because a nonconstant \emph{real}-valued function of
a complex variable is \emph{not} (complex) analytic and therefore is
\emph{not} differentiable  in the standard textbook
complex-variables sense.

\hfill
\newpage

Nonetheless, the same real-valued function viewed as a function of
the real-valued real and imaginary \emph{components} of the complex
variable can have a (real) gradient when partial derivatives
 are
taken with respect to those two (real) components. In this way we
can shift from viewing the real-valued function as a
non-differentiable mapping between $\mathbb{C}$ and $\mathbb{R}$ to
treating it as a differentiable mapping between $\mathbb{R}^2$ and
$\mathbb{R}$. Indeed, the modern graduate-level textbook  in complex
variables theory by Remmert \cite{RR91} continually and easily
shifts back and forth between the real function $\mathbb{R}^2
\rightarrow \mathbb{R} \, \text{\footnotesize or} \, \mathbb{R}^2$
perspective and the complex function $\mathbb{C} \rightarrow
\mathbb{C}$ perspective of a  complex or real scalar-valued
function,
\[ f(z) = f(r) = f(x,y), \] of a complex variable $z = x + j \, y$,
\[ z \in \mathbb{C} \Leftrightarrow r =
\begin{pmatrix} x \\ y \end{pmatrix} \in \mathbb{R}^2 . \]
In particular,  when optimizing a real-valued function of a complex
variable $z = x + j \, y$ one can work with the equivalent real
gradient of the function viewed as a mapping from $\mathbb{R}^2$ to
$\mathbb{R}$ in lieu of a nonexistent  complex derivative
\cite{Br83}. However, because the
 real gradient perspective
arises within a complex variables framework, a direct reformulation
of the problem to the real domain is awkward. Instead, it greatly
simplifies derivations if one can represent the real gradient as a
redefined, new \emph{complex gradient} operator. As we shall see
below, the complex gradient is an extension of the standard complex
derivative to non-complex analytic functions.

\medskip
Confusing the issue is the fact that there is no one unique way to
consistently define a ``complex gradient'' which applies to
(necessarily non-complex-analytic) real-valued functions of a complex variable, and authors
do not uniformly adhere to the same definition. Thus it is often
difficult to resolve questions about the nature or derivation of the
complex gradient by comparing authors. Given the additional fact
that typographical errors seem to be rampant these days, it is
therefore reasonable to be skeptical of the algorithms provided in
many textbooks--especially if one is a novice in these matters.

\medskip
An additional source of confusion arises from the fact that the
derivative of a function with respect to a vector can be
alternatively represented as a row vector or as a column vector when
a space is Cartesian,\footnote{I.e., is Euclidean with identity
metric tensor.} and both representations can be found in the
literature.
In this note we carefully distinguish
between the complex \emph{cogradient} operator (covariant derivative operator \cite{TF01}), which is a \emph{row
vector} operator, and the associated \emph{complex gradient}
operator which is a \emph{vector} operator which gives the direction
of steepest ascent of a real scalar-valued function.

\medskip
Because of the constant back-and-forth shift between a real function
(``{\footnotesize $\mathbb{R}$}-calculus'') perspective and a
complex function (``{\footnotesize $\mathbb{C}$}-calculus'')
perspective which a careful analysis of nonanalytic complex
functions requires \cite{RR91}, we refer to the mathematics
framework underlying the derivatives given in this note as a
``{{\footnotesize{$\mathbb{CR}$}}-calculus}.'' In the following, we
start by reviewing some of the properties of standard univariate
analytic functions, describe the
{{\footnotesize{$\mathbb{CR}$}}-calculus} for univariate nonanalytic
functions, and then develop a multivariate second order
{{\footnotesize{$\mathbb{CR}$}}-calculus} appropriate for
optimizing scalar real-valued cost functions of a complex
parameter vector. We end the note with some examples.

\section{The Derivative of a Holomorphic Function}

Let $z = x + j y$, for $x,y$ real, denote a complex number and let
\[ f(z) = u(x,y) + j\,  v(x,y)\] be a general complex-valued function
of the complex number $z$.\footnote{\label{foot:a}Later, in Section
\ref{sec:cpd}, we will interchangeably alternate between this
notation and the more informative notation $f(z,\bar{z})$. Other
useful representations are $f(u,v)$ and $f(x,y)$. \emph{In this
section} we look for the (strong) conditions for which $f:z \mapsto
f(z) \in \mathbb{C}$ is differentiable as a mapping $\mathbb{C}
\rightarrow \mathbb{C}$ (in which case we say that $f$ is
\emph{{\footnotesize $\mathbb{C}$}-differentiable}), but in
subsequent sections we will admit the \emph{weaker condition} that
$f: (x,y) \mapsto (u,v)$ be differentiable as a mapping
$\mathbb{R}^2 \rightarrow \mathbb{R}^2$ (in which case we say that
$f$ is \emph{{\footnotesize $\mathbb{R}$}-differentiable}); see
Remmert \cite{RR91} for a discussion of these different types of
differentiability.} In standard complex variables courses it is
emphasized that for the complex derivative,
\[ f'(z) = \lim\limits_{\Delta z \rightarrow 0} \frac{f(z + \Delta z) -
f(z)}{\Delta z} , \] to exist in a meaningful way it must be
\emph{independent} of the direction with which $\Delta z$ approaches
zero in the complex plane. \emph{This is a very strong condition} to
be placed on the function $f(z)$. As noted in an introductory
comment from the textbook by Flanigan \cite{Fl72}:

\begin{quote}\small\textsf
You will learn to appreciate the difference between a complex
analytic function (roughly a complex-valued function $f(z)$ having a
complex derivative $f'(z)$) and the real functions $y=f(x)$ which
you differentiated in calculus. Don't be deceived by the similarity
of the notations $f(z)$, $f(x)$. The complex analytic function
$f(z)$ turns out to be much more special, enjoying many beautiful
properties not shared by the run-of-the-mill function from ordinary
real calculus.  The reason [ $\cdots$ ] is that $f(x)$ is merely
$f(x)$ whereas the complex analytic function $f(z)$ can be written
as
\[ f(z) = u(x,y) + i v(x,y),\]
where $z = x + i y$ and $u(x,y)$, $v(x,y)$ are each real-valued
\emph{harmonic} functions related to each other in a very strong
way: the Cauchy-Riemann equations
\begin{equation}\label{eCRconds}
\frac{\partial u}{\partial x} = \frac{\partial v}{\partial y} \qquad
\frac{\partial v}{\partial x} = - \frac{\partial u}{\partial y} .
\end{equation}
In summary, the deceptively simple hypothesis that
\[ f'(z) \quad \text{exists} \]
forces a great deal of structure on $f(z)$; moreover, this structure
mirrors the structure of the harmonic $u(x,y)$ and $v(x,y)$,
functions of \emph{two} real variables.\footnote{Quoted from page 2
of reference \cite{Fl72}. Note that in the quote $i = \sqrt{-1}$
whereas in this note we take $j = \sqrt{-1}$ following standard
electrical engineering practice.}
\end{quote}

In particular the following conditions are equivalent statements
about a complex function $f(z)$ on an open set containing $z$ in the
complex plane \cite{Fl72}:
\begin{itemize}
\item The derivative $f'(z)$ exists and is continuous.

\item The function $f(z)$ is \emph{holomorphic} (i.e, \emph{complex}-analytic in $z$).\footnote{A function is
\emph{analytic} on some domain if it can be expanded in a convergent
power series on that domain. Although this condition implies that
the function has derivatives of all orders, analyticity is a
stronger condition than infinite differentiability as there exist
functions which have derivatives of all orders but which cannot be
expressed as a power series. For a complex-valued function of a
complex variable, the term (complex) analytic has been replaced in modern
mathematics by the entirely synonymous term \emph{holomorphic}. Thus
\emph{real-valued} power-series-representable functions of a
\emph{real-variable} are analytic (real-analytic), while
\emph{complex-valued} power-series-representable functions of a
\emph{complex-variable} are \emph{holomorphic} (complex-analytic).  We can now appreciate the merit of
distinguishing between holomorphic and (real) analytic functions---a
function can be nonholomorphic (i.e. non-complex-analytic) in the
\emph{complex variable} $z = x + j \, y$ yet still be (real) analytic in
the \emph{real variables} $x$ and $y$.}

\item The function $f(z)$ satisfies the \emph{Cauchy-Riemann
conditions} (\ref{eCRconds}).

\item All derivatives of the function $f(z)$ exist and $f(z)$ has a
convergent power series.
\end{itemize}
Furthermore, it is a simple consequence of the Cauchy-Riemann
conditions that \[ f(z) = u(x,y) + j \, v(x,y) \] is holomorphic
only if the functions $u(x,y)$ and $v(x,y)$ both satisfy Laplace's
equation
\[  \frac{\partial^2 u(x,y)}{\partial x^2} + \frac{\partial^2 u(x,y)}{\partial
y^2}= 0 \quad \text{and} \quad  \frac{\partial^2 v(x,y)}{\partial
x^2} + \frac{\partial^2 v(x,y)}{\partial y^2} =0 .
\]
Such functions are known as \emph{harmonic} functions. Thus if
either $u(x,y)$ or $v(x,y)$ fail to be harmonic, the function $f(z)$
is not differentiable.\footnote{Because a harmonic function on
$\mathbb{R}^2$ satisfies the partial differential equation known as
Laplace's equation, by existence and uniqueness of the solution to
this partial differential equation its value is completely
determined at a point in the interior of \emph{any} simply connected
region which contains that point once the values on the boundary
(boundary conditions) of that region are specified. This is the
reason that contour integration of a complex-analytic (holomorphic) function
works and that we have the freedom to select that contour to make
the integration as easy as possible. On the other hand, there is, in
general, no equivalent to contour integration for an arbitrary
function on $\mathbb{R}^2$. See the excellent discussion in Flanigan
\cite{Fl72}.}

\medskip
Although many important complex functions are holomorphic, including
the functions $z^n$, $\text{e}^z$, $\ln (z)$, $\sin(z)$, and
$\cos(z)$, and hence differentiable in the standard complex
variables sense, there are commonly encountered useful functions
which are not:
\begin{itemize}
\item The function $f(z) = \bar{z}$, where `$\bar{z}$' denotes
complex conjugation, fails to satisfy the Cauchy-Riemann conditions.

\item The functions $f(z) = \text{Re}(z)= \frac{z + \bar{z}}{2} =
x $ and $g(z) = \text{Im}(z) = \frac{z - \bar{z}}{2j} = y$ fail the
Cauchy-Riemann conditions.

\item The function $f(z) = \left| z \right|^2 = \bar{z} z = x^2 +
y^2$ is not harmonic.

\item Any \emph{nonconstant purely real-valued} function $f(z)$
(for which it must be the case that  $v(z,y) \equiv 0$) fails the
Cauchy-Riemann condition. In particular the real function $f(z) =
\left| z \right| =  \sqrt{\bar{z} z} = \sqrt{x^2 + y^2}$ is not
differentiable.\footnote{Thus we have the classic result that the
only holomorphic real-valued functions are the constant real-valued
functions.}
\end{itemize}

Note in particular, the implication of the above for the problem of
minimizing the real-valued squared-error loss functional
\begin{equation}\label{eq:e0}
  \ell(a) = \text{E} \left\{ \left| \eta_k - \bar{a} \xi_k \right|^2 \right\}
  = \text{E} \left\{ \overline{(\eta_k - \bar{a} \xi_k )}(\eta_k -\bar{a} \xi_k) \right\}
  \triangleq \text{E} \left\{ \bar{e}_k e_k \right\}
\end{equation}
for finite second-order moments stationary scalar complex random
variables $\xi_k$ and $\eta_k$, and unknown complex constant $a =
a_x + j a_y$. Using the theory of optimization in Hilbert spaces,
the minimization can be done by invoking the \emph{projection
theorem} (which is equivalent to the \emph{orthogonality
principle}) \cite{Lu69}. Alternatively, the minimization can be
performed by completing the square. Either procedure will result in
the Wiener-Hopf equations, which can then be solved for the optimal
complex coefficient variable $a$.

\medskip
However, if a gradient procedure for determining the optimum is
desired, we are immediately stymied by the fact that the
\emph{purely real} nonconstant function $\ell(a)$ is \emph{not}
complex-analytic (holomorphic) and therefore its derivative  with respect to $a$
\emph{does not exist in the conventional sense} of a complex
derivative \cite{Ch96}-\cite{Fl72}, which applies only to
holomorphic functions of $a$. A way to break this impasse will be
discussed in the next section. Meanwhile note that \emph{all} of the
real-valued nonholomorphic functions shown above can be viewed as
functions of both $z$ \emph{and} its complex conjugate $\bar{z}$, as
this fact will be of significance in the following discussion.

\section{\boldmath Extensions of the Complex Derivative -- The $\mathbb{CR}$-Calculus}
\label{sec:cpd}

In this section we continue to focus on functions of a \emph{single}
complex variable $z$. The primary references for the material
developed here are Nehari \cite{ZN61}, Remmert \cite{RR91}, and
Brandwood \cite{Br83}.

\subsection{A Possible Extension of the Complex Derivative.} As we have seen, in order for
the complex derivative of a function of $z = x + j\, y$,
\[ f(z) =
u(x,y) + j \, v(x,y),
\] to exist in the standard holomorphic sense, the real partial
derivatives of $u(x,y)$ and $v(x,y)$ must not only exist, they must
\emph{also} satisfy the Cauchy-Riemann conditions (\ref{eCRconds}). As
noted by Flanigan \cite{Fl72}: ``This is much stronger than the mere
\emph{existence} of the partial derivatives.''  However, the ``mere
existence'' of the (real) partial derivatives \emph{is} necessary
and sufficient for a stationary point of a (necessarily
nonholomorphic) non-constant \emph{real-valued} functional $f(z)$ to
exist when $f(z)$ is viewed \emph{as a differentiable function of
the real and imaginary parts of $z$,} i.e., as a function over
$\mathbb{R}^2$,
\begin{equation}\label{eq:e1}
f(z) = f(x,y): \mathbb{R}^2 \rightarrow \mathbb{R} \, .
\end{equation}
Thus the trick is to exploit the real $\mathbb{R}^2$ vector space
structure which underlies $\mathbb{C}$ when performing
gradient-based optimization. In essence, the remainder of this note
is concerned with a thorough discussion of this ``trick.''

\medskip
Towards this end, it is convenient to define a generalization or
extension of the standard partial derivative to nonholomorphic
functions of $z= x + j \, y$ that are nonetheless differentiable
with respect to $x$ and $y$ and which incorporates the real gradient
information directly within the complex variables framework. After
Remmert \cite{RR91}, we will call this the \emph{real-derivative},
or \emph{{\footnotesize $\mathbb{R}$}-derivative,} of a possibly
nonholomorphic function in order to avoid confusion with the
standard \emph{complex-derivative,} or \emph{{\footnotesize
$\mathbb{C}$}-derivative}, of a holomorphic function which was
presented and discussed in the previous section. Furthermore,
\emph{we would like the real-derivative to reduce to the standard
complex derivative when applied to holomorphic functions}.

\medskip
Note that if one rewrites the real-valued loss function
(\ref{eq:e0}) in terms of purely real quantities, one obtains
(temporarily suppressing the time dependence, $k$)
\begin{equation}\label{eq:e3}
 \ell(a) =  \ell(a_x,a_y) = \text{E} \left\{ e_x^2 + e_y^2 \right\}
  = \text{E} \left\{ \left( \eta_x - a_x \xi_x - a_y \xi_y \right)^2
  + \left( \eta_y + a_y \xi_x - a_x \xi_y \right)^2 \right\}
   \, .
\end{equation}
 From this we can easily determine that
\[
  \frac{\partial \ell(a_x,a_y)}{\partial a_x} = -2 \, \text{E} \left\{ e_x \xi_x
  +
  e_y \xi_y \right\}
  \, ,
\]
and
\[
  \frac{\partial \ell(a_x,a_y)}{\partial a_y} =   -2 \, \text{E} \left\{  e_x \xi_y
  -
  e_y \xi_x \right\} \, .
\]
Together these can be written as
\begin{equation}\label{eq:e2}
  \left( \frac{\partial}{\partial a_x} + j \frac{\partial}{\partial
  a_y}\right) \ell(a) =
  \frac{\partial \ell(a_x,a_y)}{\partial a_x} + j \frac{\partial \ell(a_x,a_y)}{\partial a_y}
   = -2 \, \text{E} \left\{ \xi_k \bar{e}_k \right\}
\end{equation}
which looks very similar to the standard result for the real case.

\medskip
Indeed, equation (\ref{eq:e2}) is the definition of the generalized
complex partial derivative often given in engineering textbooks,
including references \cite{GS01}-\cite{MM80}. However, this is
\emph{not} the definition used in this note, which instead follows
the formulation presented in \cite{Wi27}-\cite{BS92}. We do not use
the definition (\ref{eq:e2}) because it \emph{does not} reduce to
the standard $\mathbb{C}$-derivative for the case when a function
$f(a)$ \emph{is} a holomorphic function of the complex variable $a$.
For example, take the simplest case of $f(a) = a$, for which the
standard derivative yields $\frac{d}{da} f(a) = 1$. In this case,
the definition (\ref{eq:e2}) applied to $f(a)$ unfortunately results
in the value $0$. Thus we will \emph{not} view the definition
(\ref{eq:e2}) as an admissible \emph{generalization} of the standard
complex partial derivative, although it does allow the determination
of the stationary points of $\ell(a)$.\footnote{In fact, it is a
scaled version of the conjugate {\footnotesize
$\mathbb{R}$}-derivative discussed in the next subsection.}

\subsection{The $\mathbb{R}$-Derivative and Conjugate
$\mathbb{R}$-Derivative.} There are a variety of ways to develop the
formalism discussed below (see \cite{ZN61}-\cite{Br83}). Here, we
roughly follow the development given in Remmert \cite{RR91} with
additional material drawn from Brandwood \cite{Br83} and Nehari
\cite{ZN61}.

\medskip
Note that the \emph{nonholomorphic} (nonanalytic in the complex variable
$z$) functions given as examples
in the previous section can all be written in the form
$f(z,\bar{z})$, where they \emph{are} holomorphic in $z = x + j \,
y$ \emph{for fixed $\bar{z}$} and holomorphic in $\bar{z} = x - j \,
y$ \emph{for fixed $z$.}\footnote{That is, if we make the
substitution $w = \bar{z}$, they are analytic in $w$ for fixed $z$,
and analytic in $z$ for fixed $w$. This simple insight underlies the
development given in Brandwood \cite{Br83} and Remmert \cite{RR91}.}
It can be shown that \emph{this fact is true in general} for any
complex- or real-valued function
\begin{equation}\label{eq:zzzbarxy}
 f(z) = f(z,\bar{z}) = f(x,y) = u(x,y) + j\, v(x,y)
\end{equation}
of a complex variable for which the real-valued functions $u$ and
$v$ are differentiable as functions of the real variables $x$ and
$y$. This fact underlies the development of the so-called
\emph{Wirtinger calculus} \cite{RR91} (or, as we shall refer to it
later, the \emph{{\footnotesize $\mathbb{CR}$}-calculus.}) In
essence, the so-called \emph{conjugate coordinates,}

\begin{equation}\label{eq:firstconj}
\boxed{\ \ \text{\sf Conjugate Coordinates:} \quad c \triangleq
(z,\bar{z})^T \in \mathbb{C} \times \mathbb{C} \, , \quad z = x +
j\, y  \quad \text{and} \quad \bar{z} = x - j \, y \ \ }
\end{equation}

\noindent can serve as a formal substitute for the real $r =
(x,y)^T$ representation of the point $z = x + j \, y  \in
\mathbb{C}$ \cite{RR91}.\footnote{\textbf{Warning!} The interchangeable use
of the various notational forms of $f$ implicit in the statement
$f(z) = f(z,\bar{z})$ can lead to confusion. To minimize this
possibility we define the term ``{$f(z)$ ($z$-only)}'' to mean that
$f(z)$ is independent of $\bar{z}$ (and hence is holomorphic) and
the term ``$f(\bar{z})$ ($\bar z$ only)'' to mean that $f(z)$ is a
function of $\bar{z}$ only. Otherwise there are no restrictions on
$f(z) = f(z,\bar{z})$.} According to Remmert \cite{RR91}, the
calculus of complex variables utilizing this perspective was
initiated by Henri Poincar\'{e} (over 100 years ago!) and further
developed by Wilhelm Wirtinger in the 1920's \cite{Wi27}. Although
this methodology has been fruitfully exploited by the
German-speaking engineering community (see, e.g., references
\cite{Fi02} or \cite{Fr97}), it has not generally been appreciated
by the English speaking engineering community until relatively
recently.\footnote{An important exception is Brandwood \cite{Br83}
and the work that it has recently influenced such as
\cite{VT02,Ka93,Ha96}. However, these latter references do not seem
to fully appreciate the clarity and ease of computation that the
Wirtinger calculus ({\footnotesize $\mathbb{CR}$}-calculus) can
provide to the problem of differentiating nonholomorphic function
and optimizing real-valued functions of complex variables. Perhaps
this is do to the fact that \cite{Br83} did \emph{not} reference the
Wirtinger calculus as such, nor cite the rich body of work which had
already existed in the mathematics community
(\cite{ZN61,LH90,RR91}).}

\medskip For a general complex- or real-valued function $f(c) = f(z,\bar{z})$ consider the
\emph{pair} of partial derivatives of $f(c)$
\emph{formally}\footnote{These statements are \emph{formal} because
one cannot truly vary $z = x + j \, y$ while keeping $\bar{z} = x -j
\, y$ constant, and vice versa.} defined by

{\footnotesize
\begin{equation} \label{eq:e4} \boxed{\ \ \left.
\underline{\text{\sf $\mathbb{R}$-Derivative of $f(c)$}} \triangleq
\frac{\partial f(z,\bar{z})}{\partial z} \right|_{\bar{z} = \,
\text{const.}}
 \ \text{and} \ \ \
\underline{\text{\sf Conjugate $\mathbb{R}$-Derivative of $f(c)$}}
\triangleq \left. \frac{\partial f(z,\bar{z})}{\partial \bar{z}}
\right|_{z = \, \text{const.}} \ \ }
\end{equation}}

\noindent where the formal partial derivatives are taken to be
standard complex partial derivatives ($\mathbb{C}$-derivatives)
taken with respect to $z$ in the first case and with respect to
$\bar{z}$ in the second.\footnote{\label{foot:b}A careful and
rigorous analysis of these formal partial derivatives can be found
in Remmert \cite{RR91}. In \cite{RR91}, a differentiable complex
function $f$ is called \emph{$\mathbb{C}$-differentiable} while if
$f$ is differentiable as a mapping from $\mathbb{R}^2 \rightarrow
\mathbb{R}^2$, it is said to be \emph{real-differentiable}
($\mathbb{R}$-differentiable) (See Footnote \ref{foot:a}). It is
shown in \cite{RR91} that the partial derivatives (\ref{eq:e4})
exist if and only if $f$ is $\mathbb{R}$-differentiable. As
discussed further below, throughout this note we assume that all
functions are globally \emph{real-analytic} ($\mathbb{R}$-analytic),
which is a sufficient condition for a function to be globally
$\mathbb{R}$-differentiable.} For example, with $f(z,\bar{z}) = z^2
\bar{z}$ we have
\[ \frac{\partial f}{\partial z} = 2 z \bar{z} \quad \text{and}
\quad \frac{\partial f}{\partial \bar{z}} = z^2 \, . \]
 As denoted in (\ref{eq:e4}), we
call the first expression the \emph{{\footnotesize
$\mathbb{R}$}-derivative} (the \emph{real-derivative}) and the
second expression the \emph{conjugate {\footnotesize
$\mathbb{R}$}-derivative} (or {\footnotesize $\overline{\,
\mathbb{R}}$}-\emph{derivative}).

\medskip
It is proved in \cite{ZN61,Br83,RR91} that the {\footnotesize
$\mathbb{R}$}-derivative and {\footnotesize $\overline{\,
\mathbb{R}}$}-derivative formally defined by (\ref{eq:e4}) can be
equivalently written as\footnote{Recall the representation $f =
f(x,y) = u(x,y) + j\, v(x,y)$. Note that the relationships
(\ref{eq:equivdef}) make it clear why the partial derivatives
(\ref{eq:e4}) exist if and only if $f$ is
$\mathbb{R}$-differentiable. (See footnotes \ref{foot:a} and
\ref{foot:b}).}
\begin{equation}\label{eq:equivdef}
\frac{\partial f}{\partial z} = \frac{1}{2} \, \left( \frac{\partial
f}{\partial x} - j \frac{\partial f}{\partial
  y}\right) \quad \text{and} \quad \frac{\partial f}{\partial \bar{z}} =  \frac{1}{2} \,
  \left( \frac{\partial f}{\partial x} + j \frac{\partial f}{\partial
  y}\right) \,
\end{equation}
where the partial derivatives with respect to $x$ and $y$ are
\emph{true} (i.e., non-formal) partial derivatives of the function
$f(z) = f(x,y)$, which is always assumed in this note to be
differentiable with respect to $x$ and $y$ (i.e., to be
{\footnotesize $\mathbb{R}$}-differentiable). Thus it is the
\emph{right-hand-sides} of the expressions given in
(\ref{eq:equivdef}) which make rigorous the formal definitions of
(\ref{eq:e4}).

\medskip
Note that from equation (\ref{eq:equivdef}) that we immediately have
the properties
\begin{equation}\label{eq:derivfacts1}
\frac{\partial z}{\partial z} = \frac{\partial \bar{z}}{\partial
\bar{z}} = 1 \quad \text{and} \quad \frac{\partial z}{\partial
\bar{z}} = \frac{\partial \bar{z}}{\partial z} = 0 \, .
\end{equation}

\medskip
\noindent{\bf Comments:}
\begin{enumerate}

\item  The condition $\frac{\partial f}{\partial \bar{z}} = 0$ is
true for an {\footnotesize $\mathbb{R}$}-differentiable function $f$
if and only the Cauchy-Riemann conditions are satisfied (see
\cite{ZN61,Br83,RR91}). \emph{Thus a function $f$ is holomorphic
(complex-analytic in $z$) if and only if it does not depend on the complex
conjugated variable $\bar{z}$. I.e., if and only if $f(z) = f(z)$ ($z$
only)}.\footnote{This obviously provides a simple and powerful
characterization of holomorphic and nonholomorphic functions and
shows the elegance of the Wirtinger calculus formulation based on
the use of conjugate coordinates $(z,\bar{z})$. Note that the two
Cauchy-Riemann conditions are replaced by the single condition
$\frac{\partial f}{\partial \bar{z}} = 0 $. The reader should
reexamine the nonholomorphic (nonanalytic in $z$) functions
discussed in the previous section in the light of this condition.}

\item The $\mathbb{R}$-derivative, $\frac{\partial f}{\partial z}$,
of an {\footnotesize $\mathbb{R}$}-differentiable function $f$ is
equal to the standard $\mathbb{C}$-derivative, $f'(z)$, when
$f(z,\bar{z})$ is independent of $\bar{z}$, i.e., when $f(z) = f(z)$
($z$ only).

\item An {\footnotesize $\mathbb{R}$}-differentiable
 function $f$ is holomorphic in $\bar{z}$ (complex-analytic in $\bar{z}$)
if and only if it does not depend on the variable $z$, $f(z,\bar{z})
= f(\bar{z})$ ($\bar z$ only), which is true if and only if
$\frac{\partial f}{\partial z} = 0$.
\end{enumerate}
To summarize, an {\footnotesize $\mathbb{R}$}-differentiable
function $f$ is holomorphic (complex-analytic in $z$) if and only if $f(z) =
f(z)$ ($z$ only), which is true if and only if $\frac{\partial
f}{\partial \bar{z}} = 0$, in which case the {\footnotesize
$\mathbb{R}$}-derivative coincides with the standard {\footnotesize
$\mathbb{C}$}-derivative, $\frac{\partial f}{\partial z} = f'(z)$.
We call the \emph{single condition} $\frac{\partial f}{\partial
\bar{z}} = 0$ the \emph{Cauchy-Riemann condition} for $f$ to be
holomorphic:
\begin{equation}
\boxed{\ \text{\sf Cauchy Riemann Condition:} \quad \frac{\partial
f}{\partial \bar{z}} = 0 \ }
\end{equation}

\paragraph{Real-Analytic Complex Functions.}
Throughout the discussion given above we have been making the
assumption that a complex function $f$ is real differentiable
({\footnotesize $\mathbb{R}$}-differentiable). We henceforth make
the stronger assumption that complex functions over $\mathbb{C}$ are
globally \emph{real-analytic} ({\footnotesize
$\mathbb{R}$}-analytic) over $\mathbb{R}^2$. As discussed above, and
rigorously proved in Remmert \cite {RR91}, {\footnotesize
$\mathbb{R}$}-analytic functions are {\footnotesize
$\mathbb{R}$}-differentiable and {\footnotesize $\overline{\,
\mathbb{R}}$}-differentiable.

\medskip
A function $f(z)$ has a power series expansion in the complex
variable $z$,
\[ f(z) = f(z_0) + f'(z_0) (z -z_0) + \frac{1}{2}f''(z_0) (z - z_0)^2 + \cdots
+ \frac{1}{n!} f^{(n)}(z_0) (z-z_0)^n + \cdots \] where the complex
coefficient $f^{(n)}(z_0)$ denotes an $n$-times {\footnotesize
$\mathbb{C}$}-derivative of $f(z)$ evaluated at the point $z_0$, if
and only if it is holomorphic in an open neighborhood of $z_0$. If the function $f(z)$ is not
holomorphic over $\mathbb{C}$, so that the above expansion does not
exist, but is nonetheless still {\footnotesize
$\mathbb{R}$}-analytic as a mapping from $\mathbb{R}^2$ to
$\mathbb{R}^2$, then the real and imaginary parts of $f(z) = u(x,y)
+ j \, v(x,y)$, $z = x+ j\, y$, can be expanded in terms of the real
variables $r = (x,y)^T$,
\begin{eqnarray*}
u(r) & = & u(r_0) + \frac{\partial u(r_0)}{\partial r} (r - r_0)+
(r-r_0)^T \frac{\partial}{\partial r} \left(\frac{\partial
u(r_0)}{\partial
r}\right)^T (r- r_0) + \cdots \\
v(r) & = & v(r_0) + \frac{\partial v(r_0)}{\partial r} (r - r_0)+
(r-r_0)^T \frac{\partial}{\partial r} \left(\frac{\partial
v(r_0)}{\partial r}\right)^T (r- r_0) + \cdots
\end{eqnarray*}
Note that if the {\footnotesize $\mathbb{R}$}-analytic function is
\emph{purely real,} then $f(z) = u(x,y)$ and we have
\[ f(r) =  f(r_0) + \frac{\partial f(r_0)}{\partial r} (r - r_0) +
(r-r_0)^T \frac{\partial}{\partial r} \left(\frac{\partial
f(r_0)}{\partial r}\right)^T (r- r_0) + \cdots  \]

\newpage
\paragraph{Properties of the {$\mathbb{R}$}- and {$\overline{\,
\mathbb{R}}$}-Derivatives.} The {\footnotesize
$\mathbb{R}$-derivative} and {\footnotesize $\overline{\,
\mathbb{R}}$}-derivative are both \emph{linear operators} which obey
the \emph{product rule} of differentiation. The following important
and useful properties also hold (see references
\cite{ZN61,RR91}).\footnote{In the following for $z = x + j\, y$ we
define $dz = dx + j\, dy$ and $d\bar{z} = d x - j\, dy$, while $h(g)
= h \circ g$ denotes the composition of the two function $h$ and
$g$.}

\medskip
\bigskip \noindent{\large\bf Complex Derivative Identities:}

\begin{eqnarray}
\frac{\partial \bar{f}}{\partial \bar{z}} & = &
\overline{\left(\frac{\partial f}{\partial z}\right)}
\label{eq:conjusa}
\\
\frac{\partial \bar{f}}{\partial z} & = &
\overline{\left(\frac{\partial f}{\partial \bar{z}}\right)} \label{eq:conjusb}\\
df & = & \frac{\partial f}{\partial z}\,  dz + \frac{\partial
f}{\partial \bar{z}}\, d \bar{z} \qquad \quad \ \, \text{\sf Differential Rule} \label{eq:drule} \\
 \frac{\partial
h(g)}{\partial z} & = & \frac{\partial h}{\partial g}\,
\frac{\partial g}{\partial z} + \frac{\partial h}{\partial
\bar{g}}\, \frac{\partial \bar{g}}{\partial z} \qquad \quad
\text{\sf Chain Rule} \label{eq:ch1}
\\
\frac{\partial h(g)}{\partial \bar{z}} & = & \frac{\partial
h}{\partial g}\,  \frac{\partial g}{\partial \bar{z}} +
\frac{\partial h}{\partial \bar{g}}\, \frac{\partial
\bar{g}}{\partial \bar{z}} \qquad \quad \text{\sf Chain Rule}
\label{eq:ch2}
\end{eqnarray}
As a simple consequence of the above, note that if $f(z)$ is
real-valued then $\bar{f}(z) = f(z)$ so that we have the additional
very important identity that
\begin{equation} \label{eq:reallossresult}
f(z) \in \mathbb{R}\  \Rightarrow \ \overline{\left(\frac{\partial
f}{\partial z}\right)} = \frac{\partial f }{\partial \bar{z}}
\end{equation}

\bigskip As a simple first application of the above, note that the
{\footnotesize $\overline{\, \mathbb{R}}$}-derivative of $\ell(a)$
can be easily computed from the definition (\ref{eq:e0}) and the
above properties to be
\begin{equation}\label{eq:e5}
\frac{\partial \ell(a)}{\partial \bar{a}} = \text{E} \left\{
\bar{e}_k e_k \right\} = \text{E} \left\{ \frac{\partial
\bar{e}_k}{\partial \bar{a}}\,  e_k + \bar{e}_k \, \frac{\partial
e_k}{\partial \bar{a}} \right\}= \text{E} \left\{ 0 \cdot e_k -
\bar{e}_k\,  \xi_k \right\} = - \, \text{E} \left\{ \xi_k \,
\bar{e}_k \right\} \, .
\end{equation}
which is the same result obtained from the ``brute force'' method
based on deriving expanding the loss function in terms of the real
and imaginary parts of $a$, followed by computing (\ref{eq:e2}) and
then using the result (\ref{eq:equivdef}). Similarly, it can be
easily shown that the $\mathbb{R}$-derivative of $\ell(a)$ is given
by
\begin{equation}\label{eq:e6}
 \frac{\partial
\ell(a)}{\partial a} = - \, \text{E} \left\{ \bar{\xi}_k e_k
\right\} \, .
\end{equation}
Note that the results (\ref{eq:e5}) and (\ref{eq:e6}) are the
complex conjugates of each other, which is consistent with the
identity (\ref{eq:reallossresult}).

\medskip
We view the \emph{pair} of formal partial derivatives for a possibly
nonholomorphic function defined by (\ref{eq:e4}) as the natural
generalization of the \emph{single} complex derivative
($\mathbb{C}$-derivative) of a holomorphic function. The fact that
there are \emph{two} derivatives under general consideration does
not need to be developed in elementary standard complex analysis
courses where it is usually assumed that $f$ is always holomorphic
(complex-analytic in $z$). In the case when $f$ is holomorphic then $f$ is
independent of $\bar{z}$ and the conjugate partial derivative is
zero, while the extended derivative reduces to the standard complex
derivative.

\medskip
\paragraph{First-Order Optimality Conditions.} As mentioned in the introduction, we are
often interested in optimizing a scalar function with respect to the
real and imaginary parts $r = (x,y)^T$ of a complex number $z = x +
j \, y$. It is a standard result from elementary calculus that a
first-order necessary condition for a point $r_0= (x_0,y_0)^T$ to be
an optimum is that this point be a stationary point of the loss
function. Assuming differentiability, stationarity is equivalent to
the condition that the partial derivatives  of the loss function
with respect the parameters $r=(x,y)^T$ vanish at the point $r=
(x_0,y_0)^T$. The following fact is an easy consequence of the
definitions (\ref{eq:e4}) and is discussed in \cite{Br83}:
\begin{itemize}
\item A necessary and sufficient condition for a real-valued
function, $f(z) = f(x,y)$, $z = x + j \, y$, to have a stationary
point with respect to the real parameters $r = (x,y)^T \in
\mathbb{R}^2$ is that its {\footnotesize $\overline{\,
\mathbb{R}}$}-derivative vanishes. \quad Equivalently, a necessary
and sufficient condition for $f(z) = f(x,y)$ to have a stationary
point with respect to $r = (x,y)^T \in \mathbb{R}^2$ is that its
$\mathbb{R}$-derivative vanishes.
\end{itemize}

For example, setting either of the derivatives (\ref{eq:e5}) or
(\ref{eq:e6}) to zero results in the so-called Wiener-Hopf equations
for the optimal MMSE estimate of $a$. This result can be readily
extended to the multivariate case, as will be discussed later in
this note.

\paragraph{The Univariate $\mathbb{CR}$-Calculus.} As noted in \cite{RR91},
the approach we have been describing is known as the Wirtinger
calculus in the German speaking countries, after the pioneering work
of Wilhelm Wirtinger in the 1920's \cite{Wi27}. Because this
approach is based on being able to apply the calculus of \emph{real
variables} to make statements about functions of \emph{complex
variables}, in this note we use the term ``{\footnotesize
$\mathbb{CR}$}-calculus'' interchangeable with ``Wirtinger
calculus.''

\medskip
Despite the important insights and ease of computation that it can
provide, it is the case that the use of conjugate coordinates $z$
and $\bar z$ (which underlies the {\footnotesize
$\mathbb{CR}$}-calculus) is \emph{not} needed when developing the
classical univariate theory of holomorphic (complex-analytic in $z$)
functions.\footnote{``The differential calculus of these operations
... [is] ... largely irrelevant for classical function theory ...''
--- R.~Remmert \cite{RR91}, page 66.} It is only in the multivariate and/or nonholonomic
case that the tools of the {\footnotesize $\mathbb{CR}$}-calculus
begin to be indispensable. Therefore it is not developed in the
standard courses taught to undergraduate engineering and science
students in this country \cite{Ch96}-\cite{Fl72} which have changed
little in mode of presentation from the earliest
textbooks.\footnote{For instance, the widely used textbook by
Churchill \cite{Ch96} adheres closely to the format and topics of
its first edition which was published in 1948. The latest edition
(the 7th at the time of this writing) does appear to have one brief
homework problem on differentiating nonholomorphic functions.}

\medskip
Ironically, the elementary textbook by Nehari \cite{ZN61} was an
attempt made in 1961 (almost 50 years ago!) to integrate at least some
aspects of the {\footnotesize $\mathbb{CR}$}-calculus into the
elementary treatment of functions of a single complex
variable.\footnote{This is still an excellent textbook that is
highly recommended for an accessible introduction to the use of
derivatives based on the conjugate coordinates $z$ and $\bar{z}$.}
However, because the vast majority of textbooks treat the univariate
case, as long as the mathematics community, and most of the
engineering community, was able to avoid dealing with nonholomorphic
functions, there was no real need to bring the ideas of the
{\footnotesize $\mathbb{CR}$}-calculus into the mainstream
univariate textbooks.

\medskip
Fortunately, an excellent, sophisticated and extensive  introduction
to univariate complex variables theory and the {\footnotesize
$\mathbb{CR}$}-calculus  is available in the textbook by Remmert
\cite {RR91}, which is a translation from the 1989 German edition.
This book also details the historical development of complex
analysis. The highly recommended Remmert and Nehari texts have been
used as primary references for this note (in addition to the papers
by Brandwood \cite{Br83} and, most importantly for the second-order analysis given below, van den Bos \cite{VDB94a}).

\paragraph{\boldmath The Multivariate {$\mathbb{CR}$}-Calculus.}
Although one can forgo the tools of the {\footnotesize
$\mathbb{CR}$}-calculus in the case of univariate holomorphic
functions, this is not the situation in the multivariate holomorphic
case where mathematicians have long utilized these tools
\cite{GH84}-\cite{BS92}.\footnote{``[The {\footnotesize
$\mathbb{CR}$}-calculus] is quite indispensable  in the function
theory of several variables.''
--- R.~Remmert \cite{RR91}, page 67.} Unfortunately, multivariate
complex analysis is highly specialized and technically abstruse, and
therefore virtually all of the standard textbooks are accessible
only to the specialist or to the aspiring specialist. It is commonly
assumed in these textbooks that the reader has great facility with
differential geometry, topology, calculus on manifolds, and
differential forms, in addition to a good grasp of advanced
univariate complex variables theory.  Moreover, because the focus of
the theory of multivariate complex functions is primarily on
\emph{holomorphic} functions, whereas our concern is the essentially
ignored (in this literature) case of nonholomorphic real-valued
functionals, it appears to be true that only a very small part of
the material presented in these references is directly useful for our purposes (and primarily for
creating a rigorous and self-consistent multivariate {\footnotesize
$\mathbb{CR}$}-calculus framework based on the results given in the
papers by Brandwood \cite{Br83} and van den Bos \cite{VDB94a}).

\medskip
The clear presentation by Brandwood \cite{Br83} provides a highly
accessible aspect of the first-order multivariate {\footnotesize
$\mathbb{CR}$}-calculus as applied to the problem of finding stationary values for
real-valued functionals of complex variables.\footnote{Although, as
mentioned in an earlier footnote, Brandwood for some reason did not
cite or mention any prior work relating to the use of conjugate
coordinates or the Wirtinger calculus.} As this is the primary
interest of many engineers, this pithy paper is a very useful
presentation of just those very few theoretical and practical issues
which are needed to get a clear grasp of the problem. Unfortunately,
even twenty years after its publication, this paper still is not as
widely known as it should be. However, the recent utilization of the
Brandwood results in \cite{VT02,Fi02,Ka93,Ha96} seems to indicate a
standardization of the Brandwood presentation of the complex
gradient into the mainstream textbooks. The results given in the
Brandwood paper \cite{Br83} are particulary useful when coupled with
with the significant extension of Brandwood's results to the problem
of computing complex Hessians which has been provided by van den
Bos's paper \cite{VDB94a}.

\medskip
At this still relatively early stage in the development of a widely
accepted framework for dealing with real-valued (nonholomorphic)
functions of several complex variables, presumably even the
increasingly widely used formalism of Brandwood \cite{Br83} and van
den Bos \cite{VDB94a} potentially has some room for improvement
and/or clarification (though this is admittedly a matter of taste).
In this spirit, and mindful of the increasing acceptance of the
approach in \cite{Br83} and \cite{VDB94a}, in the remainder of this
note we develop a multivariate {\footnotesize
$\mathbb{CR}$}-calculus framework that is only slightly different
than that of \cite{Br83} and \cite{VDB94a}, incorporating insights
available from the literature on the calculus of multivariate
complex functions and complex differential manifolds
\cite{GH84}-\cite{BS92}.\footnote{Realistically, one must admit that
many, and likely most, practicing engineers will be unlikely to make the move
from the perspective and tools provided by \cite{Br83} and
\cite{VDB94a} (which already enable the engineer to solve most
problems of practical interest) to that developed in this note,
primarily because of the requirement of some familiarity of (or
willingness to learn) concepts of differential geometry  at the level of the earlier
chapters of \cite{Sh80} and \cite{TF01}).}

\section{\boldmath Multivariate {$\mathbb{CR}$}-Calculus}

The remaining sections of this note will provide an expanded
discussion and generalized presentation of the \emph{multivariate
{\footnotesize $\mathbb{CR}$}-calculus} as presented in Brandwood
\cite{Br83} and van den Bos \cite{VDB94a}.
The discussion given below also utilizes insights gained from references
\cite{GH84,LH90,SK01,BS92,Sh80,TF01}.

\subsection{\boldmath The Space $\mathcal{Z} = \mathbb{C}^n$.}

We define the $n$-dimensional column vector $\mathbf{z}$ by
\[ \mathbf{z} = \begin{pmatrix} z_1 & \cdots & z_n \end{pmatrix}^T
\in \mathcal{Z} = \mathbb{C}^n
\]
where $z_i = x_i + j \, y_i$, $i = 1, \cdots , n$, or, equivalently,
\[ \mathbf{z} = \mathbf{x} + j \, \mathbf{y} \]
with $\mathbf{x} = ( x_1 \cdots x_n)^T$ and $\mathbf{y} = ( y_1
\cdots y_n)^T$. The space $\mathcal{Z} = \mathbb{C}^n$ is a vector
space over the field of complex numbers with the standard
component-wise definitions of vector addition and scalar
multiplication. Noting the one-to-one correspondence
\[ \mathbf{z} \in \mathbb{C}^n \Leftrightarrow \mathbf{r} =
\begin{pmatrix} \mathbf{x} \\ \mathbf{y} \end{pmatrix} \in
\mathcal{R} \triangleq \mathbb{R}^{2n} = \mathbb{R}^n \times
\mathbb{R}^n
\] it is evident that there exists a natural isomorphism between
$\mathcal{Z} = \mathbb{C}^n$ and $\mathcal{R} = \mathbb{R}^{2n}$.

\medskip
The conjugate coordinates of $\mathbf{z} \in \mathbb{C}^n$ are
defined by
\[ \mathbf{\bar{z}} = \begin{pmatrix} \bar{z}_1 & \cdots & \bar{z}_n \end{pmatrix}^T
\in \mathcal{Z} = \mathbb{C}^n \,
\]
We denote the pair of conjugate coordinate vectors
$(\mathbf{z},\mathbf{\bar{z}})$ by
\[ \mathbf{c} \triangleq \begin{pmatrix} \mathbf{z} \\ \mathbf{\bar{z}}
\end{pmatrix} \in \mathbb{C}^{2n} = \mathbb{C}^n \times \mathbb{C}^n \]

Noting that $\mathbf{c}$, $(\mathbf{z},\mathbf{\bar{z}})$,
$\mathbf{z}$, $(\mathbf{x},\mathbf{y})$, and $\mathbf{r}$ are
alternative ways to denote the \emph{same point} $\mathbf{z} =
\mathbf{x} + j \, \mathbf{y}$ in $\mathcal{Z} = \mathbb{C}^n$, for a
function \[ \mathbf{f}:\mathbb{C}^n \rightarrow \mathbb{C}^m
\] throughout this note we will use the convenient (albeit abusive)
notation
\[ \mathbf{f}(\mathbf{c}) = \mathbf{f}
(\mathbf{z},\mathbf{\bar{z}}) = \mathbf{f}(\mathbf{z}) =
\mathbf{f}(\mathbf{x},\mathbf{y}) = \mathbf{f}(\mathbf{r}) \in
\mathbb{C}^m
\]
where $\mathbf{z} = \mathbf{x} + j \, \mathbf{y} \in \mathcal{Z} =
\mathbb{C}^n$.  We will have more to say about the relationships
between these representations later on in Section
\ref{sec:powerseries} below.

\medskip
We further assume that $\mathcal{Z} = \mathbb{C}^n$ is a Riemannian
manifold with a hermitian, positive-definite $n \times n$ metric
tensor $\Omega_\mathbf{z} = \Omega_\mathbf{z}^H
> 0$. This assumption makes every tangent space\footnote{A tangent space at the point $\mathbf{z}$
is the space of all differential displacements, $d\mathbf{z}$, at
the point $\mathbf{z}$ or, alternatively, the space of all velocity
vectors $\mathbf{v} = \frac{d\mathbf{z}}{dt}$ at the point
$\mathbf{z}$. These are equivalent statements because $d\mathbf{z}$
and $\mathbf{v}$ are scaled version of each other, $d\mathbf{z} =
\mathbf{v} dt$. The tangent space $\textsf{T}_z \mathcal{Z} =
\mathbb{C}_\mathbf{z}^n$ is a linear variety in the space
$\mathcal{Z} = \mathbb{C}^n$. Specifically it is a copy of
$\mathbb{C}^n$ affinely translated to the point $\mathbf{z}$,
$\mathbb{C}_\mathbf{z}^n = \{ \mathbf{z} \}+ \mathbb{C}^n$.}
$\textsf{T}_z \mathcal{Z} = \mathbb{C}^n_\mathbf{z}$ a Hilbert space
with inner product
\[ \left<\mathbf{v}_1, \mathbf{v}_2 \right> = \mathbf{v}_1^H \Omega_\mathbf{z} \mathbf{v}_2 \qquad \mathbf{v}_1,\mathbf{v}_2 \in
\mathbb{C}_\mathbf{z}^n . \]

\subsection{The Cogradient Operator and the Jacobian Matrix}

\paragraph{The Cogradient and Conjugate Cogradient.}
Define the \emph{cogradient} and \emph{conjugate cogradient
operators} respectively as the \emph{row operators}\footnote{The ``cogradient'' is a \emph{\underline{co}}variant operator \cite{TF01}. It is \emph{not} itself a gradient, but is the {\emph{\underline{co}}}\,mpanion to the gradient operator defined below.}
\begin{equation}\label{eq:extended}
 \boxed{ \ \text{\sf Cogradient Operator:} \quad
 \mathbf{\frac{\partial \ }{\partial z}} \triangleq \begin{pmatrix} \frac{\partial \  }{\partial z_1}
 & \cdots & \frac{\partial \  }{\partial z_n} \end{pmatrix} \ }
\end{equation}
\begin{equation}\label{eq:extendedb}
\boxed{ \text{\sf Conjugate cogradient Operator:} \quad
\mathbf{\frac{\partial \ }{\partial \bar{z}}} \triangleq
\begin{pmatrix} \frac{\partial \  }{\partial \bar{z}_1} & \cdots
& \frac{\partial \  }{\partial \bar{z}_n} \end{pmatrix} \ }
\end{equation}
where $(z_i,\, \bar{z}_i)$, $i = 1, \cdots, n$ are conjugate
coordinates as discussed earlier and the component operators are
$\mathbb{R}$-derivatives and {\footnotesize $\overline{\,
\mathbb{R}}$}-derivatives defined according to equations
(\ref{eq:e4}) and (\ref{eq:equivdef}),
\begin{equation}\label{eq:partialsa}
\frac{\partial }{\partial z_i} =  \frac{1}{2} \, \left(
\frac{\partial }{\partial x_i} - j \frac{\partial }{\partial
  y_i}\right) \quad \text{and} \quad \frac{\partial }{\partial \bar{z}_i} =  \frac{1}{2} \,
  \left( \frac{\partial }{\partial x_i} + j \frac{\partial }{\partial
  y_i}\right) \, ,
\end{equation}
for $i = 1, \cdots, n$.\footnote{As before the left-hand-sides of (\ref{eq:partialsa}) and
(\ref{eq:partials}) are \emph{formal} partial derivatives, while the
right-hand-sides are \emph{actual} partial derivatives.} Equivalently, we
have
\begin{equation}\label{eq:partials}
\frac{\partial }{\partial \mathbf{z}} =  \frac{1}{2} \, \left(
\frac{\partial }{\partial \mathbf{x}} - j \frac{\partial }{\partial
  \mathbf{y}}\right) \quad \text{and} \quad \frac{\partial }{\partial \mathbf{\bar{z}}} =  \frac{1}{2} \,
  \left( \frac{\partial }{\partial \mathbf{x}} + j \frac{\partial }{\partial
  \mathbf{y}}\right) \, ,
\end{equation}

\medskip
When applying the cogradient operator $\mathbf{\frac{\partial \
}{\partial z}}$, \, $\mathbf{\bar{z}}$ is formally treated as a
constant, and when applying the conjugate cogradient operator
$\mathbf{\frac{\partial \ }{\partial \bar{z}}}$, \, $\mathbf{z}$ is
formally treated as a constant. For example, consider the
scalar-valued function
\[ f (\mathbf{c}) = f(\mathbf{z},\mathbf{\bar{z}}) = z_1 \bar{z}_2 + \bar{z}_1 z_2
\, .
\] For this function we can readily determine by partial differentiation
on the $z_i$ and $\bar{z}_i$ components that
\[ \frac{\partial f(\mathbf{c})
}{\partial \mathbf{z}}  = \begin{pmatrix} \bar{z}_2 & \bar{z}_1
\end{pmatrix}\quad \text{and} \quad
\frac{\partial f(\mathbf{c}) }{\partial \mathbf{\bar{z}}}  =
\begin{pmatrix} z_2 & z_1
\end{pmatrix} \, . \]

\paragraph{The Jacobian Matrix.} Let $\mathbf{f}(\mathbf{c}) =
\mathbf{f}(\mathbf{z},\mathbf{\bar{z}}) \in \mathbb{C}^m$ be a
mapping\footnote{It will always be assumed that the components of
vector-valued functions are $\mathbb{R}$-differentiable as discussed
in footnotes (\ref{foot:a}) and (\ref{foot:b}).}
\[ \mathbf{f}: \mathcal{Z} = \mathbb{C}^n \rightarrow \mathbb{C}^m . \]
The generalization of the identity (\ref{eq:drule}) yields the
\emph{vector form of the differential rule},\footnote{At this point
in our development, the expression $\mathbf{\frac{\partial
\mathbf{f}(\mathbf{c}) }{\partial \mathbf{c}}} d \mathbf{c}$ only
has meaning as a shorthand expression for $\mathbf{\frac{\partial
\mathbf{f}(\mathbf{c}) }{\partial \mathbf{z}}} \, d \mathbf{z} +
\mathbf{\frac{\partial \mathbf{f}(\mathbf{c}) }{\partial
\mathbf{\bar{z}}}} \, d \mathbf{\bar{z}}$,  each term of which must
be interpreted formally as $\mathbf{z}$ and $\mathbf{\bar z}$ cannot
be varied independently of each other. (Later, we will examine the
very special sense in which the a derivative with respect to
$\mathbf{c}$ itself can make sense.) Also note that, unlike the real
case discussed in \cite{TF01}, the mapping $d \mathbf{z} \mapsto d
\mathbf{f}(\mathbf{c})$ is \emph{not} linear in $d \mathbf{z}$. Even
when interpreted formally, the mapping is affine in $d\mathbf{z}$,
not linear.}
\begin{equation}\label{eq:vecdiffrule1}
d \mathbf{f}(\mathbf{c}) =  \mathbf{\frac{\partial
\mathbf{f}(\mathbf{c}) }{\partial \mathbf{c}}} d \mathbf{c} =
\mathbf{\frac{\partial \mathbf{f}(\mathbf{c}) }{\partial
\mathbf{z}}} \, d \mathbf{z} + \mathbf{\frac{\partial
\mathbf{f}(\mathbf{c}) }{\partial \mathbf{\bar{z}}}} \, d
\mathbf{\bar{z}} \, , \qquad \textsf{Differential Rule}
\end{equation}
where  the $m\times n$ matrix $\mathbf{\frac{\partial \mathbf{f}
}{\partial \mathbf{z}}}$ is called the \emph{Jacobian}, or
\emph{Jacobian matrix}, of the mapping $\mathbf{f}$, and the
$m\times n$ matrix $\mathbf{\frac{\partial \mathbf{f}}{\partial
\mathbf{\bar{z}}}}$ the \emph{conjugate Jacobian} of $\mathbf{f}$.
The Jacobian of $\mathbf{f}$ is often denoted by $J_\mathbf{f}$ and
is computed by applying the cogradient operator component-wise to
$\mathbf{f}$,
\begin{equation}\label{eq:jacobiandef} J_\mathbf{f}(\mathbf{c})
= \mathbf{\frac{\partial \mathbf{f}(\mathbf{c}) }{\partial
\mathbf{z}}} =
\begin{pmatrix}
\frac{\partial f_1(\mathbf{c}) }{\partial \mathbf{z}} \\ \vdots \\
\frac{\partial f_n(\mathbf{c}) }{\partial \mathbf{z}}
\end{pmatrix} =
\begin{pmatrix}
\frac{\partial f_1(\mathbf{c}) }{\partial z_1} & \cdots & \frac{\partial f_1(\mathbf{c}) }{\partial z_n} \\ \vdots & \ddots & \vdots \\
\frac{\partial f_n(\mathbf{c}) }{\partial z_1} & \cdots &
\frac{\partial f_n(\mathbf{c}) }{\partial z_n} \end{pmatrix} \in
\mathbb{C}^{m \times n} ,
\end{equation}
and similarly the conjugate Jacobian, denoted by $J_\mathbf{f}^c$ is
computing by applying the conjugate cogradient operator
component-wise to $\mathbf{f}$,
\begin{equation}\label{eq:cojacobiandef}
J_\mathbf{f}^c(\mathbf{c}) = \mathbf{\frac{\partial
\mathbf{f}(\mathbf{c}) }{\partial \mathbf{\bar{z}}}} =
\begin{pmatrix}
\frac{\partial f_1(\mathbf{c}) }{\partial \mathbf{\bar{z}}} \\ \vdots \\
\frac{\partial f_n(\mathbf{c}) }{\partial \mathbf{\bar{z}}}
\end{pmatrix} =
\begin{pmatrix}
\frac{\partial f_1(\mathbf{c}) }{\partial \bar{z}_1} & \cdots & \frac{\partial f_1(\mathbf{c}) }{\partial \bar{z}_n} \\ \vdots & \ddots & \vdots \\
\frac{\partial f_n(\mathbf{c}) }{\partial \bar{z}_1} & \cdots &
\frac{\partial f_n(\mathbf{c}) }{\partial \bar{z}_n} \end{pmatrix}
\in \mathbb{C}^{m \times n} .
\end{equation}
With this notation we can write the differential rule as
\begin{equation}\label{eq:vecdiffrule}
d \mathbf{f}(\mathbf{c}) = J_\mathbf{f}(\mathbf{c}) \, d \mathbf{z}
+ J_\mathbf{f}^c(\mathbf{c})
 \, d \mathbf{\bar{z}}\,  .\qquad
\textsf{Differential Rule}
\end{equation}

Applying properties (\ref{eq:conjusa}) and (\ref{eq:conjusb})
component-wise yields the identities
\begin{equation}\label{eq:veconjus}
\frac{\partial \mathbf{\bar{f}}(\mathbf{c})}{\partial
\mathbf{\bar{z}}}  = \overline{\left(\frac{\partial
\mathbf{f}(\mathbf{c})}{\partial \mathbf{z}}\right)} =
\bar{J}_\mathbf{f}(\mathbf{c}) \qquad \text{and} \qquad
\frac{\partial \mathbf{\bar{f}}(\mathbf{c})}{\partial \mathbf{z}} =
\overline{\left(\frac{\partial \mathbf{f}(\mathbf{c})}{\partial
\mathbf{\bar{z}}}\right)}= \bar{J}^c_\mathbf{f}(\mathbf{c}) \, .
\end{equation}

Note from (\ref{eq:veconjus}) that,
\begin{equation}
\bar{J}_\mathbf{f}(\mathbf{c}) = \overline{\left(
\mathbf{\frac{\partial \mathbf{f}(\mathbf{c}) }{\partial
\mathbf{z}}}\right) } = \frac{\partial \mathbf{\bar{f}}(\mathbf{c})
}{\partial \mathbf{\bar{z}}} \ne J_\mathbf{f}^c(\mathbf{c}) =
\mathbf{\frac{\partial \mathbf{f}(\mathbf{c}) }{\partial
\mathbf{\bar{z}}}} .
\end{equation}
However, in {\emph{the important special case that
$\mathbf{f}(\mathbf{c})$ is real-valued}} (in which case
$\mathbf{\bar{f}}(\mathbf{c}) = \mathbf{f}(\mathbf{c})$) we have
\begin{equation}
\mathbf{f}(\mathbf{c}) \in \mathbb{R}^m \,  \Rightarrow \,
\bar{J}_\mathbf{f}(\mathbf{c}) = \overline{\mathbf{\frac{\partial
\mathbf{f}(\mathbf{c}) }{\partial \mathbf{z}}}} = \frac{\partial
\mathbf{f}(\mathbf{c}) }{\partial \mathbf{\bar{z}}} =
J_\mathbf{f}^c(\mathbf{c})  .
\end{equation}
With (\ref{eq:vecdiffrule}) this yields the following important fact
which holds for real-valued functions
$\mathbf{f}(\mathbf{c})$,\footnote{The real part of a vector (or
matrix) is the  vector (or matrix) of the real parts. Note that the
mapping $d \mathbf{z} \mapsto d \mathbf{f}(\mathbf{c})$ is not
linear.}
\begin{equation}\label{eq:conjdiffrule}
\mathbf{f}(\mathbf{c}) \in \mathbb{R}^m \,  \Rightarrow \, d
\mathbf{f}(\mathbf{c}) = J_\mathbf{f}(\mathbf{c}) \, d \mathbf{z} +
\overline{J_\mathbf{f}(\mathbf{c})
 \, d \mathbf{z}} = 2 \, \text{Re} \left\{ J_\mathbf{f}(\mathbf{c}) \, d
\mathbf{z} \right\} \,  .\qquad
\end{equation}

\medskip
Consider the composition of two mappings $\mathbf{h}: \mathbb{C}^m
\rightarrow \mathbb{C}^r$ and $\mathbf{g}: \mathbb{C}^n \rightarrow
\mathbb{C}^m$,
\[ \mathbf{h} \circ \mathbf{g} = \mathbf{h}(\mathbf{g}) :
\mathbb{C}^n \rightarrow \mathbb{C}^r \, . \] The vector extensions
of the chain rule identities (\ref{eq:ch1}) and (\ref{eq:ch2}) to
$\mathbf{h} \circ \mathbf{g}$ are
\begin{eqnarray}
\frac{\partial \mathbf{h}(\mathbf{g})}{\partial \mathbf{z}} & = &
\frac{\partial \mathbf{h}}{\partial \mathbf{g}}\, \frac{\partial
\mathbf{g}}{\partial \mathbf{z}} + \frac{\partial
\mathbf{h}}{\partial \mathbf{\bar{g}}}\, \frac{\partial
\mathbf{\bar{g}}}{\partial \mathbf{z}} \qquad \quad \text{Chain
Rule} \label{eq:vch1}
\\
\frac{\partial \mathbf{h}(\mathbf{g})}{\partial \mathbf{\bar{z}}} &
= & \frac{\partial \mathbf{h}}{\partial \mathbf{g}}\, \frac{\partial
\mathbf{g}}{\partial \mathbf{\bar{z}}} + \frac{\partial
\mathbf{h}}{\partial \mathbf{\bar{g}}}\, \frac{\partial
\mathbf{\bar{g}}}{\partial \mathbf{\bar{z}}} \qquad \quad
\text{Chain Rule} \label{eq:vch2}
\end{eqnarray}
which can be written as
\begin{eqnarray}
J_{\mathbf{h} \circ \mathbf{g}} & = & J_\mathbf{h}\, J_\mathbf{g} +
J_\mathbf{h}^c\, \bar{J}_\mathbf{g}^c   \label{eq:jch1}
\\
J^c_{\mathbf{h} \circ \mathbf{g}} & = & J_\mathbf{h}\,
J_\mathbf{g}^c + J_\mathbf{h}^c\, \bar{J}_\mathbf{g} \label{eq:jch2}
\end{eqnarray}

\paragraph{Holomorphic Vector-valued Functions.} \emph{By definition} the vector-valued function
$\mathbf{f}(\mathbf{z})$ is holomorphic (analytic in the complex
vector $\mathbf{z}$) if and only if each of its components
\[ f_i(\mathbf{c}) = f_i(\mathbf{z}, \mathbf{\bar{z}}) = f_i(z_1,
\cdots, z_n, \bar{z}_1, \cdots, \bar{z}_n) \quad i = 1, \cdots, m \]
is holomorphic separately with respect to each of the components
$z_j, j = 1, \cdots, n$. In the references
\cite{GH84,LH90,SK01,BS92} it is shown that $\mathbf{f}(\mathbf{z})$
is holomorphic on a domain if and only if it satisfies a matrix
Cauchy Riemann condition everywhere on the domain:
\begin{equation}\label{eq:vcrcond}
\boxed{\ \text{\sf  Cauchy Riemann Condition:} \quad J_\mathbf{f}^c
= \frac{\partial \mathbf{f}}{\partial \mathbf{\bar{z}}} = 0 }
\end{equation}
This shows that \emph{a vector-valued function which is holomorphic
on $\mathbb{C}^n$ must be a function of $\mathbf{z}$ only},
$\mathbf{f}(\mathbf{c}) = \mathbf{f}(\mathbf{z}, \mathbf{\bar{z}}) =
\mathbf{f}(\mathbf{z})$ ($\mathbf z$ only).

\paragraph{Stationary Points of Real-Valued Functionals.} Suppose
that $f$ is a \emph{scalar} real-valued function from $\mathbb{C}^n$
to $\mathbb{R}$,\footnote{The function $f$ is unbolded to indicate
its scalar-value status.}
\[ f:\mathbb{C}^n \rightarrow \mathbb{R} \, ; \ \mathbf{z}
\mapsto f(\mathbf{z}) \, . \] As discussed in \cite{Br83}, the
first-order differential condition for  a real-valued functional $f$
to be optimized with respect to the real and imaginary parts of
$\mathbf{z}$ at the point $\mathbf{z}_0$ is
\begin{equation}\label{eq:sta1}
\boxed{\ \text{\sf Condition I for a Stationary Point:} \quad
\frac{\partial f(\mathbf{z}_0,\mathbf{\bar{z}}_0)}{\partial
\mathbf{z}} = 0 \ }
\end{equation}
That this fact is true is straightforward to ascertain from
equations (\ref{eq:extended}) and (\ref{eq:partials}). An equivalent
first-order condition for a real-valued functional $f$ to be
stationary at the point $\mathbf{z}_0$ is given by
\begin{equation}\label{eq:sta2}
\boxed{ \ \text{\sf Condition II for a Stationary Point:} \quad
\frac{\partial f(\mathbf{z}_0,\mathbf{\bar{z}}_0)}{\partial
\mathbf{\bar{z}}} = 0 \ }
\end{equation}
The equivalence of the two conditions (\ref{eq:sta1}) and
(\ref{eq:sta2}) is a direct consequence of (\ref{eq:veconjus}) and
the fact that $f$ is real-valued.

\paragraph{Differentiation of Conjugate Coordinates?} Note that the use
of the notation $f(\mathbf{c})$ as shorthand for $f(\mathbf{z},
\mathbf{\bar z})$ appears to suggest that it is permissible to take
the complex cogradient of $f(\mathbf{c})$ with respect to the
conjugate coordinates vector $\mathbf{c}$ by treating the complex
vector $\mathbf{c}$ \emph{itself} as the variable of
differentiation. \emph{This is not correct.} Only complex
differentiation with respect to the complex vectors $\mathbf{z}$ and
$\mathbf{\bar z}$ is well-defined. Thus, from the definition
$\mathbf{c} \triangleq \text{col}(\mathbf{z},\mathbf{\bar z})\in
\mathbb{C}^{2n}$, for $\mathbf{c}$ \emph{viewed as a complex
$2n$-dimensional vector,}
 the correct interpretation of
$\frac{\partial}{\partial \mathbf{c}} f(\mathbf{c})$ is given by
\[ \frac{\partial}{\partial \mathbf{c}} f(\mathbf{c}) = \left[\frac{\partial}{\partial
\mathbf{z}} f(\mathbf{z},\mathbf{\bar z}) \, , \,
\frac{\partial}{\partial \mathbf{\bar z}} f(\mathbf{z},\mathbf{\bar
z}) \right] \] Thus, for example, we have that
\[ \frac{\partial}{\partial \mathbf{c}} \mathbf{c}^H \Omega
\mathbf{c}\,  \ \text{\boldmath $\ne$} \ \, \mathbf{c}^H \Omega \]
which would be true \textbf{if} it were permissible to take the
complex cogradient with respect to the complex vector $\mathbf{c}$
(which it isn't).

\medskip
Remarkably, however, below we will show that the $2n$-dimensional
complex vector $\mathbf{c}$ is an element of an $n$-dimensional real
vector space and that, as a consequence, it is permissible to take
the real cogradient with respect to the \emph{real} vector
$\mathbf{c}$!

\paragraph{Comments.}
\emph{With the machinery developed up to this point, one can solve
optimization problems which have closed-form solutions to the
first-order stationarity conditions.}  However, to solve general
nonlinear problems one must often resort to gradient-based iterative
methods.  Furthermore, to verify that the solutions are optimal, one
needs to check  second order conditions which require the
construction of the hessian matrix. Therefore, the remainder of this
note is primarily concerned with the development of the machinery
required to construct the gradient and hessian of a scalar-valued
functional of complex parameters.

\subsection{Biholomorphic Mappings and Change of Coordinates.}

\paragraph{Holomorphic and Biholomorphic Mappings.} A
vector-valued function $\mathbf{f}$ is holomorphic (complex-analytic) if its
components are holomorphic. In this case the function does not
depend on the conjugate coordinate $\mathbf{\bar{z}}$,
$\mathbf{f}(\mathbf{c}) = \mathbf{f}(\mathbf{z})$ ($\mathbf{z}$-only), and satisfies the
Cauchy-Riemann Condition,
\[ J_\mathbf{f}^c = \frac{\partial
\mathbf{f}}{\partial \mathbf{\bar{z}}} = 0 \, . \] As a consequence
(see (\ref{eq:vecdiffrule})),
\begin{equation}\label{eq:vecdiffrulehol}
\mathbf{f}(\mathbf{z})\text{ holomorphic }\ \Rightarrow  \ d
\mathbf{f}(\mathbf{z}) = J_\mathbf{f}(\mathbf{z}) \, d \mathbf{z} =
\mathbf{\frac{\partial \mathbf{f}(\mathbf{z}) }{\partial
\mathbf{z}}} \, d \mathbf{z} \, .
\end{equation}
Note that when $\mathbf{f}$ is holomorphic, the mapping $d
\mathbf{z} \mapsto d \mathbf{f}(\mathbf{z})$ is linear, exactly as
in the real case.

\medskip
Consider the composition of two mappings $\mathbf{h}: \mathbb{C}^m
\rightarrow \mathbb{C}^r$ and $\mathbf{g}: \mathbb{C}^n \rightarrow
\mathbb{C}^m$,
\[ \mathbf{h} \circ \mathbf{g} = \mathbf{h}(\mathbf{g}) :
\mathbb{C}^n \rightarrow \mathbb{C}^r \, , \] which are \emph{both
holomorphic}. In this case, as a consequence of the Cauchy-Riemann
condition (\ref{eq:vcrcond}), the second chain rule condition
(\ref{eq:jch2}) vanishes, $J^c_{\mathbf{h} \circ \mathbf{g}} = 0$,
and the first chain rule condition (\ref{eq:jch1}) simplifies to
\begin{equation}
\mathbf{f}\text{ and } \mathbf{g}\text{ holomorphic } \Rightarrow
J_{\mathbf{h} \circ \mathbf{g}} = J_\mathbf{h}\, J_\mathbf{g} \, .
\end{equation}

\medskip
Now consider the holomorphic mapping  $\text{\boldmath $\xi$} =
\mathbf{f}(\mathbf{z})$,
\begin{equation} d \text{\boldmath $\xi$} = d \mathbf{f}(\mathbf{z}) =
J_\mathbf{f}(\mathbf{z}) \, d \mathbf{z} \,
\end{equation}
and assume that it is invertible,
\begin{equation} \mathbf{z} =
\mathbf{g}(\text{\boldmath $\xi$}) = \mathbf{f}^{-1}
(\text{\boldmath $\xi$}) \, .
\end{equation}
If the invertible function $\mathbf{f}$ and its inverse $\mathbf{g}
= \mathbf{f}^{-1}$ are \emph{both} holomorphic, then $\mathbf{f}$
(equivalently, $\mathbf{g}$) is said to be \emph{biholomorphic}. In
this case, we have that
\begin{equation}\label{eq:biholin}
d \mathbf{z} = \frac{\partial \mathbf{g}(\text{\boldmath
$\xi$})}{\partial \text{\boldmath $\xi$}} \, d \text{\boldmath
$\xi$} = J_\mathbf{g}(\text{\boldmath $\xi$}) \,  d \text{\boldmath
$\xi$} =  J_\mathbf{f}^{-1}(\mathbf{z}) \, d \text{\boldmath $\xi$}
\, , \qquad \text{\boldmath $\xi$} = \mathbf{f}(\mathbf{z}) \, ,
\end{equation}
showing that
\begin{equation}
J_\mathbf{g}(\text{\boldmath $\xi$}) =
J_\mathbf{f}^{-1}(\mathbf{z})\, , \qquad \text{\boldmath $\xi$} =
\mathbf{f}(\mathbf{z}) \, .
\end{equation}

\paragraph{Coordinate Transformations.}  Admissible coordinates on
a space defined over a space of complex numbers are related via
biholomorphic transformations \cite{GH84,LH90,SK01,BS92}. Thus if
$\mathbf{z}$ and $\text{\boldmath $\xi$}$ are admissible coordinates
on $\mathcal{Z} = \mathbb{C}^n$, there \emph{must} exist a
biholomorphic mapping relating the two coordinates, $\text{\boldmath
$\xi$}= \mathbf{f}(\mathbf{z})$. This relationship is often denoted
in the following (potentially confusing) manner,
\begin{equation}\text{\boldmath $\xi$}=
\text{\boldmath $\xi$}(\mathbf{z}) \, , \quad d \text{\boldmath
$\xi$} = \mathbf{\frac{\partial \text{\boldmath $\xi$}(\mathbf{z})
}{\partial \mathbf{z}}} \, d \mathbf{z}= J_{\text{\boldmath
$\xi$}}(\mathbf{z})\,  d \mathbf{z} \, , \quad
\mathbf{\frac{\partial \text{\boldmath $\xi$}(\mathbf{z}) }{\partial
\mathbf{z}}} = J_{\text{\boldmath $\xi$}}(\mathbf{z}) =
J_{\mathbf{z}}^{-1}(\text{\boldmath $\xi$}) = \left( \frac{\partial
\mathbf{z}(\text{\boldmath $\xi$})}{\partial \text{\boldmath $\xi$}}
\right)^{-1}
\end{equation}
\begin{equation}
\mathbf{z} = \mathbf{z}(\text{\boldmath $\xi$}) \, , \quad d
\mathbf{z} = \frac{\partial \mathbf{z}(\text{\boldmath
$\xi$})}{\partial \text{\boldmath $\xi$}} \, d \text{\boldmath
$\xi$} = J_{\mathbf{z}}(\text{\boldmath $\xi$}) \, d \text{\boldmath
$\xi$} \, , \quad \frac{\partial \mathbf{z}(\text{\boldmath
$\xi$})}{\partial \text{\boldmath $\xi$}} =
J_{\mathbf{z}}(\text{\boldmath $\xi$}) =
 J_{\text{\boldmath
$\xi$}}^{-1}(\mathbf{z}) = \left( \mathbf{\frac{\partial
\text{\boldmath $\xi$}(\mathbf{z}) }{\partial \mathbf{z}}}
\right)^{-1} ,
\end{equation}
These equations tell us how vectors (elements of any particular
tangent space $\mathbb{C}_z^n$) properly transform under a change of
coordinates.

\medskip
In particular under the change of coordinates $\text{\boldmath
$\xi$} = \text{\boldmath $\xi$}(\mathbf{z})$, a vector $\mathbf{v}
\in \mathbb{C}^n_\mathbf{z}$ must transform to its new
representation $\mathbf{w} \in \mathbb{C}^n_{\text{\boldmath
$\xi$}(\mathbf{z})}$ according to the
\begin{equation}
\boxed{\ \text{\sf Vector Transformation Law:} \quad \mathbf{w} =
\frac{\partial \text{\boldmath $\xi$}}{\partial \mathbf{z}} \,
\mathbf{v} = J_{\text{\boldmath $\xi$}} \, \mathbf{v} \ }
\end{equation}

For the composite coordinate transformation $\text{\boldmath
$\eta$}(\text{\boldmath $\xi$}(\mathbf{z}))$, the chain rule yields
\begin{equation}
\boxed{ \ \text{\sf Transformation Chain Rule:} \quad \frac{\partial
\text{\boldmath $\eta$}}{\partial \mathbf{z}} = \frac{\partial
\text{\boldmath $\eta$}}{\partial \text{\boldmath $\xi$}} \,
\mathbf{\frac{\partial \text{\boldmath $\xi$}}{\partial \mathbf{z}}}
\quad \text{or} \quad J_{\text{\boldmath $\eta$} \circ
\text{\boldmath $\xi$}} = J_{\text{\boldmath $\eta$}} \,
J_{\text{\boldmath $\xi$}} \ }
\end{equation}

Finally, applying the chain rule to the cogradient, $\frac{\partial
\mathbf{f} }{\partial \emph{z}}$,  of a an arbitrary holomorphic
function $\mathbf{f}$ we obtain
\[ \frac{\partial \mathbf{f}
}{\partial \text{\boldmath $\xi$}} = \frac{\partial \mathbf{f}
}{\partial \mathbf{z}} \, \frac{ \partial \mathbf{z}}{\partial
\text{\boldmath $\xi$}}  \quad \text{for} \quad \text{\boldmath
$\xi$} = \text{\boldmath $\xi$}(\mathbf{z}) \, . \] This shows that
the cogradient, \emph{as an operator on holomorphic functions,}
transforms like
\begin{equation}
\boxed{ \ \text{\sf Cogradient Transformation Law:} \quad
\frac{\partial (\, \cdot \, )}{\partial \text{\boldmath $\xi$}} =
\frac{\partial (\, \cdot \, ) }{\partial \mathbf{z}} \,  \frac{
\partial \mathbf{z}}{\partial \text{\boldmath $\xi$}}=
\frac{\partial (\, \cdot \, ) }{\partial \mathbf{z}}\,
J_{\mathbf{z}} = \frac{\partial (\, \cdot \, ) }{\partial
\mathbf{z}}\, J_{\text{\boldmath $\xi$}}^{-1} \ }
\end{equation}
Note that generally the cogradient transforms quite differently than
does a vector.

\medskip Finally the transformation law for the metric tensor under
a change of coordinates can be determined from the requirement that
the inner product must be  invariant under a change of coordinates.
For arbitrary vectors $\mathbf{v}_1,\mathbf{v}_2 \in
\mathbb{C}^n_{\mathbf{z}}$ transformed as \[\mathbf{w}_i =
J_{\text{\boldmath $\xi$}} \, \mathbf{v}_i \in
\mathbb{C}^n_{\text{\boldmath $\xi$}(\mathbf{z})} \, \quad i = 1,2
\, , \] we have
\[ \left< \mathbf{w}_1, \mathbf{w}_2 \right> = \mathbf{w}_1^H\,
\Omega_{\text{\boldmath $\xi$}} \,\mathbf{w}_2 = \mathbf{v}_1^H \,
J^H_{\text{\boldmath $\xi$}} \, \Omega_{\text{\boldmath $\xi$}}\,
J_{\text{\boldmath $\xi$}} \, \mathbf{v}_2 = \mathbf{v}_1^H \,
J^{-H}_{\mathbf{z}} \, \Omega_{\text{\boldmath $\xi$}}\,
J_{\mathbf{z}} \, \mathbf{v}_2 = \mathbf{v}_1^H \,
 \Omega_{\mathbf{z}}\, \mathbf{v}_2  = \left< \mathbf{v}_1, \mathbf{v}_2 \right> \, . \]
This results in the
\begin{equation}
\boxed{ \ \text{\sf Metric Tensor Transformation Law:} \quad
\Omega_{\text{\boldmath $\xi$}} = J^{-H}_{\text{\boldmath $\xi$}}\,
\Omega_{\mathbf{z}} \, J^{-1}_{\text{\boldmath $\xi$}}
 = J^H_{\mathbf{z}}\, \Omega_{\mathbf{z}} \, J_{\mathbf{z}} \ }
\end{equation}

\section{\boldmath The Gradient Operator $\mathbf{\nabla_z}$}
\paragraph{\boldmath $1^{st}$-Order Approximation of a Real-Valued
Function.} Let $f(\mathbf{c})$ be a \emph{real-valued scalar}\footnote{And
therefore unbolded.} functional to be optimized with respect to the
real and imaginary parts of the vector $\mathbf{z} \in \mathcal{Z} =
\mathbb{C}^n$,
\[ f:  \mathbb{C}^n \rightarrow \mathbb{R} \, . \]
As a \emph{real-valued} function, $f(\mathbf{c})$ does not satisfy
the Cauchy-Riemann condition (\ref{eq:vcrcond}) and is therefore not
holomorphic.

\medskip
From (\ref{eq:conjdiffrule}) we have (with $f(\mathbf{z}) =
f(\mathbf{z}, \mathbf{\bar z}) = f(\mathbf{c})$) that
\begin{equation}\label{eq:sconjdiffrule}
d f (\mathbf{z}) =  2 \, \text{Re} \left\{ J_f(\mathbf{z}) \, d
\mathbf{z} \right\} = 2\,  \text{Re} \left\{ \frac{\partial
f(\mathbf{z})}{\partial \mathbf{z}} \, d \mathbf{z} \right\}\, .
\end{equation}
This yields the first order relationship
\begin{equation}\label{eq:dlinexpansion}
f(\mathbf{z} + d \mathbf{z}) = f(\mathbf{z}) + 2\, \text{Re} \left\{
\frac{\partial f(\mathbf{z})}{\partial \mathbf{z}} \, d \mathbf{z}
\right\}
\end{equation}
and the corresponding first-order power series approximation
\begin{equation}\label{eq:deltalinexpansion}
f(\mathbf{z} + \Delta \mathbf{z}) \approx f(\mathbf{z}) + 2\,
\text{Re} \left\{ \frac{\partial f(\mathbf{z})}{\partial \mathbf{z}}
\, \Delta \mathbf{z} \right\}
\end{equation}
which will be rederived by other means in Section
\ref{sec:powerseries} below.

\medskip
\paragraph{The Complex Gradient of a Real-Valued Function.} The
relationship (\ref{eq:sconjdiffrule}) defines a \emph{nonlinear}
functional, $d f_{\mathbf{c}}(\cdot)$, on the tangent space
$\mathbb{C}^n_\mathbf{z}$,\footnote{Because this operator is
\emph{nonlinear} in $d \mathbf{z}$, unlike the real vector-space
case
\cite{TF01}, we will avoid calling it a ``differential
operator.''.}
\begin{equation}\label{eq:diff1}
d f_{\mathbf{c}}(\mathbf{v}) =  2\,  \text{Re} \left\{
\frac{\partial f(\mathbf{c})}{\partial \mathbf{z}} \, \mathbf{v}
\right\} \, , \qquad \mathbf{v} \in \mathbb{C}^n_\mathbf{z} \, , \
\mathbf{c} = (\mathbf{z}, \mathbf{\bar{z}}) \, .
\end{equation}
Assuming the existence of a metric tensor $\Omega_{\mathbf{z}}$ we
can write
\begin{equation}\label{eq:diff2}  \frac{\partial f}{\partial
\mathbf{z}} \, \mathbf{v} = \left[ \Omega_{\mathbf{z}}^{-1}
\left(\frac{\partial f}{\partial \mathbf{z}}\right)^H \right]^H \,
\Omega_{\mathbf{z}} \, \mathbf{v} = \left(\nabla_{\mathbf{z}}
f\right)^H \Omega_{\mathbf{z}} \, \mathbf{v} = \left<
\nabla_{\mathbf{z}} f,  \, \mathbf{v} \right> \, ,
\end{equation}
where $\nabla_{\mathbf{z}} f$ is the \emph{gradient} of $f$, defined
as
\begin{equation}\label{eq:grad1} \boxed{ \ \text{\sf Gradient of
$f$:} \quad \nabla_{\mathbf{z}} f \triangleq
\Omega_{\mathbf{z}}^{-1} \left(\frac{\partial f}{\partial
\mathbf{z}}\right)^H \ }
\end{equation}
Consistent with this definition, the \emph{gradient operator} is defined as
\begin{equation}\label{eq:grad2}
\boxed{ \ \text{\sf Gradient Operator:} \quad \nabla_{\mathbf{z}}
(\, \cdot \, ) \triangleq \Omega_{\mathbf{z}}^{-1}
\left(\frac{\partial (\, \cdot \, )}{\partial \mathbf{z}}\right)^H \
}
\end{equation}
Note the relationships between the gradients and the cogradients.
One can show from the coordinate transformation laws for cogradients
and metric tensors that \emph{the gradient $\nabla_{\mathbf{z}} f$
transforms like a vector} and therefore \emph{is} a vector,
\[ \nabla_{\mathbf{z}} f \in \mathbb{C}^n_{\mathbf{z}} \, .
\]

Equations (\ref{eq:diff1}) and (\ref{eq:diff2}) yield,
\[ d f_{\mathbf{c}}(\mathbf{v}) =  2\, \text{Re} \left\{ \left<
\nabla_{\mathbf{z}} f,  \, \mathbf{v} \right> \right\} \, . \]
Keeping $\| \mathbf{v} \| = 1$ we want to find the directions
$\mathbf{v}$ of steepest increase in the value of $ \left| d
f_{\mathbf{c}}(\mathbf{v}) \right|$. We have as a consequence of the
Cauchy-Schwarz inequality that for all unit vectors $v \in
\mathbb{C}^n_{\mathbf{z}}$,
\[ \left| d
f_{\mathbf{c}}(\mathbf{v}) \right| =  2\, \left| \text{Re} \left\{
\left< \nabla_{\mathbf{z}} f,  \, \mathbf{v} \right> \right\}
\right| \le 2 \, \left| \left< \nabla_{\mathbf{z}} f,  \, \mathbf{v}
\right> \right| \le 2 \, \| \nabla_{\mathbf{z}} f \| \, \|
\mathbf{v} \| = 2 \, \| \nabla_{\mathbf{z}} f \| \, . \] This upper
bound is attained if and only if $\mathbf{v} \propto
\nabla_{\mathbf{z}} f$, showing that the gradient gives the
directions of steepest increase, with $+ \nabla_{\mathbf{z}} f$
giving the direction of \emph{steepest ascent} and
$-\nabla_{\mathbf{z}} f$ giving the direction of \emph{steepest
descent.} The result (\ref{eq:grad2}) is derived in \cite{Br83} for
the special case that the metric is Euclidean $\Omega_{\mathbf{z}} =
I$.\footnote{Therefore one must be careful to ascertain when a
result derived in \cite{Br83} holds in the general case. Also note
the corresponding notational difference between this note and \cite{Br83}. We have
$\nabla_\mathbf{z}$ denoting the gradient operator for the general case $\Omega_\mathbf{z} \ne I$ while \cite{Br83}
denotes the gradient operator as $\nabla_\mathbf{\bar{z}}$ for the special case
$\Omega_\mathbf{z} = I$.}

\medskip
Note that the first-order necessary conditions for a stationary
point to exist is given by $\nabla_{\mathbf{z}} f = 0$, but that it
is much easier to apply the simpler condition $\frac{\partial
f}{\partial \mathbf{z}} = 0$ which does not require knowledge of the
metric tensor. Of course this distinction vanishes when
$\Omega_\mathbf{z} = I$ as is the case in \cite{Br83}.

\paragraph{\boldmath Comments on Applying the Multivariate  {\footnotesize $\mathbb{CR}$}-Calculus.} Because
the components of the cogradient and conjugate cogradient operators
(\ref{eq:extended}) and (\ref{eq:extendedb}) formally behave like
partial derivatives of functions over real vectors, to use them does
\emph{not} require the development of additional vector partial-derivative
identities over and above those that already exist for the real
vector space case. Real vector space identities and procedures
for vector partial-differentiation carry over without change, \emph{provided one first
carefully distinguishes between those variables which are to be
treated like constants and those variables which are to be formally
differentiated.}

\medskip
Thus, although a variety of complex derivative identities are given
in various references \cite{Br83,Ka93,Ha96}, there is actually \emph{no
need} to memorize or look up additional ``complex derivative
identities'' if one already knows the real derivative identities. In
particular, the derivation of the complex derivative identities
given in  references \cite{Br83,Ka93,Ha96}  is trivial if one
already knows the standard real-vector derivative identities.
 For example, it is
\emph{obviously} the case that
\[ \frac{\partial \ }{\partial \mathbf{\bar{z}}} \left( \mathbf{a}^H \mathbf{z} \right) =
\mathbf{a}^H \, \frac{\partial \mathbf{z}}{\partial
\mathbf{\bar{z}}} = 0 \, ,
\] \emph{as $\mathbf{z}$ is to be treated as a constant} when taking partial derivatives
with respect to $\mathbf{\bar{z}}$.  Therefore the fact that $\frac{\partial
\ }{\partial \mathbf{\bar{z}}}  \mathbf{a}^H \mathbf{z}  = 0$ does
\emph{not} have to be memorized as a special complex derivative
identity.

\medskip
To reiterate, if one already knows the standard gradient
identities for real-valued functions of real variables, \emph{there
is no need to memorize additional complex derivative
identities}.\footnote{This extra emphasis is made because virtually
all of the textbooks (even the exemplary text \cite{Ka93}) provide
such extended derivative identities and use them to derive results.
This sends the message that unless such identities are at hand, one
cannot solve problems. Also, it places one at the mercy of
typographical errors which may occur when identities are printed in
the textbooks.} Instead, one can merely use the regular real
derivative identities
 \emph{while keeping track of which complex variables are to be treated
as constants.}\footnote{Thus, in the real case, $\mathbf{x}$ is the
variable to be differentiated in $\mathbf{x}^T \mathbf{x}$ and we
have $\frac{\partial}{\partial \mathbf{x}} \mathbf{x}^T \mathbf{x} =
2 \mathbf{x}^T$, while in the complex case, if we take $\mathbf{\bar
z}$ \emph{to be treated as constant} and $\mathbf{z}$ to be the
differentiated variable, we have $\frac{\partial}{\partial
\mathbf{z}} \mathbf{z}^H \mathbf{z} = \mathbf{z}^H
\frac{\partial}{\partial \mathbf{z}}  \mathbf{z} = \mathbf{z}^H$.
Note that in both cases we use the differentiation rules for vector
differentiation which are developed initially for the purely real
case once we have decided \emph{which variables are to be treated as
constant.}} This is the approach used to easily derive the complex
LMS algorithm in the applications section at the end of this note.

\medskip
To implement a true gradient descent algorithm, one needs to know
the metric tensor. The correct gradient, which depends on the metric
tensor, is called the ``natural gradient'' in \cite{Am98} where it
is argued that superior performance of gradient descent algorithms
in certain statistical parameter estimation problems occurs when the
natural gradient is used in lieu of the standard ``naive''
gradient usually used in such algorithms (where ``naive'' corresponds to assuming that $\Omega_z = I$ even if that is not the case). However, the determination of the metric tensor for a
specific application can be highly nontrivial and the resulting
algorithms significantly more complex, as discussed in \cite{Am98},
although there are cases where the application of the natural
gradient methodology is surprisingly straightforward.

\medskip
To close this section, we mention that interesting and useful
applications of the {\footnotesize $\mathbb{CR}$}-calculus as
developed in \cite{Br83} and \cite{VDB94a} can be found in
references \cite{Fi02}, \cite{VDB94b}-\cite{AJ04}, and \cite{TA91},
in addition to the plentiful material to be found in the textbooks
\cite{VT02}, \cite{Ka93}, \cite{Ha96}, and \cite{AS03}.

\section{\boldmath $2^{nd}$-Order Expansions of a Real-Valued Function on $\mathbb{C}^n$}\label{sec:powerseries}

It is common to numerically optimize cost functionals using
iterative gradient descent-like techniques.
Determination of the gradient of a real-valued loss function via
equation (\ref{eq:grad1}) allows the use of elementary gradient
descent optimization, while the linear approximation of a
biholomorphic mapping $\mathbf{g}(\text{\boldmath $\xi$})$ via
(\ref{eq:biholin}) enables optimization of the nonlinear
least-squares problem using the Gauss-Newton
algorithm.\footnote{Recall that the Gauss-Newton algorithm is based
on iterative re-linearization of a nonlinear model $\mathbf{z}
\approx \mathbf{g}(\text{\boldmath $\xi$})$.}

\medskip
Another commonly used iterative algorithm is the Newton method,
which is based on the repeated computation and optimization of the
quadratic approximation to the loss function as given by a power
series expansion to second order. Although the first-order approximation to the loss function given by
(\ref{eq:deltalinexpansion}) was relatively straight-forward to
derive, it is somewhat more work to determine the second order
approximation, which is the focus of this section and which will be
attacked using the elegant approach of van den Bos
\cite{VDB94a}.\footnote{A detailed exposition of the second order
case is given by Abatzoglou, Mendel, \& Harada in \cite{TA91}. See
also \cite{GY00}. The references \cite{TA91}, \cite{VDB94a} and
\cite{GY00} all develop the complex Newton algorithm, although with
somewhat different notation.} Along the way we will rederive the
first order approximation (\ref{eq:deltalinexpansion}) and the
Hessian matrix of second partial derivatives of a real scalar-valued
function which is needed to verify the optimality of a solution
solving the first order necessary conditions.

\subsection{\boldmath Alternative Coordinate Representations of $\mathcal{Z}
= \mathbb{C}^n.$}

\paragraph{\boldmath Conjugate Coordinate Vectors $\mathbf{c} \in \mathcal{C}$ Form a Real Vector Space.} The complex space, $\mathbb{C}^n$, of dimension $n$ naturally has
the structure of a real space, $\mathbb{R}^{2n}$, of dimension $2n$,
$\mathbb{C}^n \approx \mathbb{R}^{2n}$, as a consequence of the
equivalence
\[ \mathbf{z} = \mathbf{x} + j \, \mathbf{y} \in \mathcal{Z} = \mathbb{C}^n
\Leftrightarrow \mathbf{r} = \begin{pmatrix} \mathbf{x} \\
\mathbf{y} \end{pmatrix} \in \mathcal{R} \triangleq \mathbb{R}^{2n}
.
\]

\medskip
Furthermore, as noted earlier, an alternative representation is
given by the set of conjugate coordinate vectors
\[ \mathbf{c} = \begin{pmatrix} \mathbf{z} \\ \mathbf{\bar z}
\end{pmatrix} \in \mathcal{C} \subset \mathbb{C}^{2n}\approx \mathbb{R}^{4n} \, , \]
\emph{where $\mathcal C$ is defined to be the collection of all such
vectors $\mathbf{c}$.} Note that the set $\mathcal{C}$ is obviously
a subset (and \emph{not} a vector subspace)\footnote{It is, in fact,
a $2n$ dimensional submanifold of the space $\mathbb{C}^{2n} \approx
\mathbb{R}^{4n}$.} of the $4n$ dimensional complex vector space $\mathbb{C}^{2n}$.
Remarkably, \emph{it is also a $2n$ dimensional vector space over
the field of {real} numbers!}

\medskip
This is straightforward to show. First, in the obvious manner, one
can define vector addition of any two elements of $\mathcal{C}$. To
show closure under scalar multiplication by a \emph{real} number
$\alpha$ is also straight forward,
\[ \mathbf{c} = \begin{pmatrix} \mathbf{z} \\ \mathbf{\bar z}
\end{pmatrix} \in \mathcal{C} \Rightarrow \alpha \, \mathbf{c} =
\begin{pmatrix} \alpha \,  \mathbf{z} \\ \mathbf{\overline{\alpha \, z}}
\end{pmatrix} \in \mathcal{C} \, . \] Note that this homogeneity property
obviously fails when $\alpha$ is complex.

\medskip
To demonstrate that $\mathcal{C}$ is $2n$ dimensional, we will
construct below the one-to-one transformation, $\mathsf J$, which
maps $\mathcal{C}$ onto $\mathcal{R}$, and vice versa, thereby
showing that $\mathcal{C}$ and $\mathcal{R}$ are isomorphic,
$\mathcal{C} \simeq \mathcal{R}$. In this manner $\mathcal{C}$ and
$\mathcal{R}$ are shown to be alternative, \emph{but entirely
equivalent} (including their dimensions), real coordinate
representations for $\mathcal{Z} = \mathbb{C}^n$. The coordinate
transformation $\mathsf J$ is a linear mapping, and therefore also
corresponds to the Jacobian of the transformation between the
coordinate system $\mathcal{R}$ and the coordinate system
$\mathcal{C}$.

\medskip
In summary, we have available three \emph{vector space} coordinate
representations for representing complex vectors $\mathbf{z} =
\mathbf{x} + j \, \mathbf{y}$. The first is the canonical
$n$-dimensional vector space of complex vectors $\mathbf{z} \in
\mathcal{Z} = \mathbb{C}^n$ itself. The second is the canonical
$2n$-dimensional real vector space of vectors $\mathbf{r} =
\text{col}(\mathbf{x},\mathbf{y}) \in \mathcal{R} =
\mathbb{R}^{2n}$, which arises from the natural correspondence
$\mathbb{C}^n \approx \mathbb{R}^{2n}$. The third is the
$2n$-dimensional real vector space of vectors $\mathbf{c}  \in
\mathcal{C} \subset \mathbb{C}^{2n},$ $\mathcal{C} \approx
\mathbb{R}^{2n}$.

\medskip
Because $\mathcal{C}$ can be alternatively viewed as a complex
subset of $\mathbb{C}^{2n}$ or as a real vector space isomorphic to
$\mathbb{R}^{2n}$, we actually have a fourth ``representation'';
namely the \emph{non-vector space} complex-vector perspective of
elements of $\mathcal{C}$ as elements of the space
$\mathbb{C}^{2n}$, $\mathbf{c} = \text{col}(\mathbf{z},\mathbf{\bar
z})$.\footnote{Since when viewed as a subset of $\mathbb{C}^{2n}$
the set $\mathcal{C}$ is \emph{not} a subspace, this view of
$\mathcal{C}$ does not result in a true \emph{coordinate}
representation.} This perspective is just the
$(\mathbf{z},\mathbf{\bar z})$ perspective used above to analyze
general, possibly nonholomorphic, functions $f(\mathbf{z}) =
f(\mathbf{z},\mathbf{\bar z})$.

\medskip
In order to avoid confusion, we will refer to these two alternative
interpretations of $\mathbf{c} \in \mathcal{C} \subset
\mathbb{C}^{2n}$ as the $\mathbf{c}$-real case (respectively, the
$\mathcal{C}$-real case) for when we consider the vector $\mathbf{c}
\in \mathcal{C} \approx \mathbb{R}^{2n}$ (respectively, the real
vector space $\mathcal{C} \approx \mathbb{R}^{2n}$), and the
$\mathbf{c}$-complex case (respectively, the $\mathcal C$-complex
case) when we consider a vector $\mathbf{c} \in \mathcal{C} \subset
\mathbb{C}^{2n}$ (respectively, the complex subset $\mathcal{C}
\subset \mathbb{C}^{2n}$).\footnote{In the latter case $\mathbf{c} =
\text{col}(\mathbf{z},\mathbf{\bar{z}})$ is understood in terms of
the behavior and properties of its components, especially for
differentiation purposes because, as mentioned earlier, in the
\emph{complex case} the derivative $\frac{\partial}{\partial
\mathbf{c}}$ is not well-defined in itself, but is defined in terms
of the formal derivatives with respect to $\mathbf z$ and
$\mathbf{\bar z}$. As we shall discover below, in the
$\mathbf{c}$-real case, the derivative $\frac{\partial}{\partial
\mathbf{c}}$ is a true real derivative which is well understood in
terms of the behavior of the derivative $\frac{\partial}{\partial
\mathbf{r}}$. } These two different perspectives of $\mathcal{C}$
are used throughout the remainder of this note.

\paragraph{Coordinate Transformations and Jacobians.} From the fact that
\[ \mathbf{z} = \mathbf{x} + j \, \mathbf{y} \quad \text{and} \quad
\mathbf{\bar z} = \mathbf{x} - j \, \mathbf{y} \] it is easily shown
that
\[ \begin{pmatrix} \mathbf{z} \\ \mathbf{\bar z} \end{pmatrix} =
 \begin{pmatrix} I & \quad j\, I \\ I & - j \, I \end{pmatrix}
 \begin{pmatrix} \mathbf{x} \\ \mathbf{y} \end{pmatrix} \]
 where $I$ is the $n \times n$ identity matrix.
 Defining\footnote{Except for a trivial reordering of the elements
 of $\mathbf{r} = (\mathbf{x}^T \, \mathbf{y}^T )^T$, this is the transformation proposed and utilized by van den Bos
 \cite{VDB94a}, who claims in \cite{VDB96} to have been inspired to
 do so by Remmert. (See, e.g., the discussion on page 87 of
 \cite{RR91}.)}

\begin{equation}\label{eq:bigj}
\mathsf{J} \triangleq \begin{pmatrix} I & \quad j\, I \\ I & - j \,
I \end{pmatrix}
\end{equation}
then results in the mapping
\begin{equation}\label{eq:mapforward}
\mathbf{c} = \mathbf{c}(\mathbf{r}) = \mathsf{J} \, \mathbf{r}\, .
\end{equation}
 It is
easily determined that
\begin{equation}\label{eq:bigjinv}
\mathsf{J}^{-1} = \frac{1}{2} \mathsf{J}^H \,
\end{equation}
so that we have the inverse mapping
\begin{equation}
\label{eq:mapback} \mathbf{r} = \mathbf{r}(\mathbf{c}) =
\mathsf{J}^{-1} \mathbf{c} = \frac{1}{2} \mathsf{J}^H \mathbf{c} \,
.
\end{equation}

\medskip
Because the mapping between $\mathcal{R}$ and $\mathcal{C}$ is
linear, one-to-one, and onto, both of these spaces are obviously
isomorphic real vector spaces of dimension $2n$. The mappings
(\ref{eq:mapforward}) and (\ref{eq:mapback}) therefore correspond to
an admissible coordinate transformation between the $\mathbf{c}$ and
$\mathbf{r}$ representations of $\mathbf{z} \in \mathcal{Z}$.
Consistent with this fact, we henceforth assume that the real vector
calculus (including all of the vector derivative identities) apply to real-valued functions over $\mathcal{C}$.

\medskip
Note that for the coordinate transformation $\mathsf{c} =
\mathbf{c}(\mathsf{r}) = \mathsf{J}  \mathbf{r}$ we have the
Jacobian
\begin{equation}
J_{\mathbf{c}} \triangleq \frac{\partial}{\partial \mathbf{r}}
\,\mathbf{c}(\mathsf{r}) = \frac{\partial}{\partial \mathbf{r}} \,
\mathsf{J} \mathbf{r} = \mathsf{J}
\end{equation}
showing that $\mathsf{J}$ is also the Jacobian of the coordinate
transformation from $\mathcal{R}$ to $\mathcal{C}$.\footnote{We have
just proved, of course, the general property of linear operators
that they are their own Jacobians.} The Jacobian of the inverse
transformation $\mathbf{r} = \mathbf{r}(\mathbf{c})$ is given by
\begin{equation}\label{eq:invjacobian} J_{\mathbf{r}} = J^{-1}_\mathbf{c} = \mathsf{J}^{-1} =
\frac{1}{2} \mathsf{J}^H .
\end{equation}
Of course, then, we have the differential relationships
\begin{equation}\label{eq:dorder}
d\mathbf{c} = \frac{\partial \mathbf{c}}{\partial \mathbf{r}}\, d
\mathbf{r} = J_{\mathbf{c}}\, d\mathbf{r} = \mathsf{J} d\mathbf{r}
\quad \text{and} \quad d\mathbf{r} = \frac{\partial
\mathbf{r}}{\partial \mathbf{c}}\,  d \mathbf{c} = J_{\mathbf{r}}\,
d\mathbf{c} = \frac{1}{2} \mathsf{J}^H d \mathbf{c}
\end{equation}
which correspond to the first-order relationships\footnote{For a
general, \emph{nonlinear,} coordinate transformation these finite-difference (non-infinitesimal)
first-order relationships would be \emph{approximate.} However,
because the coordinate transformation considered here happens to be
\emph{linear,} the relationships are \emph{exact.}}
\begin{equation}\label{eq:deltaorder} \boxed{\ \ \text{\sf 1st-Order Relationships:} \quad  \Delta \mathbf{c}
= J_{\mathbf{c}}\,  \Delta \mathbf{r} = \mathsf{J} \Delta \mathbf{r}
\quad \text{and} \quad \Delta \mathbf{r} = J_{\mathbf{r}}\, \Delta
\mathbf{c} = \frac{1}{2} \mathsf{J}^H \Delta \mathbf{c}\ \ }
\end{equation}
where the Jacobian $\mathsf{J}$ is given by (\ref{eq:bigjinv}) and
\begin{equation}\label{eq:deltaequals}
\Delta \mathbf{c} = \begin{pmatrix} \Delta \mathbf{z} \\ \Delta
\mathbf{\bar z} \end{pmatrix} \quad \text{and} \quad \Delta
\mathbf{r} = \begin{pmatrix} \Delta \mathbf{x} \\ \Delta \mathbf{y}
\end{pmatrix}
\end{equation}

\paragraph{The Cogradient with respect to the Real Conjugate Coordinates
Vector $\mathbf{c}.$} The reader might well wonder why we didn't
just point out that (\ref{eq:dorder}) and (\ref{eq:deltaorder}) are
merely simple consequences of the linear nature of the coordinate
transformations (\ref{eq:mapforward}) and (\ref{eq:mapback}), and
thereby skip the intermediate steps given above. The point is
that once we have identified the
Jacobian of a coordinate transformation over a real manifold, we can
readily transform between different coordinate representations of
\emph{all} vector-like (contravariant) objects, such as the gradient
of a functional, and between \emph{all} covector-like (covariant)
objects, \emph{such as the cogradient of a functional,} over that
manifold. Indeed, as a consequence of this fact we immediately have
the important cogradient operator transformations
\begin{equation}\label{eq:cogradtransfs}
\boxed{\ \ \text{\sf Cogradient Transf's:} \quad \frac{\partial
(\cdot)}{\partial \mathbf{c}} = \frac{\partial (\cdot)}{\partial
\mathbf{r}}\,  J_r = \frac{1}{2} \frac{\partial (\cdot)}{\partial
\mathbf{r}}\,  \mathsf{J}^H \quad \text{and} \quad \frac{\partial
(\cdot)}{\partial \mathbf{r}} = \frac{\partial (\cdot)}{\partial
\mathbf{c}}\,  J_c = \frac{\partial (\cdot)}{\partial \mathbf{c}}\,
\mathsf{J}\ \ }
\end{equation}
with the Jacobian $\mathsf{J}$ given by (\ref{eq:bigj}) and
$J_{\mathbf{r}} = J_{\mathbf{c}}^{-1}$.

\medskip
Equation (\ref{eq:cogradtransfs}) is \emph{very important as it
allows us to easily, yet rigorously, define the cogradient taken
with respect to $\mathbf{c}$ as a true (nonformal) differential
operator provided that we view $\mathbf{c}$ as an element of the
real coordinate representation space $\mathcal{C}$.} The cogradient
$\frac{\partial (\cdot)}{\partial \mathbf{c}}$ is  well-defined in
terms of the cogradient $\frac{\partial (\cdot)}{\partial
\mathbf{r}}$ and the ``pullback'' transformation,
\[ \frac{\partial
(\cdot)}{\partial \mathbf{c}} =  \frac{1}{2} \frac{\partial
(\cdot)}{\partial \mathbf{r}}\,  \mathsf{J}^H . \] This shows that
$\frac{\partial (\cdot)}{\partial \mathbf{c}}$, which was originally
defined in terms of the cogradient and conjugate cogradients taken
with respect to $\mathbf{z}$ (the \emph{$\mathbf{c}$-complex
interpretation} of $\frac{\partial (\cdot)}{\partial \mathbf{c}}$),
can be treated as \emph{a real differential operator} with respect
to the ``real'' vector $\mathbf{c}$ (the \emph{$\mathbf{c}$-real
interpretation} of $\frac{\partial (\cdot)}{\partial
\mathbf{c}}$).\footnote{Thus we can directly differentiate an
expression like $\mathbf{c}^T \Omega \mathbf{c}$ with respect to
$\mathbf{c}$ using the standard identities of real vector calculus.
(The fact that these identities hold for the $\mathbf{r}$ calculus
and be used to prove their validity for the $\mathbf{c}$-real
calculus.) More problematic is an expression like $\mathbf{c}^H
\Omega \mathbf{c}$. It is not appropriate to take the complex
derivative of this expression with respect to the complex vector
$\mathbf{c}$ because $\mathbf{c}$, as an element of $\mathbb{C}^n$
is subject to constraints amongst its components. Instead (see immediately below) one can
use the identity $\mathbf{\bar c} = \mathbf{\tilde c} = S
\mathbf{c}$ to obtain $\mathbf{c}^H \Omega \mathbf{c} = \mathbf{c}^T
S \Omega \mathbf{c}$ which can then be differentiated with respect
to $\mathbf{c}$. Of course, this latter approach can fail if $\mathbf{c}^T
S \Omega \mathbf{c}$ cannot be interpreted in some appropriate sense in the field of real
numbers. Note that real versus complex differentiation of
$\mathbf{c}^H \Omega \mathbf{c}$ with respect to $\mathbf{c}$ would
differ by a factor of 2.}

\paragraph{Complex Conjugation.} It is easily determined that the
operation of complex conjugation, $\mathbf{z} \rightarrow
\mathbf{\bar{z}}$, is a nonlinear mapping on $\mathcal{Z} =
\mathbb{C}^n$. Consider a \emph{general} element ${\text{\boldmath ${\zeta}$}} \in
\mathbb{C}^{2n}$ written as
\[ {\text{\boldmath ${\zeta}$}} = \binom{{\text{\boldmath
${\zeta}$}}_{\text{\scriptsize top}}}{{\text{\boldmath
${\zeta}$}}_{\text{\scriptsize bottom}}}\in \mathbb{C}^{2n} =
\mathbb{C}^n \times \mathbb{C}^n \quad \text{with} \quad
{\text{\boldmath ${\zeta}$}}_{\text{\scriptsize top}} \in
\mathbb{C}^n \quad \text{and} \quad {\text{\boldmath
${\zeta}$}}_{\text{\scriptsize bottom}} \in \mathbb{C}^n \, .
\]
Of course the operation of complex conjugation on $\mathbb{C}^{2n}$,
$\text{\boldmath ${\zeta}$} \rightarrow \text{\boldmath
$\bar{\zeta}$}$, is, in general, a nonlinear mapping.

\medskip
Now consider the \emph{linear} operation of swapping the top and
bottom elements of $\text{\boldmath ${\zeta}$}$, $\text{\boldmath
${\zeta}$} \rightarrow \text{\boldmath $\tilde{\zeta}$}$, defined as
\[ {\text{\boldmath ${\zeta}$}} = \binom{{\text{\boldmath
${\zeta}$}}_{\text{\scriptsize top}}}{{\text{\boldmath
${\zeta}$}}_{\text{\scriptsize bottom}}} \rightarrow
{\text{\boldmath $\bar{\zeta}$}} = \binom{{\text{\boldmath
${\zeta}$}}_{\text{\scriptsize bottom}}}{{\text{\boldmath
${\zeta}$}}_{\text{\scriptsize top}}} = \begin{pmatrix} 0 & I \\ I &
0 \end{pmatrix} \binom{{\text{\boldmath
${\zeta}$}}_{\text{\scriptsize top}}}{{\text{\boldmath
${\zeta}$}}_{\text{\scriptsize bottom}}} = S  \text{\boldmath
${\zeta}$} \] where \[ S \triangleq \begin{pmatrix} 0 & I \\ I & 0
\end{pmatrix} \] is the \emph{swap operator} on $\mathbb{C}^{2n}$
which obeys the properties \[ S = S^T = S^{-1} \, , \] showing that
$S$ is symmetric and its own inverse, $S^2 = I$. Note that, in
general, swapping is \emph{not} equal to complex conjugation,
$\text{\boldmath $\tilde{\zeta}$} \ne \text{\boldmath
$\bar{\zeta}$}$.

\medskip
The swap operator $S$ will be used extensively throughout the
remainder of this note, so it is important to become comfortable
with its use and manipulation. The swap operator is a \emph{block
permutation matrix} which permutes (swaps)\footnote{``Permutation''
is just a fancy term for ``swapping.''} blocks of rows or blocks of
columns depending on whether $S$ premultiplies or postmultiplies a
matrix. Specifically, let a $2n \times 2n$ matrix $A$ be block
partitioned as \[ A = \begin{pmatrix} A_{11} & A_{12} \\ A_{21} &
A_{22}
\end{pmatrix} . \]
Then premultiplication by $S$ results in a block swap of the top $n$
rows \emph{en masse} with the bottom $n$ rows,\footnote{Matrix
premultiplication of $A$ by \emph{any} matrix always yields a row
operation.}
\[ S A = \begin{pmatrix}  A_{21} &  A_{22} \\ A_{11}  &  A_{12}
\end{pmatrix}. \]
Alternatively, postmultiplication by $S$ results in a block swap of
the first $n$ columns  with the last $n$ columns,\footnote{Matrix
postmultiplication  of $A$ by any matrix always yields a column
operation. The fact that pre- and postmultiplication yield different
actions on $A$ is an interesting and illuminating way to interpret
the fact that matrix multiplication is noncommutative, $MA \ne AM$.}
\[ A S = \begin{pmatrix}  A_{12} &  A_{11} \\ A_{22}  &  A_{21}
\end{pmatrix}. \]
It is also useful to note the result of a ``sandwiching'' by $S$,
\[ S A S = A = \begin{pmatrix} A_{22} & A_{21} \\ A_{12} & A_{11}
\end{pmatrix} . \]
Because $S$ permutes $n$ rows (or columns), it is a product of $n$
elementary permutation matrices, each of which is known to have a
determinant which evaluates to $-1$. As an easy consequence of this,
we have
\[ \det S = (-1)^n . \]
Other important properties of the swap operator $S$ will be
developed as we proceed.

\medskip
Now note that the subset $\mathcal{C} \in \mathbb{C}^{2n}$ contains
precisely those elements of $\mathbb{C}^{2n}$ for which the
operations of swapping and complex conjugation coincide,
\[ \mathcal{C} = \left\{ \text{\boldmath ${\zeta}$} \in \mathbb{C}^{2n} \, \left| \, \, \text{\boldmath
$\bar{\zeta}$} = \text{\boldmath $\tilde{\zeta}$} \right\} \right.
\subset \mathbb{C}^{2n} \, , \] and thus it is true by construction
that $\mathbf{c} \in \mathcal{C}$ obeys $\mathbf{\bar c} =
\mathbf{\tilde c}$, \emph{even though swapping and complex
conjugation are different operations} on $\mathbb{C}^{2n}$. Now
although $\mathcal{C}$ is not a subspace of the \emph{complex} vector space
$\mathbb{C}^{2n}$, it is a \emph{real} vector space in its own
right. We see that the \emph{linear} operation of component swapping
on the $\mathcal{C}$-space coordinate representation of $\mathcal{Z}
= \mathbb{C}^n$ is exactly equivalent to the \emph{nonlinear}
operation of complex conjugation on $\mathcal{Z}$. It is important
to note that complex conjugation and coordinate swapping represent
different operations on a vector $\mathbf{c}$ when $\mathbf{c}$ is
viewed as an element of $\mathbb{C}^{2n}$.\footnote{As mentioned
earlier, $\mathbf{c}$, in a sense, does ``double duty'' as a
representation for $\mathbf z$; once as a (true coordinate)
representation of $\mathbf z$ in the real vector space
$\mathcal{C}$, and alternatively as a ``representation'' of $z$ in
the ``doubled up'' complex space $\mathbb{C}^{2n} = \mathbb{C}^n
\times \mathbb{C}^n$. In the development given below, we will switch
between these two perspectives of $\mathbf{c}$.}

\medskip
We can view the linear swap mapping $S: \mathcal{C} \rightarrow
\mathcal{C}$ as a coordinate transformation (a coordinate
``reparameterization''), $\mathbf{\bar c} = \mathbf{\tilde c} = S
\mathbf{c}$, on $\mathcal{C}$. Because $S$ is linear, the Jacobian
of this transformation is just $S$ itself. Thus from the cogradient
transformation property we obtain the useful identity
\begin{equation}\label{eq:tildecogradid}
\frac{\partial (\cdot)}{\partial \mathbf{\bar c}} \, S =
\frac{\partial (\cdot)}{\partial \mathbf{\tilde c}} \, S =
\frac{\partial (\cdot)}{\partial \mathbf{c}}
\end{equation}
It is also straightforward to show that
\begin{equation}\label{eq:ij}
I = \frac{1}{2} \mathsf{J}^T S \mathsf{J}
\end{equation}
for $\mathsf{J}$ given by (\ref{eq:bigj})

\medskip
Let us now turn to the alternative coordinate representation given
by vectors $\mathbf{r}$ in the space $\mathcal{R} =
\mathbb{R}^{2n}$. Specifically, consider the $\mathcal{R}$
coordinate vector $\mathbf{r}$ corresponding to the change of
coordinates $\mathbf{r} = \frac{1}{2} \mathsf{J}^H \mathbf{c}$.
Since the vector $\mathbf{r}$ is real, it is its own complex
conjugate, $\mathbf{\bar r} = \mathbf{r}$.\footnote{Note that our
theoretical developments are consistent with this requirement, as \[
\mathbf{\bar r} = \frac{1}{2}\overline{(\mathsf{J}^H \mathbf{c})}=
\frac{1}{2} \mathsf{J}^T \mathbf{\bar c} = \frac{1}{2} \mathsf{J}^T
\mathbf{\tilde c} = \frac{1}{2} \mathsf{J}^T S \mathbf{c} =
\frac{1}{2} \mathsf{J}^T S \mathsf{J} \mathbf{r} = I \mathbf{r} =
\mathbf{r} \, . \]} Complex conjugation of $\mathbf{z}$
 is the \emph{nonlinear mapping} in $\mathbb{C}^n$
\[ \mathbf{z} = \mathbf{x} + j \, \mathbf{y} \rightarrow
\mathbf{\bar z} = \mathbf{x} + j \, (- \mathbf{y}) \, , \] and
corresponds in the representation space $\mathcal{R}$ to the
\emph{linear mapping}\footnote{We refer to $\mathbf{\check{r}}$ as
``$\mathbf{r}$-check.''}
\[ \mathbf{r} = \binom{\mathbf{x}}{\mathbf{y}} \rightarrow
\mathbf{\check{r}} \triangleq \binom{\mathbf{x}}{- \mathbf{y}} =
\begin{pmatrix} I & 0 \\ 0 & -I \end{pmatrix}
\binom{\mathbf{x}}{\mathbf{y}} = C \mathbf{r}
\]
where $C$ is the conjugation matrix
\begin{equation}
C \triangleq \begin{pmatrix} I & 0 \\ 0 & -I \end{pmatrix}.
\end{equation}
Note that
\[ C = C^T = C^{-1} \, , \]
i.e., that $C$ is symmetric, $C = C^T$, and its own inverse, $C^2 =
I$. It is straightforward to show that
\begin{equation}\label{eq:cj}
C = \frac{1}{2} \mathsf{J}^H S \mathsf{J}
\end{equation}
which can be compared to (\ref{eq:ij}). Finally, it is
straightforward to show that
\begin{equation}
\mathbf{c} = \mathsf{J} \mathbf{r} \Leftrightarrow \mathbf{\bar c} =
\mathbf{\tilde c} = \mathsf{J} \mathbf{\check{r}} \, .
\end{equation}

\medskip
To summarize, we can represent the complex vector $\mathbf{z}$ by
either $\mathbf{c}$ or $\mathbf{r}$, where $\mathbf{c}$ has two
interpretations (as a complex vector, ``$\mathbf c$-complex'', in
$\mathbb{C}^{2n}$, or as an element, ``$\mathbf c$-real'',  of the
real vector space $\mathcal{C} \approx \mathbb{R}^{2n}$), and we can
represent the complex conjugate $\mathbf{\bar z}$ by $\mathbf{\bar
c}$, $\mathbf{\tilde c}$, or $\mathbf{\check r}$. And complex
conjugation, which is a nonlinear operation in $\mathbb{C}^n$,
corresponds to linear operators in the $2n$-dimensional isomorphic
real vector spaces $\mathcal{C}$ and $\mathcal{R}$.

\subsection{Low Order Series Expansions of a
Real-Valued Scalar Function.} By noting that a real-valued scalar
function of complex variables can be viewed as a function of either
$\mathbf{r}$ or $\mathbf{c}$-real or  $\mathbf{c}$-complex or
$\mathbf{z}$,
\[ f(\mathbf{r}) = f(\mathbf{c}) = f(\mathbf{z}) \, ,
\] it is evident that one should be able to represent $f$ as a power
series in any of these representations.  Following the line of
attack pursued by van den Bos in \cite{VDB94a}, by exploiting the relationships
(\ref{eq:deltaorder}) and (\ref{eq:cogradtransfs}) we will readily
show the equivalence up to second order in a power series expansion
of $f$.

\medskip
Up to second order, the multivariate power series expansion of the
real-valued function $f$ viewed as an analytic function of vector
$\mathbf{r} \in \mathcal{R}$ is given as, {\small
\begin{equation}\label{eq:real2ndorder} \boxed{\ \ \text{2nd-Order
Expansion in $\mathbf{r}$:} \ \ f(\mathbf{r} + \Delta \mathbf{r}) =
f(\mathbf{r}) + \frac{\partial f(\mathbf{r})}{\partial \mathbf{r}}
\, \Delta \mathbf{r} + \frac{1}{2} \Delta \mathbf{r}^T \,
\mathcal{H}_{\mathbf{r}\mathbf{r}}(\mathbf{r}) \, \Delta \mathbf{r}
+ \text{h.o.t.}\ \ }
\end{equation}}
where\footnote{When no confusion can arise, one usually drops the
subscripts on the Hessian and uses the simpler notation
$\mathcal{H}(\text{\boldmath $\rho$}) = \mathcal{H}_{\mathbf{r}
\mathbf{r}}(\text{\boldmath $\rho$})$.  Note that the Hessian is the matrix of second partial
derivatives of a \emph{real-valued} scalar function.}
\begin{equation}\label{eq:hessian}
 \mathcal{H}_{\mathbf{r} \mathbf{r}}(\text{\boldmath $\rho$}) \triangleq
\frac{\partial}{\partial \mathbf{r}} \left( \frac{\partial
f(\text{\boldmath $\rho$})}{\partial \mathbf{r}} \right)^T \quad
\text{for} \quad \text{\boldmath $\rho$}, \mathbf{r} \in \mathcal{R}
\end{equation}
is the real $\mathbf{r}$-Hessian matrix of second partial
derivatives of the real-valued function $f(\mathbf{r})$ with respect
to the components of $\mathbf{r}$. It is well known that a real
Hessian is symmetric, \[ \mathcal{H}_{\mathbf{r} \mathbf{r}} =
\mathcal{H}_{\mathbf{r} \mathbf{r}}^T \, . \] However, there is no
general guarantee that the Hessian will be a positive definite or
positive semidefinite matrix.

\medskip
It is assumed that the terms $f(\mathbf{r})$ and $f(\mathbf{r} +
\Delta \mathbf{r})$ be readily expressed in terms of $\mathbf{c}$
and $\mathbf{c} + \Delta \mathbf{c}$ or $\mathbf{z}$ and $\mathbf{z}
+ \Delta \mathbf{z}$. Our goal is to determine the proper expression
of the linear and quadratic terms of  (\ref{eq:real2ndorder}) in
terms of $\mathbf{c}$ and $\Delta \mathbf{c}$ or $\mathbf{z}$ and
$\Delta \mathbf{z}$.

\paragraph{\boldmath Scalar Products and Quadratic Forms on the Real Vector
Space $\mathcal{C}$.} Consider two vectors $\mathbf{c} =
\text{col}(\mathbf{z},\mathbf{\bar z}) \in \mathcal{C}$ and
$\mathbf{s} = \text{col}(\text{\boldmath ${\xi}$,$\bar \xi$}) \in
\mathcal{C}$. The scalar product for any two such vectors in
$\mathcal{C}$-real (i.e., in the \emph{real vector space}
$\mathcal{C}\approx \mathbb{R}^{2n}$) is defined by
\[ \left< \mathbf{c},\mathbf{s} \right> \triangleq \mathbf{c}^T S \,
\mathbf{s} = \mathbf{\bar c}^T \mathbf{s} = \mathbf{c}^H \mathbf{s}
= \mathbf{z}^H \text{\boldmath ${\xi}$} + \mathbf{\bar z}^H
\text{\boldmath ${\bar \xi}$} =
\mathbf{z}^H \text{\boldmath ${\xi}$} + \overline{\mathbf{z}^H
\text{\boldmath ${\xi}$}}  = 2\, \text{Re} \, \mathbf{z}^H
\text{\boldmath ${\xi}$} \, . \] The row vector  $\mathbf{c}^T S =
\mathbf{c}^H$ is a linear functional which maps the elements of
$\mathcal{C}$-real into the real numbers. The set of all such linear
functionals is a vector space itself and is known as the \emph{dual
space,} $\mathcal{C}^*$, of $\mathcal C$ \cite{Lu69,NS82}. The
elements of $\mathcal{C}^*$ are known as \emph{dual vectors} or
\emph{covectors}, and the terms ``dual vector'', ``covector'', and
``linear functional'' should all be taken to be synonymous. Given a
vector $\mathbf{c} \in \mathcal{C}$, there is a natural one-to-one
mapping between $\mathbf{c}$ and a corresponding dual vector,
$\mathbf{c}^*$ in $\mathcal{C}^*$ defined by\footnote{Warning! Do
not confuse the dual vector (linear functional) $\mathbf{c}^*$ with
an adjoint operator, which is often also denoted using the ``star''
notation.}
\[ \mathbf{c}^* \triangleq \mathbf{c}^T  S = \mathbf{c}^H. \]

Henceforth it is understood that scalar-product expressions like
\[ \mathbf{a}^H \mathbf{s} \quad \text{or} \quad \mathbf{c}^H \mathbf{b} \]
where $\mathbf{s} \in \mathcal{C}$ and $\mathbf{c} \in \mathcal{C}$
are known to be elements of $\mathcal{C}$ \emph{are only
{meaningful} if $\mathbf{a}$ and $\mathbf{b}$ are also elements of
$\mathcal{C}$}. Thus,  \emph{it must be the case} that \emph{both}
vectors in a scalar product must belong to $\mathcal C$ if it is the
case that one of them does, otherwise we view the resulting
numerical value as nonsensical.

\medskip
Thus, for a real-valued function of up to quadratic order in a
vector $\mathbf{c} \in \mathcal{C}$,
\begin{equation}\label{eq:gen2fun}
f(\mathbf{c}) = a + \mathbf{b}^H \mathbf{c} + \frac{1}{2}
\mathbf{c}^H  M \mathbf{c} = a + \mathbf{b}^H \mathbf{c} +
\frac{1}{2} \mathbf{c}^H  \mathbf{s}, \quad \mathbf{s} = M
\mathbf{c},
\end{equation}
to be well-posed, it \emph{must} be the case that $a \in \mathbb{R}$,
$\mathbf{b} \in \mathcal{C}$,\footnote{I.e., that $\mathbf{b}^H$ be
a bona fide linear functional on $\mathcal{C}$, $\mathbf{b}^H =
\mathbf{b}^* \in \mathcal{C}^*$.} and $\mathbf{s} = M \mathbf{c} \in
\mathcal{C}$.\footnote{I.e., because $\mathbf{c}^H = \mathbf{c}^*
\in \mathcal{C}^*$, is a linear functional on $\mathcal{C}$, it must
have a legitimate object $\mathbf{s}$ to operate on, namely an
element $\mathbf{s} = M \mathbf{c} \in \mathcal{C}$.} Thus, as we
proceed to derive various first and second order functions of the
form (\ref{eq:gen2fun}), \emph{we will need to check for these conditions.}
If the conditions are met, we will say that vector $\textbf{b}$ and the operator $M$; the terms $\mathbf{b}^H \mathbf{c}$ and $\mathbf{c}^H  M \mathbf{c}$; and the entire
quadratic form itself, are \emph{\underline{admissible}} \, (or \emph{{meaningful}}).

\medskip
Thus $\textbf{b}$ is admissible if and only if $\textbf{b}\in \mathcal{C}$, and $M$ is admissible if and only if $M$ is a linear mapping from $\mathcal{C}$ to $\mathcal{C}$, $M \in \mathcal{L}(\mathcal{C},\mathcal{C})$.

\medskip
To test whether a vector $\mathbf{b} \in \mathbb{C}^{2n}$ belongs to
$\mathcal{C}$ is straightforward:
\begin{equation}\label{eq:vecttest}
\mathbf{b} \in \mathcal{C} \Leftrightarrow \mathbf{\bar b} = S
\mathbf{b} .
\end{equation}

\medskip
It is somewhat more work to develop a test to determine if a  matrix
$M \in \mathbb{C}^{2n \times 2n}$ has the property that it is a
linear mapping from $\mathcal{C}$ to $\mathcal{C}$,
\[ M \in \mathcal{L}(\mathcal{C},\mathcal{C}) = \left\{ M \left| \,
\,  M \mathbf{c} \in \mathcal{C}, \ \forall \mathbf{c}\in
\mathcal{C}  \ \ \text{and} \ \  M \ \text{is linear} \, \right\}
\right. \subset \mathcal{L}(\mathbb{C}^{2n},\mathbb{C}^{2n}) =
\mathbb{C}^{2n \times 2n} .
\] Note that the fact that $\mathcal{L}(\mathcal{C},\mathcal{C}) \subset
\mathcal{L}(\mathbb{C}^{2n},\mathbb{C}^{2n})$ is just the statement
that any matrix which maps from $\mathcal{C} \subset
\mathbb{C}^{2n}$ to $\mathcal{C} \subset \mathbb{C}^{2n}$ is also
obviously a linear mapping from $\mathbb{C}^{2n}$ to
$\mathbb{C}^{2n}$. However, \emph{this is just a subset statement};
it is \emph{not} a subspace statement. This is because
$\mathcal{L}(\mathcal{C},\mathcal{C})$ is a \emph{real} vector space
of linear operators,\footnote{I.e., a vector space over the field of
real numbers.} while $\mathcal{L}(\mathbb{C}^{2n},\mathbb{C}^{2n})$
is a complex vector space of linear operators.\footnote{I.e., a
vector space over the field of complex numbers.} Because they are
vector spaces over \emph{different fields}, they cannot have a
vector-subspace/vector-parent-space relationship to each other.

\medskip
To determine necessary and sufficient conditions for a matrix $M \in
\mathbb{C}^{2n \times 2n}$ to be an element of
$\mathcal{L}(\mathcal{C},\mathcal{C})$ suppose that the vector
$\mathbf{c} = \text{col}(\mathbf{z},\mathbf{\bar z}) \in
\mathcal{C}$ always maps to a vector $\mathbf{s} =
\text{col}(\text{\boldmath $\xi$},\text{\boldmath $\bar \xi$}) \in
\mathcal{C}$ under the action of $M$, $\mathbf{s} = M \mathbf{c}$.
Expressed in block matrix form, this relationship is
\[ \binom{\text{\boldmath $\xi$}}{\text{\boldmath $\bar \xi$}} =
\begin{pmatrix} M_{11} & M_{12} \\ M_{21} & M_{22}
\end{pmatrix}
\binom{\mathbf{z}}{\mathbf{\bar z}} .\] The first block row of this
matrix equation yields the conditions \[ \text{\boldmath $\xi$} =
M_{11} \mathbf{z} + M_{12} \mathbf{\bar z} \] while the complex
conjugate of the second block row yields
\[ \text{\boldmath $\xi$} = \bar{M}_{22}\mathbf{z} +
\bar{M}_{21}\mathbf{\bar z} \] and subtracting these two sets of
equations results in the following condition on the block elements
of $M$,
\[ (M_{11} - \bar{M}_{22}) \mathbf{z} + (M_{12}- \bar{M}_{21})
\mathbf{\bar z} = 0 \, . \] With $\mathbf{z} = \mathbf{x} + j \,
\mathbf{y}$, this splits into the two sets of conditions,
\[ [(M_{11} - \bar{M}_{22}) +  (M_{12}- \bar{M}_{21})]
\mathbf{x} = 0 \] and
\[ [(M_{11} - \bar{M}_{22}) -  (M_{12}- \bar{M}_{21})]
\mathbf{y} = 0 . \] Since these equations must hold for any
$\mathbf{x}$ and $\mathbf{y}$, they are equivalent to
\[ (M_{11} - \bar{M}_{22}) +  (M_{12}- \bar{M}_{21})
 = 0 \] and
\[ (M_{11} - \bar{M}_{22}) -  (M_{12}- \bar{M}_{21})
 = 0 . \]
 Finally, adding and subtracting these two equations yields the
 necessary and sufficient conditions for $M$ to admissible (i.e., to be a mapping from
 $\mathcal{C}$ to $\mathcal{C}$),
 \begin{equation}\label{eq:ccmappingcond1}
 M = \begin{pmatrix} M_{11} & M_{12} \\ M_{21} & M_{22}
\end{pmatrix} \in \mathbb{C}^{2n\times 2n} \ \ \text{is an element
of $\mathcal{L}(\mathcal{C},\mathcal{C})$ iff} \ \ M_{11} =
\bar{M}_{22} \ \ \text{and} \ \ M_{12} = \bar{M}_{21} \, .
\end{equation}
This necessary and sufficient admissibility condition is more conveniently
expressed in the following equivalent form,
\begin{equation}\label{eq:ccmappingcond2}
M \in \mathcal{L}(\mathcal{C},\mathcal{C}) \Leftrightarrow M = S
\bar{M} S \Leftrightarrow \bar{M} = S {M} S
\end{equation}
which  is straightforward to verify.

\medskip
Given an arbitrary matrix $M \in \mathbb{C}^{2n \times 2n}$, we can
define a natural mapping of $M$ into
$\mathcal{L}(\mathcal{C},\mathcal{C}) \subset \mathbb{C}^{2n \times
2n}$ by
\begin{equation}\label{eq:Cproj}
\mathbf{P}(M) \triangleq \frac{M + S \bar{M} S}{2} \in
\mathcal{L}(\mathcal{C},\mathcal{C})\, ,
\end{equation}
in which case the admissibility condition (\ref{eq:ccmappingcond2}) has an
equivalent restatement as
\begin{equation}\label{eq:ccmappingcond3}
M \in \mathcal{L}(\mathcal{C},\mathcal{C}) \Leftrightarrow
\mathbf{P}(M) = M \, .
\end{equation}
 It is
straightforward to demonstrate that
\begin{equation}\label{eq:Pprops}
  \forall M \in \mathbb{C}^{2n
\times 2n} ,  \ \ \ \mathbf{P}(\mathbf{P}(M)) =
\mathbf{P}(M) \, . \end{equation} I.e., $\mathbf P$ is an idempotent
mapping of $\mathbb{C}^{2n\times 2n}$ onto
$\mathcal{L}(\mathcal{C},\mathcal{C})$, \, $\mathbf{P}^2 = \mathbf{P}$.
However,  as things currently stand $\mathbf P$ is \emph{not} a
linear operator (the action of complex conjugation precludes this)
nor a projection operator in the conventional sense of projecting
onto a lower dimensional \emph{subspace} as its range space is
\emph{not} a subspace of its domain space. (However, with some additional work, one can reasonably
interpret $\mathbf{P}$ as a projector of the \emph{manifold}
$\mathbb{C}^{2n}$ onto the \emph{submanifold} $\mathcal{C} \subset
\mathbb{C}^{2n}$ in some sense.\footnote{With $\mathbb{C}^{2n\times 2n} \approx
\mathbb{R}^{4n \times 4n} \approx \mathbb{R}^{16n^2}$ and
$\mathcal{L}(\mathcal{C},\mathcal{C}) \approx
\mathcal{L}(\mathbb{R}^{2n},\mathbb{R}^{2n}) \approx \mathbb{R}^{2n
\times 2n} \approx \mathbb{R}^{4n^2}$, it is reasonable to view
$\mathbf{P}$ as a linear projection operator from the \emph{real vector
space} $\mathbb{R}^{16 n^2}$ onto the \emph{real vector subspace}
$\mathbb{R}^{4 n^2}$ of $\mathbb{R}^{4 n}$. This allows us to
interpret $\mathbf{P}$ as a projection operator from the
\emph{manifold} $\mathbb{C}^{2n}$ onto the \emph{submanifold}
$\mathcal{C} \subset \mathbb{C}^{2n}$. Once we know that
$\mathbf{P}$ is a \emph{linear} mapping from $\mathbb{C}^{2n}$ into
$\mathbb{C}^{2n}$, we can then compute its adjoint operator,
$\mathbf{P}^*$, and then test to see if its self-adjoint. If it is,
then the projection operator $\mathbf{P}$ is, in fact, an orthogonal
projection operator.})

\medskip
A final important fact is that if $M \in \mathbb{C}^{2n \times 2n}$
 is invertible, then $M \in \mathcal{L}(\mathcal{C},\mathcal{C})$ if
 and only if
$M^{-1} \in \mathcal{L}(\mathcal{C},\mathcal{C})$, which we state
equivalently as
\begin{equation}\label{eq:invtoo}
\text{Let $M$ be invertible, then $\mathbf{P}(M) = M$ iff
$\mathbf{P}(M^{-1}) = M^{-1}$.}
\end{equation}
I.e., if an invertible matrix $M$ is admissible, then $M^{-1}$ is
admissible. The proof is straightforward:
\begin{eqnarray*}
M  =  S \bar{M} S  & \text{and} & M \ \text{invertible} \\
\Leftrightarrow M^{-1} & = & \left( S \bar{M} S \right)^{-1} \\
& = & S (\bar{M})^{-1} S \\
& = & S \overline{M^{-1}} S \, .
\end{eqnarray*}

\paragraph{First Order Expansions.}
Up to first order, the power series expansion of the real-valued
function $f$ viewed as a function of $\mathbf{r} \in \mathcal{R}$ is
\begin{equation}\label{eq:real1ndorder}
\boxed{\ \ \text{First-Order Expansion in $\mathbf{r}$:} \qquad
f(\mathbf{r} + \Delta \mathbf{r}) = f(\mathbf{r}) + \frac{\partial
f(\mathbf{r})}{\partial \mathbf{r}} \, \Delta \mathbf{r} +
\text{h.o.t.}\ \ }
\end{equation}

\medskip
Focussing our attention first on the linear term $\frac{\partial
f(\mathbf{r})}{\partial \mathbf{r}} \, \Delta \mathbf{r}$, and using
the $\mathbf{c}$-real vector space interpretation of $\mathbf{c}$,
namely that $\mathbf{c} \in \mathcal{C}$ where, as discussed above,
$\mathcal{C}$ is a $2n$-dimensional coordinate space isomorphic to
$\mathbb{R}^{2n}$, we have
\begin{eqnarray*}
\frac{\partial f}{\partial \mathbf{r}} \, \Delta \mathbf{r} & = &
\frac{\partial f}{\partial \mathbf{r}} \, J^{-1}_{\mathbf{c}} \,
\Delta \mathbf{c} \qquad \ \text{(from equation
(\ref{eq:deltaorder}))} \\
& = & \frac{\partial f}{\partial \mathbf{c}} \, \Delta \mathbf{c}
\qquad \qquad \text{(from equation (\ref{eq:cogradtransfs}))}
\end{eqnarray*}
which yields the first order expansion of $f$ in terms of the
parameterization in $\mathbf{c}$,
\begin{equation}\label{eq:lininc}
\boxed{\ \ \text{First-Order Expansion in $\mathbf{c}$:} \quad
f(\mathbf{c} + \Delta \mathbf{c}) = f(\mathbf{c}) + \frac{\partial
f(\mathbf{c})}{\partial \mathbf{c}} \, \Delta \mathbf{c} +
\text{h.o.t.}\ \ }
\end{equation}
 Note that $\frac{\partial f(\mathbf{c})}{\partial
\mathbf{c}}\, \Delta \mathbf{c}$ is real valued. Furthermore, as a
consequence of the fact that with $f(\mathbf{c})$ real-valued we
have
\[ \overline{\left( \frac{\partial f(\mathbf{c})}{\partial
\mathbf{c}} \right)^H} = \left( \frac{\partial
f(\mathbf{c})}{\partial \mathbf{\bar c}} \right)^H = S
\left(\frac{\partial f(\mathbf{c})}{\partial \mathbf{c}}\right)^H ,
\] the quantity $\left(\frac{\partial f(\mathbf{c})}{\partial \mathbf{c}}\right)^H$ satisfies the necessary and sufficient condition given in
(\ref{eq:vecttest}) that
\[ \left(\frac{\partial f(\mathbf{c})}{\partial \mathbf{c}} \right)^H \in \mathcal{C} \, . \]
Thus $\frac{\partial f(\mathbf{c})}{\partial \mathbf{c}} \in
\mathcal{C}^*$ and the term $\frac{\partial f(\mathbf{c})}{\partial
\mathbf{c}}\, \Delta \mathbf{c}$ is admissible in the sense defined
earlier. Note that an equivalent condition for the term
$\frac{\partial f(\mathbf{c})}{\partial \mathbf{c}}\, \Delta
\mathbf{c}$ to be admissible is that
\[  S \left(\frac{\partial f(\mathbf{c})}{\partial \mathbf{c}} \right)^T \in \mathcal{C} , \]
which is true if and only if
\[   \left(\frac{\partial f(\mathbf{c})}{\partial \mathbf{c}} \right)^T \in \mathcal{C} . \]
This shows a simple inspection of $\frac{\partial
f(\mathbf{c})}{\partial \mathbf{c}}$ itself can be performed to test
for admissibility of the first-order term.\footnote{In this note, the
first order expansion (\ref{eq:lininc}) is doing double duty in that
it is simultaneously standing for the $\mathbf{c}$-real expansion
and the $\mathbf{c}$-complex expansion. A more careful development
would make this distinction explicit, in which case one would more
carefully explore the distinction between  $\left(\frac{\partial
f(\mathbf{c})}{\partial \mathbf{c}} \right)^T$ versus
$\left(\frac{\partial f(\mathbf{c})}{\partial \mathbf{c}} \right)^H$
in the first-order  term. Because this note has already become rather
notationally tedious, this option for greater precision has been
declined. However, greater care must therefore be made when
switching between the $\mathcal{C}$-real and $\mathcal{C}$-complex
perspectives.}

\medskip
As discussed above, to be meaningful as a true derivative, the
derivative with respect to $\mathbf{c}$ has to be interpreted as a
real derivative. This is provided by the $\mathbf{c}$-real interpretation of
(\ref{eq:lininc}). In addition, (\ref{eq:lininc}) has a
$\mathbf{c}$-complex interpretation for which the partial derivative
with respect to $\mathbf{c}$ is not well-defined as a complex
derivative as it stands, but rather only makes sense as a shorthand
notation for simultaneously taking the complex derivatives with
respect to $\mathbf{z}$ and $\mathbf{\bar z}$, \[
\frac{\partial}{\partial \mathbf{c}} = \left(
\frac{\partial}{\partial \mathbf{z}}\, , \, \frac{\partial}{\partial
\mathbf{\bar z}} \right) . \] Thus, to work in the domain of complex
derivatives, we must move to the $\mathbf{c}$-complex perspective
$\mathbf{c} = \text{col}(\mathbf{z},\mathbf{\bar z})$, and then
break $\mathbf{c}$ apart so that we can work with expressions
explicitly involving $\mathbf{z}$ and $\mathbf{\bar z}$, exploiting
the fact that the formal partial derivatives with respect to
$\mathbf{z}$ and $\mathbf{\bar{z}}$ are well defined.

\medskip
Noting that
\[ \frac{\partial}{\partial \mathbf{c}} =
  \begin{pmatrix} \frac{\partial}{\partial \mathbf{z}} & \frac{\partial}{\partial \mathbf{\bar
  z}}
  \end{pmatrix} \quad \text{and} \quad \Delta \mathbf{c} =
  \begin{pmatrix} \Delta \mathbf{z} \\ \Delta \mathbf{\bar z}
  \end{pmatrix} \]
we obtain
\begin{eqnarray*}
\frac{\partial f(\mathbf{c})}{\partial \mathbf{c}} \, \Delta
\mathbf{c} & = & \frac{\partial f}{\partial \mathbf{z}} \, \Delta
\mathbf{z} + \frac{\partial f}{\partial \mathbf{\bar z}} \, \Delta
\mathbf{\bar z} \\
& = & \frac{\partial f}{\partial \mathbf{z}} \, \Delta \mathbf{z} +
\overline{\frac{\partial f}{\partial \mathbf{z}} \, \Delta
\mathbf{z}} \qquad \text{($f$ is real-valued)} \\
& = & 2\,  \text{Re} \left\{ \frac{\partial f}{\partial \mathbf{z}}
\, \Delta \mathbf{z} \right\}
\end{eqnarray*}
which yields the first-order expansion of $f$ in terms of the
parameterization in $\mathbf{z}$,
\begin{equation}\label{eq:lininz}
\boxed{\ \ \text{First-Order Expansion in $\mathbf{z}$:} \quad
f(\mathbf{z} + \Delta \mathbf{z}) = f(\mathbf{z})+  2\,  \text{Re}
\left\{ \frac{\partial f}{\partial \mathbf{z}} \, \Delta \mathbf{z}
\right\}+ \text{h.o.t.} \ \ }
\end{equation}
This is the rederivation of (\ref{eq:deltalinexpansion}) promised
earlier. Note that (\ref{eq:lininz}) makes \emph{explicit} the
relationship which is \emph{implied} in the $\mathbf{c}$-complex
interpretation of (\ref{eq:lininc}).

\medskip
We also summarize our intermediate  results concerning the linear
term in a power series expansion using the $\mathbf{r}$,
$\mathbf{c}$ or $\mathbf{z}$ representations,
\begin{equation}\label{eq:rcz1terms}
\boxed{\ \ \text{Linear-Term Relationships:} \quad \frac{\partial
f}{\partial \mathbf{r}} \, \Delta \mathbf{r} = \frac{\partial
f}{\partial \mathbf{c}} \, \Delta \mathbf{c} = 2\, \text{Re} \left\{
\frac{\partial f}{\partial \mathbf{z}} \, \Delta \mathbf{z}
\right\}\ \ }
\end{equation}
The derivative in the first expression is a real derivative. The
derivative in the second expression is interpreted as a real
derivative (the $\mathbf{c}$-real interpretation). The derivative in
the last expression is a complex derivative; it corresponds to the
$\mathbf{c}$-complex interpretation of the second term in
(\ref{eq:rcz1terms}). Note that all of the linear terms are real
valued.

\medskip
We now have determined the first-order expansion of $f$ in terms of
$\mathbf{r}$, $\mathbf{c}$, and $\mathbf{z}$. To construct the
second-order expansion it remains to examine the second-order term
in (\ref{eq:real2ndorder}) and some of the properties of the real
Hessian matrix (\ref{eq:hessian}) which completely specifies that
term.

\paragraph{Second Order Expansions.} Note from (\ref{eq:real2ndorder}) that knowledge of the
real Hessian matrix $\mathcal{H}_{\mathbf{r} \mathbf{r}}$ completely
specifies the second order term in the real power series expansion
of $f$ with respect to $\mathbf{r}$. The goal which naturally
presents itself to us at this point is now to reexpress this
quadratic-order term in terms of $\mathbf{c}$, which we indeed
proceed to do. However, because the canonical coordinates vector
$\mathbf{c}$ \emph{has two interpretations,} one as a shorthand for
the pair $(\mathbf{z},\mathbf{\bar z})$ (the $\mathbf{c}$-complex
perspective) and the other as an element of a real vector space (the
$\mathbf{c}$-real perspective), we will rewrite the second order
term in two different forms, one (the $\mathbf{c}$-complex form)
involving \emph{the $\mathbf{c}$-complex Hessian matrix}
\begin{equation}\label{eq:chessian}
\mathcal{H}_{\mathbf{c} \mathbf{c}}^{\text{\boldmath \tiny
$\mathbb{C}$}} (\text{\boldmath $\upsilon$}) \triangleq
\frac{\partial}{\partial \mathbf{c}} \left( \frac{\partial
f(\text{\boldmath $\upsilon$)}}{\partial \mathbf{c}} \right)^H \quad
\text{for} \quad \text{\boldmath $\upsilon$}, \mathbf{c} \in
\mathcal{C} \subset \mathbb{C}^{2n}
\end{equation}
and the other (the $\mathbf{c}$-real form) involving \emph{the
$\mathbf{c}$-real Hessian matrix}
\begin{equation}\label{eq:chessianr}
\mathcal{H}_{\mathbf{c} \mathbf{c}}^{\text{\boldmath \tiny
$\mathbb{R}$}} (\text{\boldmath $\upsilon$}) \triangleq
\frac{\partial}{\partial \mathbf{c}} \left( \frac{\partial
f(\text{\boldmath $\upsilon$)}}{\partial \mathbf{c}} \right)^T \quad
\text{for} \quad \text{\boldmath $\upsilon$}, \mathbf{c} \in
\mathcal{C} \approx \mathbb{R}^{2n} .
\end{equation}
In (\ref{eq:chessian}), the derivative with respect to $\mathbf{c}$
only has meaning as a short-hand for $\left(
\frac{\partial}{\partial \mathbf{z}}\, , \, \frac{\partial}{\partial
\mathbf{\bar z}} \right)$. In (\ref{eq:chessianr}), the derivative
with respect to $\mathbf{c}$ is well-defined via the
$\mathbf{c}$-real interpretation.

\medskip
It is straightforward to show a relationship between the real
Hessian $\mathcal{H}_{\mathbf{r} \mathbf{r}}$ and the
$\mathbf{c}$-complex Hessian $\mathcal{H}_{\mathbf{c}
\mathbf{c}}^{\text{\boldmath \tiny $\mathbb{C}$}}$,
\begin{eqnarray*}
\mathcal{H}_{\mathbf{r} \mathbf{r}} & \triangleq &
\frac{\partial}{\partial \mathbf{r}} \left( \frac{\partial f
}{\partial \mathbf{r}} \right)^T \\ & = & \frac{\partial}{\partial
\mathbf{r}} \left( \frac{\partial f}{\partial \mathbf{r}} \right)^H
\\ & = & \frac{\partial}{\partial
\mathbf{r}} \left( \frac{\partial f}{\partial \mathbf{c}}\,
\mathsf{J} \right)^H \qquad \qquad \quad \text{(from equation
(\ref{eq:cogradtransfs}))}
\\ & = & \frac{\partial}{\partial
\mathbf{r}}\left\{ \mathsf{J}^H \left( \frac{\partial f}{\partial
\mathbf{c}} \right)^H \right\}
\\ & = & \frac{\partial}{\partial
\mathbf{c}}\left\{ \mathsf{J}^H \left( \frac{\partial f}{\partial
\mathbf{c}} \right)^H \right\} \mathsf{J}  \qquad \ \text{(from
equation (\ref{eq:cogradtransfs}))}\\ & = & \mathsf{J}^H
\frac{\partial}{\partial \mathbf{c}} \left( \frac{\partial
f}{\partial \mathbf{c}} \right)^H  \mathsf{J} \\
& = & \mathsf{J}^H \, \mathcal{H}_{\mathbf{c}
\mathbf{c}}^{\text{\boldmath \tiny $\mathbb{C}$}}\, \mathsf{J} \, .
\end{eqnarray*}
The resulting important relationship
\begin{equation}\label{eq:hrrhcc}
\mathcal{H}_{\mathbf{r} \mathbf{r}} = \mathsf{J}^H \,
\mathcal{H}_{\mathbf{c} \mathbf{c}}^{\text{\boldmath \tiny
$\mathbb{C}$}}\, \mathsf{J}
\end{equation}
between the real and $\mathbf{c}$-complex Hessians was derived in
\cite{VDB94a} based on the there unjustified (but true) assumption
that the second-order terms of the powers series expansions of $f$
in terms of $\mathbf{r}$ and $\mathbf{c}$-complex must be equal.
Here, we reverse this order of reasoning, and will show below the
equality of the second order terms in the $\mathbf{c}$-complex and
$\mathbf{r}$ expansions as a \emph{consequence} of
(\ref{eq:hrrhcc}).

\medskip
Note from (\ref{eq:bigjinv}) that
\begin{equation}\label{eq:hcchrr}
\mathcal{H}_{\mathbf{c} \mathbf{c}}^{\text{\boldmath \tiny
$\mathbb{C}$}} = \frac{1}{4} \, \mathsf{J} \,
\mathcal{H}_{\mathbf{r} \mathbf{r}}\, \mathsf{J}^H .
\end{equation}
Recalling that the Hessian $\mathcal{H}_{\mathbf{r} \mathbf{r}}$ is
a symmetric matrix,\footnote{In the real case, this is a general
property of the matrix of second partial derivatives of a scalar
function.} it is evident from (\ref{eq:hcchrr}) that
$\mathcal{H}_{\mathbf{c} \mathbf{c}}^{\text{\boldmath \tiny
$\mathbb{C}$}}$ is \emph{Hermitian}\footnote{As expected, as this is
a general property of the matrix of partial derivatives
$\frac{\partial}{\partial \mathbf{z}} \left( \frac{\partial
f(\mathbf{z})}{\partial \mathbf{z}} \right)^H$ of any
\emph{real-valued} function $f(\mathbf{z})$.}
\[ \mathcal{H}_{\mathbf{c} \mathbf{c}}^{\text{\boldmath \tiny
$\mathbb{C}$}} = \left(\mathcal{H}_{\mathbf{c}
\mathbf{c}}^{\text{\boldmath \tiny $\mathbb{C}$}}\right)^H \] (and
hence, like $\mathcal{H}_{\mathbf{r} \mathbf{r}}$, has real
eigenvalues), and positive definite (semidefinite) if and only
$\mathcal{H}_{\mathbf{r} \mathbf{r}}$ is positive definite
(semidefinite).

\medskip
As noted by van den Bos \cite{VDB94a}, one can now readily relate
the values of the eigenvalues of $\mathcal{H}_{\mathbf{c}
\mathbf{c}}^{\text{\boldmath \tiny $\mathbb{C}$}}$ and
$\mathcal{H}_{\mathbf{r} \mathbf{r}}$ from the fact, which follows
from (\ref{eq:bigjinv}) and (\ref{eq:hcchrr}), that
\begin{eqnarray*}
\mathcal{H}_{\mathbf{c} \mathbf{c}}^{\text{\boldmath \tiny
$\mathbb{C}$}} - \lambda I  =  \frac{1}{4} \, \mathsf{J} \,
\mathcal{H}_{\mathbf{r} \mathbf{r}}\, \mathsf{J}^H -
\frac{\lambda}{2} \mathsf{J} \mathsf{J}^H = \frac{1}{4} \mathsf{J}
\left( \mathcal{H}_{\mathbf{r} \mathbf{r}} - 2 \lambda I
\right)\mathsf{J}^H \, .
\end{eqnarray*}
This shows that the eigenvalues of the real Hessian matrix are twice
the size of the eigenvalues of the complex Hessian matrix (and, as a
consequence, must share the same condition number).\footnote{For a
Hermitian matrix, the singular values are the absolute values of the
(real) eigenvalues. Therefore the condition number, which is the
ratio of the largest to the smallest eigenvalue (assuming a full
rank matrix) is given by the ratio of the largest to smallest
eigenvalue magnitude.}

\medskip
Focussing our attention now on the second order term of
(\ref{eq:real2ndorder}), we have
\begin{eqnarray*}
\frac{1}{2} \Delta \mathbf{r}^T \,
\mathcal{H}_{\mathbf{r}\mathbf{r}} \, \Delta \mathbf{r} & = &
\frac{1}{2} \Delta \mathbf{r}^H \,
\mathcal{H}_{\mathbf{r}\mathbf{r}} \, \Delta \mathbf{r} \\
& = & \frac{1}{2} \Delta \mathbf{r}^H \, \mathsf{J}^H \,
\mathcal{H}_{\mathbf{c} \mathbf{c}}^{\text{\boldmath \tiny
$\mathbb{C}$}}\, \mathsf{J} \, \Delta
\mathbf{r} \qquad \ \ \, \text{(From equation (\ref{eq:hrrhcc}))} \\
& = & \frac{1}{2} \Delta \mathbf{c}^H \,
\mathcal{H}_{\mathbf{c}\mathbf{c}}^{\text{\boldmath \tiny
$\mathbb{C}$}} \, \Delta \mathbf{c}  \, , \qquad \qquad \text{(From
equation (\ref{eq:deltaorder}))}
\end{eqnarray*}
thereby showing the equality of the second order terms in an
expansion of a real-valued function $f$ either in terms of
$\mathbf{r}$ or $\mathbf{c}$-complex,\footnote{And thereby providing
a proof of this assumed equality in \cite{VDB94a}.}
\begin{equation}\label{eq:rc2terms}
\frac{1}{2} \Delta \mathbf{r}^T \,
\mathcal{H}_{\mathbf{r}\mathbf{r}} \, \Delta \mathbf{r} =
\frac{1}{2} \Delta \mathbf{c}^H \,
\mathcal{H}_{\mathbf{c}\mathbf{c}}^{\text{\boldmath \tiny
$\mathbb{C}$}} \, \Delta \mathbf{c} \, .
\end{equation}
Note that both of these terms are real valued.

\medskip
With the proof of the equalities \ref{eq:rcz1terms} and
\ref{eq:rc2terms}, we have (almost) completed a derivation of the
\begin{equation}\footnotesize \label{eq:creal2ndorder}
\boxed{\ \ \text{2nd-Order Expansion in $\mathbf{c}$-Complex:} \ \
f(\mathbf{c} + \Delta \mathbf{c}) = f(\mathbf{c}) + \frac{\partial
f(\mathbf{c})}{\partial \mathbf{c}} \, \Delta \mathbf{c} +
\frac{1}{2} \Delta \mathbf{c}^H \,
\mathcal{H}_{\mathbf{c}\mathbf{c}}^{\text{\boldmath \tiny
$\mathbb{C}$}}(\mathbf{c}) \, \Delta \mathbf{c} + \text{h.o.t.}\ \ }
\end{equation}
where the $\mathbf c$-complex Hessian
$\mathcal{H}_{\mathbf{c}\mathbf{c}}^{\text{\boldmath \tiny
$\mathbb{C}$}}$ is given by equation (\ref{eq:chessian}) and is
related to the real hessian $\mathcal{H}_{\mathbf{r}\mathbf{r}}$ by
equations (\ref{eq:hrrhcc}) and (\ref{eq:hcchrr}). Note that all of
the terms in (\ref{eq:creal2ndorder}) are real valued. The
derivation has not been fully completed because we have not verified
that $\Delta \mathbf{c}^H \,
\mathcal{H}_{\mathbf{c}\mathbf{c}}^{\text{\boldmath \tiny
$\mathbb{C}$}}(\mathbf{c}) \, \Delta \mathbf{c}$ is admissible in
the sense defined above. The derivation will be fully completed once
we have verified that
$\mathcal{H}_{\mathbf{c}\mathbf{c}}^{\text{\boldmath \tiny
$\mathbb{C}$}} \in \mathcal{L}(\mathcal{C},\mathcal{C})$, which we
will do below.

\medskip
The $\mathbf{c}$-complex expansion (\ref{eq:creal2ndorder}) is
\emph{not} differentiable with respect to $\mathbf{c}$-complex
\emph{itself,} which is not well defined, but, if differentiation is
required, should be instead interpreted as a short-hand, or
implicit, statement involving $\mathbf{z}$ and $\mathbf{\bar z}$,
for which derivatives are well defined. To explicitly show the
second order expansion of the real-valued function $f$ in terms of
the complex vectors $\mathbf{z}$ and $\mathbf{\bar z}$, it is
convenient to define the quantities{\small
\begin{equation}\label{eq:zhdefs}
\mathcal{H}_{\mathbf{z}\mathbf{z}} \triangleq
\frac{\partial}{\partial \mathbf{z}} \left( \frac{\partial
f}{\partial \mathbf{z}} \right)^H , \quad \mathcal{H}_{\mathbf{\bar
z}\mathbf{z}} \triangleq \frac{\partial}{\partial \mathbf{\bar z}}
\left( \frac{\partial f}{\partial \mathbf{z}} \right)^H , \quad
\mathcal{H}_{\mathbf{z}\mathbf{\bar z}} \triangleq
\frac{\partial}{\partial \mathbf{z}} \left( \frac{\partial
f}{\partial \mathbf{\bar z}} \right)^H, \quad \text{and} \quad
\mathcal{H}_{\mathbf{\bar z}\mathbf{\bar z}} \triangleq
\frac{\partial}{\partial \mathbf{\bar z}} \left( \frac{\partial
f}{\partial \mathbf{\bar z}} \right)^H .
\end{equation}}
With $\frac{\partial}{\partial \mathbf{c}} =
(\frac{\partial}{\partial \mathbf{z}} \, , \,
\frac{\partial}{\partial \mathbf{\bar z}})$, we also have from
(\ref{eq:chessian}) and the definitions (\ref{eq:zhdefs}) that
\begin{equation}\label{eq:blockhs}
\mathcal{H}_{\mathbf{c} \mathbf{c}}^{\text{\boldmath \tiny
$\mathbb{C}$}} = \begin{pmatrix} \mathcal{H}_{\mathbf{z}\mathbf{z}}
& \mathcal{H}_{\mathbf{\bar z}\mathbf{z}} \\
\mathcal{H}_{\mathbf{z}\mathbf{\bar z}} & \mathcal{H}_{\mathbf{\bar
z}\mathbf{\bar z}} \end{pmatrix} .
\end{equation}
Thus, using the earlier proven property that
$\mathcal{H}_{\mathbf{c} \mathbf{c}}^{\text{\boldmath \tiny
$\mathbb{C}$}}$ is Hermitian, $\mathcal{H}_{\mathbf{c}
\mathbf{c}}^{\text{\boldmath \tiny $\mathbb{C}$}} = \left(
\mathcal{H}_{\mathbf{c} \mathbf{c}}^{\text{\boldmath \tiny
$\mathbb{C}$}} \right)^H$, we immediately have from
(\ref{eq:blockhs}) the \emph{Hermitian conjugate conditions}
\begin{equation}\label{eq:hermitianequalities}
\mathcal{H}_{\mathbf{z}\mathbf{z}} =
\mathcal{H}_{\mathbf{z}\mathbf{z}}^H \quad \text{and} \quad
\mathcal{H}_{\mathbf{\bar z}\mathbf{z}} =
\mathcal{H}_{\mathbf{z}\mathbf{\bar z}}^H
\end{equation}
which also hold for $\mathbf{z}$ and $\mathbf{\bar z}$ replaced by
$\mathbf{\bar z}$ and $\mathbf{z}$ respectively.

\medskip
Some additional useful properties can be shown to be true for the
block components of (\ref{eq:blockhs}) defined in (\ref{eq:zhdefs}).
First note that as a consequence of $f$ being a real-valued
function, it is straightforward to show the validity of the
\emph{conjugation conditions}
\[ \overline{\mathcal{H}_{\mathbf{c} \mathbf{c}}^{\text{\boldmath \tiny
$\mathbb{C}$}}}  = \mathcal{H}_{\mathbf{\bar c} \mathbf{\bar
c}}^{\text{\boldmath \tiny $\mathbb{C}$}} \] or, equivalently,
\begin{equation}\label{eq:realimag}
 \mathcal{H}_{\mathbf{\bar z}\mathbf{\bar z}} =
\overline{\mathcal{H}_{\mathbf{z}\mathbf{z}}} \quad \text{and} \quad
\mathcal{H}_{\mathbf{\bar z}\mathbf{z}} =
\overline{\mathcal{H}_{\mathbf{z}\mathbf{\bar z}}} \, ,
\end{equation}
which also hold for $\mathbf{z}$ and $\mathbf{\bar z}$ replaced by
$\mathbf{\bar z}$ and $\mathbf{z}$ respectively.  It is also
straightforward to show that
\[ \mathcal{H}_{\mathbf{c}
\mathbf{c}}^{\text{\boldmath \tiny $\mathbb{C}$}}  = S
\mathcal{H}_{\mathbf{\bar c} \mathbf{\bar c}}^{\text{\boldmath \tiny
$\mathbb{C}$}} S = S \, \overline{\mathcal{H}_{\mathbf{c} \mathbf{
c}}^{\text{\boldmath \tiny $\mathbb{C}$}}} \ S \, , \] for $S = S^T
= S^{-1}$ (showing that $\mathcal{H}_{\mathbf{c}
\mathbf{c}}^{\text{\boldmath \tiny $\mathbb{C}$}}$ and
$\mathcal{H}_{\mathbf{\bar c} \mathbf{\bar c}}^{\text{\boldmath
\tiny $\mathbb{C}$}}$ are related by a similarity transformation and
therefore share the same eigenvalues\footnote{Their eigenvectors are
complex conjugates of each other, as reflected in the similarity
transformation being given by the swap operator $S$}), which is
precisely the necessary and sufficient condition
(\ref{eq:ccmappingcond2}) that the Hessian matrix $\mathcal{H}_{\mathbf{c}
\mathbf{c}}^{\text{\boldmath \tiny $\mathbb{C}$}}$ is admissible, $\mathcal{H}_{\mathbf{c}
\mathbf{c}}^{\text{\boldmath \tiny $\mathbb{C}$}} \in
\mathcal{L}(\mathcal{C},\mathcal{C})$. This verifies that the term
$\Delta \mathbf{c}^H \mathcal{H}_{\mathbf{c}
\mathbf{c}}^{\text{\boldmath \tiny $\mathbb{C}$}}\Delta \mathbf{c}$
is admissible and provides the completion of the proof of the
validity of (\ref{eq:creal2ndorder}) promised earlier. Finally, note
that properties (\ref{eq:realimag}) and
(\ref{eq:hermitianequalities}) yield the \emph{conjugate symmetry
conditions,}
\begin{equation}\label{eq:symconds}
\mathcal{H}_{\mathbf{z}\mathbf{z}} = \mathcal{H}_{\mathbf{\bar
z}\mathbf{\bar z}}^T \quad \text{and} \quad
\mathcal{H}_{\mathbf{z}\mathbf{\bar z}} = \mathcal{H}_{\mathbf{
z}\mathbf{\bar z}}^T \, ,
\end{equation}
which also hold for $\mathbf{z}$ and $\mathbf{\bar z}$ replaced by
$\mathbf{\bar z}$ and $\mathbf{z}$ respectively.

\medskip
From equations (\ref{eq:deltaequals}), (\ref{eq:rc2terms}), and
(\ref{eq:blockhs}) we can now expand the second order term in
(\ref{eq:real2ndorder}) as follows
\begin{eqnarray*}
\frac{1}{2}\,  \Delta \mathbf{r}^T \,
\mathcal{H}_{\mathbf{r}\mathbf{r}} \, \Delta \mathbf{r} & = &
\frac{1}{2}\, \Delta \mathbf{c}^H \,
\mathcal{H}_{\mathbf{c}\mathbf{c}}^{\text{\boldmath \tiny
$\mathbb{C}$}} \, \Delta \mathbf{c} \\
& = & \frac{1}{2} \, \left( \Delta \mathbf{z}^H
\mathcal{H}_{\mathbf{z}\mathbf{z}} \Delta \mathbf{z} + \Delta
\mathbf{ z}^H \mathcal{H}_{\mathbf{\bar z}\mathbf{z}}\Delta
\mathbf{\bar z} + \Delta \mathbf{\bar z}^H
\mathcal{H}_{\mathbf{z}\mathbf{\bar z}} \Delta \mathbf{z} + \Delta
\mathbf{\bar z}^H \mathcal{H}_{\mathbf{\bar z}\mathbf{\bar z}}
\Delta \mathbf{\bar z} \right) \\
& = & \text{Re} \left\{ \Delta \mathbf{z}^H
\mathcal{H}_{\mathbf{z}\mathbf{z}} \Delta \mathbf{z} + \Delta
\mathbf{ z}^H \mathcal{H}_{\mathbf{\bar z}\mathbf{z}}\Delta
\mathbf{\bar z} \right\}
\end{eqnarray*}
where the last step follows as a consequence of
(\ref{eq:realimag}).\footnote{Alternatively, the last step also
follows as a consequence of (\ref{eq:hermitianequalities}).} Thus,
we have so-far determined that
\begin{equation}\label{eq:2ndorderterms}  \frac{1}{2} \Delta \mathbf{r}^T \,
\mathcal{H}_{\mathbf{r}\mathbf{r}} \, \Delta \mathbf{r} =
\frac{1}{2} \Delta \mathbf{c}^H \,
\mathcal{H}_{\mathbf{c}\mathbf{c}}^{\text{\boldmath \tiny
$\mathbb{C}$}} \, \Delta \mathbf{c} = \text{Re} \left\{ \Delta
\mathbf{z}^H \mathcal{H}_{\mathbf{z}\mathbf{z}} \Delta \mathbf{z} +
\Delta \mathbf{ z}^H \mathcal{H}_{\mathbf{\bar z}\mathbf{z}}\Delta
\mathbf{\bar z} \right\} .
\end{equation}
Combining the results given in ({\ref{eq:real2ndorder}}),
(\ref{eq:rcz1terms}), and (\ref{eq:2ndorderterms}) yields the
desired expression for the second order expansion of $f$ in terms of
$\mathbf z$, {\footnotesize
\begin{equation}\label{eq:2ndorderinz}
\boxed{\ \ \text{$2^{nd}$-Order Exp.\ in $\mathbf{z}$:} \quad
f(\mathbf{z} + \Delta \mathbf{z}) = f(\mathbf{z})+  2\,  \text{Re}
\left\{ \frac{\partial f}{\partial \mathbf{z}} \, \Delta
\mathbf{z}\right\} +  \text{Re} \left\{\Delta \mathbf{z}^H
\mathcal{H}_{\mathbf{z}\mathbf{z}} \Delta \mathbf{z} + \Delta
\mathbf{ z}^H \mathcal{H}_{\mathbf{\bar z}\mathbf{z}}\Delta
\mathbf{\bar z} \right\}+ \text{h.o.t.} \ \ }
\end{equation}}
\hspace{-.15in} We note in passing that Equation (\ref{eq:2ndorderinz}) is exactly
the same expression given as Equation (A.7) of reference \cite{TA91}
and Equation (8) of reference \cite{GY00}, which were both derived
via an alternative procedure.

\medskip
The $\mathbf{c}$-complex expansion shown in Equation
(\ref{eq:creal2ndorder}) is one of two possible alternative
second-order representations in $\mathbf{c}$ for $f(\mathbf{c})$
(the other being the $\mathbf c$-real expansion), and was used as
the starting point of the theoretical developments leading to the
$\mathbf{z}$-expansion (\ref{eq:2ndorderinz}). We now turn to the
development of the $\mathbf{c}$-real expansion of $f(\mathbf{c})$,
which will be accomplished by writing the second order term of the
quadratic expansion in terms of the $\mathbf{c}$-real Hessian
$\mathcal{H}_{\mathbf{c}\mathbf{c}}^{\text{\boldmath \tiny
$\mathbb{R}$}}$.

\medskip
From the definitions (\ref{eq:chessianr}), (\ref{eq:chessian}), and
(\ref{eq:zhdefs}), and  using the fact that
$\frac{\partial}{\partial \mathbf{c}} = (\frac{\partial}{\partial
\mathbf{z}} \, , \, \frac{\partial}{\partial \mathbf{\bar z}})$, it
is straightforward to show that
\begin{equation}\label{eq:hrhc1}
\mathcal{H}_{\mathbf{c} \mathbf{c}}^{\text{\boldmath \tiny
$\mathbb{R}$}} = \begin{pmatrix} \mathcal{H}_{\mathbf{z}\mathbf{\bar
z}}
& \mathcal{H}_{\mathbf{\bar z}\mathbf{\bar z}} \\
\mathcal{H}_{\mathbf{z}\mathbf{z}} & \mathcal{H}_{\mathbf{\bar
z}\mathbf{z}} \end{pmatrix} = S \begin{pmatrix}
\mathcal{H}_{\mathbf{z}\mathbf{z}}
& \mathcal{H}_{\mathbf{\bar z}\mathbf{z}} \\
\mathcal{H}_{\mathbf{z}\mathbf{\bar z}} & \mathcal{H}_{\mathbf{\bar
z}\mathbf{\bar z}} \end{pmatrix}
\end{equation}
or\footnote{Alternative derivations are possible. For example,
$\mathcal{H}_{\mathbf{c} \mathbf{c}}^{\text{\boldmath \tiny
$\mathbb{C}$}} = \frac{\partial}{\partial
\mathbf{c}}\left(\frac{\partial f}{\partial \mathbf{c}} \right)^H =
\frac{\partial}{\partial \mathbf{c}}\left(\frac{\partial f}{\partial
\mathbf{\bar c}} \right)^T = \frac{\partial}{\partial
\mathbf{c}}\left(\frac{\partial f}{\partial \mathbf{c}} S \right)^T
= \frac{\partial}{\partial \mathbf{c}} S \left(\frac{\partial
f}{\partial \mathbf{c}} \right)^T = S \frac{\partial}{\partial
\mathbf{c}}\left(\frac{\partial f}{\partial \mathbf{c}} \right)^T =
S \mathcal{H}_{\mathbf{c} \mathbf{c}}^{\text{\boldmath \tiny
$\mathbb{R}$}} \Rightarrow \mathcal{H}_{\mathbf{c}
\mathbf{c}}^{\text{\boldmath \tiny $\mathbb{R}$}} = S
\mathcal{H}_{\mathbf{c} \mathbf{c}}^{\text{\boldmath \tiny
$\mathbb{C}$}}$, noting that $S = S^T = S^{-1}$.}
\begin{equation}\label{eq:hrhc2}
\mathcal{H}_{\mathbf{c} \mathbf{c}}^{\text{\boldmath \tiny
$\mathbb{R}$}} =  \mathcal{H}_{\mathbf{c} \mathbf{\bar
c}}^{\text{\boldmath \tiny $\mathbb{C}$}}= S \mathcal{H}_{\mathbf{c}
\mathbf{c}}^{\text{\boldmath \tiny $\mathbb{C}$}}  =
\mathcal{H}_{\mathbf{\bar c} \mathbf{\bar c}}^{\text{\boldmath \tiny
$\mathbb{C}$}} S .
\end{equation}
Note from the first equality in (\ref{eq:hrhc1}) and the conjugate
symmetry conditions (\ref{eq:symconds}) that the $\mathbf c$-real
Hessian is \emph{symmetric}
\begin{equation}\label{eq:hrsymmetric}
\mathcal{H}_{\mathbf{c} \mathbf{c}}^{\text{\boldmath \tiny
$\mathbb{R}$}} = \left( \mathcal{H}_{\mathbf{c}
\mathbf{c}}^{\text{\boldmath \tiny $\mathbb{R}$}} \right)^T.
\end{equation}
Equivalently,
\begin{equation}\label{eq:hrsymmetric-b}
S \mathcal{H}_{\mathbf{c} \mathbf{c}}^{\text{\boldmath \tiny
$\mathbb{C}$}} = \left( S \mathcal{H}_{\mathbf{c}
\mathbf{c}}^{\text{\boldmath \tiny $\mathbb{C}$}} \right)^T.
\end{equation}

\medskip
Let the Singular Value Decomposition (SVD) of $\mathcal{H}_{\mathbf{c} \mathbf{c}}^{\text{\boldmath
\tiny $\mathbb{C}$}}$ be
\[ \mathcal{H}_{\mathbf{c} \mathbf{c}}^{\text{\boldmath \tiny
$\mathbb{C}$}} = U \Sigma V^H \] then from (\ref{eq:hrhc2}) the SVD
of $\mathcal{H}_{\mathbf{c} \mathbf{c}}^{\text{\boldmath \tiny
$\mathbb{R}$}}$ is given by
\[ \mathcal{H}_{\mathbf{c} \mathbf{c}}^{\text{\boldmath \tiny
$\mathbb{R}$}} = U' \Sigma V^H \, , \quad U' = S U \] showing that
$\mathcal{H}_{\mathbf{c} \mathbf{c}}^{\text{\boldmath \tiny
$\mathbb{C}$}}$ and $\mathcal{H}_{\mathbf{c}
\mathbf{c}}^{\text{\boldmath \tiny $\mathbb{R}$}}$ share the same
singular values, and hence the same condition number (which is given
by the ratio of the largest to smallest singular value). The three
Hessian matrices $\mathcal{H}_{\mathbf{r} \mathbf{r}}$,
$\mathcal{H}_{\mathbf{c} \mathbf{c}}^{\text{\boldmath \tiny
$\mathbb{R}$}}$, and $\mathcal{H}_{\mathbf{c}
\mathbf{c}}^{\text{\boldmath \tiny $\mathbb{C}$}}$ are essentially
equivalent for investigating numerical issues and for testing
whether a proposed minimizer of the second order expansion of
$f(\mathbf{r}) = f(\mathbf{c})$ is a local (or even global) minimum.
Thus, one can choose to work with the Hessian matrix which is
easiest to compute and analyze. This is usually the $\mathbf
c$-complex Hessian $\mathcal{H}_{\mathbf{c}
\mathbf{c}}^{\text{\boldmath \tiny $\mathbb{C}$}}$, and it is often
most convenient to determine numerical stability and optimality
using $\mathcal{H}_{\mathbf{c} \mathbf{c}}^{\text{\boldmath \tiny
$\mathbb{C}$}}$ even when the algorithm is being developed from one
of the alternative perspectives (i.e., the real $\mathbf r$ or the
$\mathbf c$-real second order expansion).

\medskip
Now note that from (\ref{eq:hrhc2}) we immediately and easily have
\[\small \frac{1}{2}\, \Delta \mathbf{c}^T \,
\mathcal{H}_{\mathbf{c}\mathbf{c}}^{\text{\boldmath \tiny
$\mathbb{R}$}} \, \Delta \mathbf{c} = \frac{1}{2}\, \Delta
\mathbf{c}^T \,S \,
\mathcal{H}_{\mathbf{c}\mathbf{c}}^{\text{\boldmath \tiny
$\mathbb{C}$}} \, \Delta \mathbf{c} = \frac{1}{2}\, ( S \Delta
\mathbf{c})^T \, \mathcal{H}_{\mathbf{c}\mathbf{c}}^{\text{\boldmath
\tiny $\mathbb{C}$}} \, \Delta \mathbf{c} = \frac{1}{2}\,
\overline{(\Delta \mathbf{c})}^T \,
\mathcal{H}_{\mathbf{c}\mathbf{c}}^{\text{\boldmath \tiny
$\mathbb{C}$}} \, \Delta \mathbf{c} = \frac{1}{2}\, \Delta
\mathbf{c}^H \, \mathcal{H}_{\mathbf{c}\mathbf{c}}^{\text{\boldmath
\tiny $\mathbb{C}$}} \, \Delta \mathbf{c} \] showing the equivalence
of the $\mathbf c$-real and $\mathbf c$-complex second order terms
in the expansion of $f(\mathbf{c})$.\footnote{One can show that the
term $\Delta \mathbf{c}^T \,
\mathcal{H}_{\mathbf{c}\mathbf{c}}^{\text{\boldmath \tiny
$\mathbb{R}$}} \, \Delta \mathbf{c}$ is admissible if and only if
$\mathcal{H}_{\mathbf{c}\mathbf{c}}^{\text{\boldmath \tiny
$\mathbb{R}$}} = S M$ for $M \in
\mathcal{L}(\mathcal{C},\mathcal{C})$, which is the case here.}
Combining this result with (\ref{eq:2ndorderterms}), we have shown
the following equivalences between the second order terms in the
various expansions of $f$ under consideration in this note:
\begin{equation}\footnotesize \label{eq:2ndorderterms2} \boxed{\  \text{2nd-Order
Terms:} \ \  \frac{1}{2} \Delta \mathbf{r}^T \,
\mathcal{H}_{\mathbf{r}\mathbf{r}} \, \Delta \mathbf{r} =
\frac{1}{2}\, \Delta \mathbf{c}^T \,
\mathcal{H}_{\mathbf{c}\mathbf{c}}^{\text{\boldmath \tiny
$\mathbb{R}$}} \, \Delta \mathbf{c} = \frac{1}{2} \Delta
\mathbf{c}^H \, \mathcal{H}_{\mathbf{c}\mathbf{c}}^{\text{\boldmath
\tiny $\mathbb{C}$}} \, \Delta \mathbf{c} = \text{Re} \left\{ \Delta
\mathbf{z}^H \mathcal{H}_{\mathbf{z}\mathbf{z}} \Delta \mathbf{z} +
\Delta \mathbf{ z}^H \mathcal{H}_{\mathbf{\bar z}\mathbf{z}}\Delta
\mathbf{\bar z} \right\}  \  }
\end{equation}
where the second order expansion in $\mathbf{r}$ is given by
(\ref{eq:real2ndorder}), the $\mathbf{c}$-complex expansion by
(\ref{eq:creal2ndorder}), the expansion in terms of $\mathbf z$ by
(\ref{eq:2ndorderinz}), and the $\mathbf{c}$-real expansion by
\begin{equation}\small \label{eq:creal2ndorder2}
\boxed{\ \ \text{2nd-Order Expansion in $\mathbf{c}$-Real:} \ \
f(\mathbf{c} + \Delta \mathbf{c}) = f(\mathbf{c}) + \frac{\partial
f(\mathbf{c})}{\partial \mathbf{c}} \, \Delta \mathbf{c} +
\frac{1}{2} \Delta \mathbf{c}^T \,
\mathcal{H}_{\mathbf{c}\mathbf{c}}^{\text{\boldmath \tiny
$\mathbb{R}$}}(\mathbf{c}) \, \Delta \mathbf{c} + \text{h.o.t.}\ \ }
\end{equation}
Note that all of the terms in (\ref{eq:2ndorderterms2}) and
(\ref{eq:creal2ndorder2}) are real valued.

\medskip
The expansion in of $f(\mathbf{c})$ in terms of
\emph{$\mathbf{c}$-complex} shown in (\ref{eq:creal2ndorder}) is
\emph{not} differentiable with respect to $\mathbf c$ as differentiation with respect to \emph{$\mathbf{c}$-complex} is not defined. (Recall, though, that we can differentiate the \emph{$\mathbf{c}$-real} expansion with respect to \emph{$\mathbf{c}$-real}.) However,
(\ref{eq:creal2ndorder}) \emph{is} differentiable with respect to $\mathbf
z$ and $\mathbf{\bar z}$ and can be viewed as a short-hand
equivalent to the full $(\mathbf{z},\mathbf{\bar z})$ expansion
provided by (\ref{eq:2ndorderinz}). Therefore, it is Equation
(\ref{eq:2ndorderinz}) which is the natural form for optimization
with respect to $\mathbf{c}$-complex via a derivative-based
approach, because only differentiation with respect to the
components $(\mathbf{z},\mathbf{\bar z})$ of $\mathbf c$-complex is
well-posed. On the other hand, differentiation with respect to
$\mathbf{c}$-real is well-posed, so that one can optimize
(\ref{eq:creal2ndorder2}) by taking derivatives of
(\ref{eq:creal2ndorder2}) with respect to $\mathbf{c}$-real itself.

\medskip
 Note that (\ref{eq:real2ndorder}), (\ref{eq:creal2ndorder}), and
(\ref{eq:creal2ndorder2}) are the natural forms to use for
optimization via ``completing the square.'' This is
because the expansions in terms of $\mathbf r$,
$\mathbf{c}$-complex, and $\mathbf{c}$-real are less awkward for
completing-the-square purposes than the expansion in $\mathbf z$
provided by (\ref{eq:2ndorderinz}).\footnote{Although
(\ref{eq:2ndorderinz}) can also be optimized by completing the
square.} Note, further that the expansions (\ref{eq:real2ndorder}) and
(\ref{eq:creal2ndorder}) both have a form amenable to
optimization by completing the square \emph{and}  both are differentiable with respect to the
expansion variable itself.

\medskip
The various second order expansions developed above can be found in
references \cite{TA91}, \cite{VDB94a} and \cite{GY00}. In
\cite{VDB94a}, van den Bos shows the equality of the first, second,
and third second-order terms shown in equation
(\ref{eq:2ndorderterms}) but does not mention the fourth (which,
anyway, naturally follows from the third term in
(\ref{eq:2ndorderterms}) via a simple further expansion in terms of
$\mathbf z$ and $\mathbf{\bar z}$). Indeed, the approach used in this note
is a more detailed elaboration of the derivations presented by van den Bos in
\cite{VDB94a}. In reference \cite{GY00} Yan and Fan show the
equality of the first and last terms in (\ref{eq:2ndorderterms}),
but, while they cite the results of van den Bos \cite{VDB94a}
regarding the middle terms in (\ref{eq:2ndorderterms}), do not
appear to have appreciated that the fourth term in
(\ref{eq:2ndorderterms}) is an \emph{immediate consequence} of the
second or third terms, and instead derived it from scratch using an
alternative, ``brute force'' approach.

\paragraph{Quadratic Minimization and the Newton Algorithm.} The
Newton algorithm for minimizing a scalar function $f(z)$ exploits
the fact that it is generally straightforward to minimize the
quadratic approximations provided by  second order expansions such
as (\ref{eq:real2ndorder}), (\ref{eq:creal2ndorder}),
(\ref{eq:2ndorderinz}), and (\ref{eq:creal2ndorder2}). The Newton
method starts with an initial estimate of the optimal solution, say
$\mathbf{\hat c}$, then expands $f(\mathbf{c})$ about the estimate
$\mathbf{\hat c}$ to second order in $\Delta \mathbf{c} = \mathbf{c}
- \mathbf{\hat c}$, and then minimizes the resulting second order
approximation of $f(\mathbf{c})$ with respect to $\Delta
\mathbf{c}$. Having determined an estimated update $\widehat{\Delta
\mathbf{c}}$ in this manner, one updates the original estimate
$\mathbf{\hat c} \leftarrow \mathbf{\hat c} + \alpha \widehat{\Delta
\mathbf{c}}$, for some small ``stepsize'' $\alpha
> 0$, and then starts the optimization cycle all over again. For
appropriate choices of the stepsize $\alpha$, this iterative
approximate quadratic optimization algorithm can  result in a
sequence of estimates $\mathbf{\hat c}_0$, $\mathbf{\hat c}_1$,
$\mathbf{\hat c}_2$, $\cdots$, which converges to the true optimal
solution extremely quickly \cite{Lu69}.

\medskip
Note that the optimal solution to the quadratic approximations
provided by (\ref{eq:real2ndorder}), (\ref{eq:creal2ndorder}),  and
(\ref{eq:creal2ndorder2}) can be \emph{immediately} written down
using the ``completing-the-square'' procedure assuming that the relevant Hessians are all
invertible:
\begin{eqnarray}
\widehat{\Delta \mathbf{r}} \, \, & = & - \text{\large
$\left(\mathcal{H}_{\mathbf{r}\mathbf{r}}\right)^{-1}$}
\text{\footnotesize $\left(\frac{\partial f(\mathbf{r})}{\partial
\mathbf{r}}\right)^T$}
\qquad \ \left(\text{from the $\mathbf{r}$ expansion (\ref{eq:real2ndorder})}\right) \label{eq:solnr} \\
\widehat{\Delta \mathbf{c}^{\text{\boldmath \tiny $\mathbb{C}$}}} &
= & - \text{\large
$\left(\mathcal{H}_{\mathbf{c}\mathbf{c}}^{\text{\boldmath \tiny
$\mathbb{C}$}}   \right)^{-1}$} \text{\footnotesize
$\left(\frac{\partial f(\mathbf{c})}{\partial \mathbf{c}}\right)^H$}
\qquad \left(\text{from the $\mathbf{c}$-complex expansion (\ref{eq:creal2ndorder})}\right) \label{eq:solnccomplex} \\
\widehat{\Delta \mathbf{c}^{\text{\boldmath \tiny $\mathbb{R}$}}} &
= & - \text{\large
$\left(\mathcal{H}_{\mathbf{c}\mathbf{c}}^{\text{\boldmath \tiny
$\mathbb{R}$}}   \right)^{-1}$} \text{\footnotesize
$\left(\frac{\partial f(\mathbf{c})}{\partial \mathbf{c}}\right)^T$}
\qquad \, \left(\text{from the $\mathbf{c}$-real expansion
(\ref{eq:creal2ndorder2})}\right) . \label{eq:solncreal}
\end{eqnarray}
Solutions (\ref{eq:solnr}) and (\ref{eq:solnccomplex})  can also be
found in van den Bos \cite{VDB94a}. Note that $\widehat{\Delta
\mathbf{c}^{\text{\boldmath \tiny $\mathbb{C}$}}}$ is an
\emph{admissible} solution, i.e., that
\[ \widehat{\Delta \mathbf{c}^{\text{\boldmath \tiny $\mathbb{C}$}}}
\in \mathcal{C} \,  \] as required for self-consistency of our
theory, as a consequence of the fact that $\text{\footnotesize
$\left(\frac{\partial f(\mathbf{c})}{\partial
\mathbf{c}}\right)^H$}$ and
$\left(\mathcal{H}_{\mathbf{c}\mathbf{c}}^{\text{\boldmath \tiny
$\mathbb{C}$}}   \right)^{-1}$ satisfy
\[ \left(\frac{\partial f(\mathbf{c})}{\partial
\mathbf{c}}\right)^H \in \mathcal{C} \quad \text{and} \quad
\text{\large
$\left(\mathcal{H}_{\mathbf{c}\mathbf{c}}^{\text{\boldmath \tiny
$\mathbb{C}$}}   \right)^{-1}$} \in
\mathcal{L}(\mathcal{C},\mathcal{C}) \, , \] with the latter
condition a consequence of property (\ref{eq:invtoo}) and the fact
that $\text{\large
$\mathcal{H}_{\mathbf{c}\mathbf{c}}^{\text{\boldmath \tiny
$\mathbb{C}$}}$} \in \mathcal{L}(\mathcal{C},\mathcal{C})$. If this
were not the case, then we generally would have the
\emph{meaningless} answer that $\widehat{\Delta
\mathbf{c}^{\text{\boldmath \tiny $\mathbb{C}$}}} \notin
\mathcal{C}$.

\medskip
The admissibility of the solution (\ref{eq:solncreal}) follows from
the admissibility of (\ref{eq:solnccomplex}). This will be evident
from the fact, as we shall show, that all of the solutions
(\ref{eq:solnr})-(\ref{eq:solncreal}) must all correspond to the
\emph{same update,} \[ \widehat{\Delta \mathbf{c}^{\text{\boldmath
\tiny $\mathbb{C}$}}} = \widehat{\Delta \mathbf{c}^{\text{\boldmath
\tiny $\mathbb{R}$}}} = \mathsf{J}\widehat{\Delta \mathbf{r}} \, .
\]

\medskip Note that
\begin{eqnarray*}
\widehat{\Delta \mathbf{c}^{\text{\boldmath \tiny $\mathbb{C}$}}} &
= & - \text{\large
$\left(\mathcal{H}_{\mathbf{c}\mathbf{c}}^{\text{\boldmath \tiny
$\mathbb{C}$}}   \right)^{-1}$} \text{\footnotesize
$\left(\frac{\partial f(\mathbf{c})}{\partial \mathbf{c}}\right)^H$}
\\ & = & - \text{
$\left(\frac{1}{4} \mathsf{J} \mathcal{H}_{\mathbf{r}\mathbf{r}}
\mathsf{J}^H \right)^{-1}$} \text{\small $\left(\frac{1}{2}
\frac{\partial f(\mathbf{r})}{\partial \mathbf{r}} \mathsf{J}^H
\right)^H$} \qquad \text{(from (\ref{eq:cogradtransfs}) and
(\ref{eq:hcchrr}))}
\\ & = & \\ & = & - \text{
$\left( \mathsf{J} \mathcal{H}_{\mathbf{r}\mathbf{r}}
\mathsf{J}^{-1} \right)^{-1}$} \text{\small $\mathsf{J} \left(
\frac{\partial f(\mathbf{r})}{\partial \mathbf{r}}  \right)^T$}
\qquad \text{(from (\ref{eq:invjacobian}))} \\ & = &
 - \mathsf{J} \text{\large
$\left(\mathcal{H}_{\mathbf{r}\mathbf{r}}\right)^{-1}$}
\text{\footnotesize $\left(\frac{\partial f(\mathbf{r})}{\partial
\mathbf{r}}\right)^T$} \\
& = & \mathsf{J}\widehat{\Delta \mathbf{r}}
\end{eqnarray*}
as required. On the other hand,
\begin{eqnarray*} \widehat{\Delta
\mathbf{c}^{\text{\boldmath \tiny $\mathbb{R}$}}} & = & -
\text{\large
$\left(\mathcal{H}_{\mathbf{c}\mathbf{c}}^{\text{\boldmath \tiny
$\mathbb{R}$}}   \right)^{-1}$} \text{\footnotesize
$\left(\frac{\partial f(\mathbf{c})}{\partial \mathbf{c}}\right)^T$}
\\ & = & -
\text{ $\left(S \mathcal{H}_{\mathbf{c}\mathbf{c}}^{\text{\boldmath
\tiny $\mathbb{C}$}}   \right)^{-1}$} \text{\footnotesize
$\left(\frac{\partial f(\mathbf{c})}{\partial \mathbf{c}}\right)^T$}
\qquad \text{(from (\ref{eq:hrhc2}))}
\\ & = & -
\text{\large $\left(
\mathcal{H}_{\mathbf{c}\mathbf{c}}^{\text{\boldmath \tiny
$\mathbb{C}$}}   \right)^{-1}$} \text{\footnotesize
$\left(\frac{\partial f(\mathbf{c})}{\partial \mathbf{c}} S
\right)^T$}
\\ & = & -
\text{\large $\left(
\mathcal{H}_{\mathbf{c}\mathbf{c}}^{\text{\boldmath \tiny
$\mathbb{C}$}}   \right)^{-1}$} \text{\footnotesize
$\left(\frac{\partial f(\mathbf{c})}{\partial \mathbf{\bar c}}
\right)^T$}
\\ & = & -
\text{\large $\left(
\mathcal{H}_{\mathbf{c}\mathbf{c}}^{\text{\boldmath \tiny
$\mathbb{C}$}}   \right)^{-1}$} \text{\footnotesize
$\left(\frac{\partial f(\mathbf{c})}{\partial \mathbf{c}}
\right)^H$} \\
& = & \widehat{\Delta \mathbf{c}^{\text{\boldmath \tiny
$\mathbb{C}$}}} .
\end{eqnarray*}
Thus, the updates (\ref{eq:solnr})-(\ref{eq:solncreal}) are indeed
equivalent.

\medskip
The updates (\ref{eq:solnr}) and (\ref{eq:solncreal}), determined
via a completing the square argument, can alternatively be obtained
by setting the (\emph{real}) derivatives of their respective
quadratically-approximated loss functions to zero, and solving the
necessary condition for an optimum. Note that if we attempt to
(erroneously) take the (\emph{complex}) derivative of
(\ref{eq:creal2ndorder}) with respect to $\mathbf{c}$-complex and
then set this expression to zero, the resulting ``solution'' will be
off by a factor of two.\footnote{In a numerical solution procedure a constant factor error in the updates can be absorbed into the update step-size factor and therefore will likely not be noticed in simulations or applications. However, the claim that a specific step-size values results in stable or unstable convergence might not be confirmed in an experiment using the correctly computed updates.} In the latter case, we must instead take the
derivatives of (\ref{eq:2ndorderinz}) with respect to $\mathbf{z}$
and $\mathbf{\bar z}$ and set the resulting expressions to zero in
order to obtain the optimal solution.\footnote{This is the procedure
used in \cite{TA91} and \cite{GY00}.}

\medskip
At convergence, the Newton algorithm will produce a solution to the
necessary first-order condition
\[ \frac{\partial f(\mathbf{\hat c})}{\partial \mathbf{c}} = 0 \, ,
\]
and this point will be a local minimum of $f(\cdot)$ if the Hessians
are strictly positive definite at this point. Typically, one would
verify positive definiteness of the $\mathbf{c}$-complex Hessian at
the solution point $\mathbf{\hat c}$,
\[ \mathcal{H}_{\mathbf{c}
\mathbf{c}}^{\text{\boldmath \tiny $\mathbb{C}$}}(\mathbf{\hat c})=
\begin{pmatrix} \mathcal{H}_{\mathbf{z}\mathbf{z}}(\mathbf{\hat c})
& \mathcal{H}_{\mathbf{\bar z}\mathbf{z}}(\mathbf{\hat c}) \\
\mathcal{H}_{\mathbf{z}\mathbf{\bar z}}(\mathbf{\hat c}) &
\mathcal{H}_{\mathbf{\bar z}\mathbf{\bar z}}(\mathbf{\hat c})
\end{pmatrix} > 0 \, . \]

\medskip
As done in \cite{TA91} and \cite{GY00}, the solution to the
quadratic minimization problem provided by
(\ref{eq:solnr})-(\ref{eq:solncreal}) can be expressed in a closed
form expression which directly produces the solution $\mathbf{\hat
z} \in \mathbb{C}^{n}$. To do so, we rewrite the solution
(\ref{eq:solnccomplex}) for the Newton update $\widehat{\Delta
\mathbf{c}}$ as
\[\text{\large $\mathcal{H}_{\mathbf{c}\mathbf{c}}^{\text{\boldmath \tiny
$\mathbb{C}$}}$} \,  \widehat{\Delta \mathbf{c}}  = -
\text{\footnotesize $\left(\frac{\partial f(\mathbf{c})}{\partial
\mathbf{c}}\right)^H$}
\]
which we then write in expanded form in terms of $\mathbf{z}$ and
$\mathbf{\bar z}$ \begin{equation}\label{eq:blockmatrix}
\text{\large $\begin{pmatrix} \mathcal{H}_{\mathbf{z}\mathbf{z}}
& \mathcal{H}_{\mathbf{\bar z}\mathbf{z}} \\
\mathcal{H}_{\mathbf{z}\mathbf{\bar z}} & \mathcal{H}_{\mathbf{\bar
z}\mathbf{\bar z}} \end{pmatrix}\begin{pmatrix} \widehat{\Delta \mathbf{z}} \\
\widehat{\Delta \mathbf{\bar z}} \end{pmatrix}$} = - \begin{pmatrix}
\text{\footnotesize
$\left(\frac{\partial f}{\partial \mathbf{z}}\right)^H$} \\
\text{\footnotesize $\left(\frac{\partial f}{\partial \mathbf{\bar
z}}\right)^H$}
\end{pmatrix}\, .
\end{equation}
Assuming that $\mathcal{H}_{\mathbf{c}\mathbf{c}}^{\text{\boldmath
\tiny $\mathbb{C}$}}$ is positive definite, then
$\mathcal{H}_{\mathbf{z}\mathbf{z}}$ is invertible and the second
block row in (\ref{eq:blockmatrix}) results in
\[ \widehat{\Delta \mathbf{\bar z}} = - \mathcal{H}_{\mathbf{\bar
z}\mathbf{\bar z}}^{-1}\mathcal{H}_{\mathbf{z}\mathbf{\bar z}}
\widehat{\Delta \mathbf{z}} - \mathcal{H}_{\mathbf{\bar
z}\mathbf{\bar z}}^{-1} \text{\footnotesize $\left(\frac{\partial
f}{\partial \mathbf{\bar z}}\right)^H$}. \] Plugging this into the
first block row of (\ref{eq:blockmatrix}) then yields the Newton
algorithm update equation
\begin{equation}\label{eq:zoptimal}
\widetilde{\mathcal{H}_{\mathbf{z}\mathbf{z}}}\,  \widehat{\Delta
\mathbf{z}} = - \text{\footnotesize $\left(\frac{\partial
f}{\partial \mathbf{z}}\right)^H$} + \mathcal{H}_{\mathbf{\bar
z}\mathbf{z}} \mathcal{H}_{\mathbf{\bar z}\mathbf{\bar z}}^{-1}
\text{\footnotesize $\left(\frac{\partial f}{\partial \mathbf{\bar
z}}\right)^H$} \, ,
\end{equation}
where \[ \widetilde{\mathcal{H}_{\mathbf{z}\mathbf{z}}} \triangleq
\mathcal{H}_{\mathbf{z}\mathbf{z}} - \mathcal{H}_{\mathbf{\bar
z}\mathbf{z}} \mathcal{H}_{\mathbf{\bar z}\mathbf{\bar
z}}^{-1}\mathcal{H}_{\mathbf{z}\mathbf{\bar z}} \] is the
\emph{Schur complement} of $\mathcal{H}_{\mathbf{z}\mathbf{z}}$ in
$\mathcal{H}_{\mathbf{c}\mathbf{c}}^{\text{\boldmath \tiny
$\mathbb{C}$}}$. Equation (\ref{eq:zoptimal}) is equivalent to the
solution given as Equation (A.12) in \cite{TA91}. Invertibility of
the Schur complement
$\widetilde{\mathcal{H}_{\mathbf{z}\mathbf{z}}}$ follows from our
assumption that $\mathcal{H}_{\mathbf{c}\mathbf{c}}^{\text{\boldmath
\tiny $\mathbb{C}$}}$ is positive definite, and
 the Newton update is therefore given by
\begin{equation}\label{eq:zoptimal2}
\widehat{\Delta \mathbf{z}}= \left(
\mathcal{H}_{\mathbf{z}\mathbf{z}} - \mathcal{H}_{\mathbf{\bar
z}\mathbf{z}} \mathcal{H}_{\mathbf{\bar z}\mathbf{\bar
z}}^{-1}\mathcal{H}_{\mathbf{z}\mathbf{\bar z}}\right)^{-1} \left\{
\mathcal{H}_{\mathbf{\bar z}\mathbf{z}} \mathcal{H}_{\mathbf{\bar
z}\mathbf{\bar z}}^{-1} \text{\footnotesize $\left(\frac{\partial
f}{\partial \mathbf{\bar z}}\right)^H$} - \text{\footnotesize
$\left(\frac{\partial f}{\partial \mathbf{z}}\right)^H$} \right\}.
\end{equation}

\medskip
The matrices $\mathcal{H}_{\mathbf{\bar z}\mathbf{\bar z}}$ and
$\widetilde{\mathcal{H}_{\mathbf{\bar z}\mathbf{\bar z}}} = \left(
\mathcal{H}_{\mathbf{z}\mathbf{z}} - \mathcal{H}_{\mathbf{\bar
z}\mathbf{z}} \mathcal{H}_{\mathbf{\bar z}\mathbf{\bar
z}}^{-1}\mathcal{H}_{\mathbf{z}\mathbf{\bar z}}\right)$ in
(\ref{eq:zoptimal}) are invertible if and only if
$\mathcal{H}_{\mathbf{c}\mathbf{c}}^{\text{\boldmath \tiny
$\mathbb{C}$}}$ is invertible. Note  that invertibility of
$\mathcal{H}_{\mathbf{z}\mathbf{z}}$ (equivalently,
$\mathcal{H}_{\mathbf{\bar z}\mathbf{\bar z}} =
\overline{\mathcal{H}_{\mathbf{z}\mathbf{z}}}$) is not a sufficient
condition for the Schur complement to be nonsingular. However, if
$\mathcal{H}_{\mathbf{\bar z}\mathbf{z}} =
\overline{\mathcal{H}_{\mathbf{z}\mathbf{\bar z}}} = 0$ then
invertibility of $\mathcal{H}_{\mathbf{z}\mathbf{z}}$ is a necessary
and sufficient condition for a solution  $\widehat{\Delta
\mathbf{z}}$ to exist.

\medskip As noted by Yan \& Fan \cite{GY00}, the need to guarantee positive
definiteness of the Schur complement
$\widetilde{\mathcal{H}_{\mathbf{\bar z}\mathbf{\bar z}}} = \left(
\mathcal{H}_{\mathbf{z}\mathbf{z}} - \mathcal{H}_{\mathbf{\bar
z}\mathbf{z}} \mathcal{H}_{\mathbf{\bar z}\mathbf{\bar
z}}^{-1}\mathcal{H}_{\mathbf{z}\mathbf{\bar z}}\right)$ is a
significant computational burden for an on-line adaptive filtering
algorithm to bear. For this reason, to improve the numerical
robustness of the Newton algorithm and to provide a substantial
simplification, they suggest making the approximation that the block
off-diagonal elements of
$\mathcal{H}_{\mathbf{c}\mathbf{c}}^{\text{\boldmath \tiny
$\mathbb{C}$}}$ are zero \[ \mathcal{H}_{\mathbf{\bar z}\mathbf{z}}
= \overline{\mathcal{H}_{\mathbf{z}\mathbf{\bar z}}} \approx 0\]
which results in the simpler approximate solution
\begin{equation}\label{eq:fansimplification} \widehat{\Delta
\mathbf{ z}} \approx - \mathcal{H}_{\mathbf{z}\mathbf{z}}^{-1}
\text{\footnotesize $\left(\frac{\partial f}{\partial
\mathbf{z}}\right)^H$} .
\end{equation}
 The argument given by Yan and Fan supporting the use of
the approximation $\mathcal{H}_{\mathbf{\bar z}\mathbf{z}} \approx
0$ is that as the Newton algorithm converges to the optimal solution
$\boldmath{\hat z} = \boldmath{z}_0$, setting
$\mathcal{H}_{\mathbf{\bar z}\mathbf{z}}$ ``to zero implies that we
will use a quadratic function to approximate the cost near
$\mathbf{z}_0$'' \cite{GY00}. However Yan and Fan do not give a
formal definition of a ``quadratic function'' and this statement is
\emph{not} \emph{generally} true as there is no a priori reason why
the off-diagonal block matrix elements of the Newton Hessian should
be zero, or approach zero, as we demonstrate in Example 2 of the
Applications section below.

\medskip
However, as we shall discuss later below, setting the block
off-diagonal elements to zero \emph{is} justifiable, but not
necessarily as an \emph{approximation} to the Newton algorithm.
Setting the block off-diagonal elements in the Newton Hessian to
zero, results in an \emph{alternative}, ``quasi-Newton'' algorithm
\emph{which can be studied in its own right} as a competitor algorithm to
the Newton algorithm, the Gauss-Newton algorithm, or the gradient
descent algorithm.\footnote{That is not to say that there can't be
conditions under which the quasi-Newton algorithm does converge to
the Newton algorithm. Just as one can give conditions for which the
Gauss-Newton algorithm converges to the Newton algorithm, one should
be able to do the same for the quasi-Newton algorithm.}

\paragraph{Nonlinear Least-Squares: Gauss vs.\ Newton.} In this
section we are interested in finding an approximate solution,
$\mathbf{\hat z}$, to the nonlinear inverse problem \[
\mathbf{g}(\mathbf{z}) \approx y
\] for known $y \in \mathbb{C}^m$ and  known real-analytic function
$\mathbf{g}: \mathbb{C}^n \rightarrow \mathbb{C}^m$. We desire a
least-squares solution, which is a solution that minimizes the
weighted least-squares loss function\footnote{The factor of
$\frac{1}{2}$ has been included for notational convenience in the
ensuing derivations. If it is removed, some of the intermediate
quantities derived subsequently (such as Hessians, etc.) will differ
by a factor of 2, although the ultimate answer is independent of any
overall constant factor of the loss function. If in your own problem
solving ventures, your intermediate quantities appear to be off by a
factor of 2 relative to the results given in this note, you should
check whether your loss function does or does not have this factor.}
\[ \ell(\mathbf{z}) = \frac{1}{2}\,  \left\| \mathbf{y} -
\mathbf{g}(\mathbf{z}) \right\|^2_W = \frac{1}{2}\, \left(\mathbf{y}
- \mathbf{g}(\mathbf{z})\right)^H W \left(\mathbf{y} -
\mathbf{g}(\mathbf{z})\right) \] where $W$ is a Hermitian
positive-definite weighting matrix. Although the nonlinear function
$\mathbf{g}$ is assumed to be real-analytic, in general it is
assumed to be \emph{not} holomorphic (i.e., $\mathbf{g}$ is
\emph{not complex-analytic in $\mathbf{z}$}).

\medskip
In the subsequent development we will analyze the problem using the
$\mathbf{c}$-real perspective developed in the preceding
discussions. Thus, the loss function is assumed to be re-expressible
in terms of $\mathbf{c}$,
\begin{equation}\label{eq:weightedlsloss} \ell(\mathbf{c}) = \frac{1}{2} \left\|
\mathbf{y} - \mathbf{g}(\mathbf{c}) \right\|^2_W = \frac{1}{2} \,
\left(\mathbf{y} - \mathbf{g}(\mathbf{c})\right)^H W
\left(\mathbf{y} - \mathbf{g}(\mathbf{c})\right)\, .
\end{equation}
Intermediate quantities produced from this perspective\footnote{Such as the
Gauss-Newton Hessian to be discussed below.} may have a different
functional form than those produced purely within the $\mathbf{z}
\in \mathcal{Z}$ perspective, but the end results will be the same.

\medskip
We will consider two iterative algorithms for minimizing the loss
function (\ref{eq:weightedlsloss}): The Newton algorithm, discussed
above, and the  Gauss-Newton algorithm  which is usually a somewhat
simpler, yet related, method for iteratively finding a solution
which minimizes a least-squares function of the form
(\ref{eq:weightedlsloss}).\footnote{The Newton algorithm is a
\emph{general method} that can be used to minimize a variety of
different loss functions, while the Gauss-Newton algorithm is a
\emph{least-squares estimation method} which is specific to the
problem of minimizing the least-squares loss function
(\ref{eq:weightedlsloss}).}

\medskip
As discussed earlier, the Newton method is based on an iterative
quadratic expansion and minimization of the loss function
$\ell(\mathbf{z})$ about a current solution estimation,
$\mathbf{\hat z}$. Specifically the Newton method minimizes an
approximation to $\ell (\mathbf{c}) = \ell(\mathbf{z})$ based on the
second order expansion of $\ell(\mathbf{c})$ in $\Delta \mathbf{c}$
about a current solution estimate $\mathbf{\hat c} =
\text{col}(\mathbf{\hat z}, \mathbf{\hat{\bar{z}}})$, \[
\ell(\mathbf{ \hat c} + \Delta \mathbf{c}) \approx \hat{\ell}(\Delta
\mathbf{c})^{\text{\tiny Newton}} \] where we define the \emph{Newton approximate loss function,}
\begin{equation}\label{eq:newtonloss} \hat{\ell}(\Delta
\mathbf{c})^{\text{\tiny Newton}} = \ell(\mathbf{\hat c}) +
\frac{\partial \ell(\mathbf{\hat c})}{\partial \mathbf{c}} \, \Delta
\mathbf{c} + \frac{1}{2} \Delta \mathbf{c}^H \,
\mathcal{H}_{\mathbf{c}\mathbf{c}}^{\text{\boldmath \tiny
$\mathbb{C}$}}(\mathbf{\hat c}) \, \Delta \mathbf{c} .
\end{equation}
Minimizing the Newton loss function $\hat{\ell}(\Delta
\mathbf{c})^{\text{\tiny Newton}}$ then results in a correction
$\widehat{\Delta {\mathbf{c}}}^{\text{\tiny Newton}}$ which is then
used to update the estimate $\mathbf{\hat c} \leftarrow \mathbf{\hat
c} + \alpha \widehat{\Delta {\mathbf{c}}}^{\text{\tiny Newton}}$ for
some stepsize $\alpha
> 0$. The algorithm then starts all over again. As mentioned above, a ``completing-the-square'' argument can be invoked to
readily show that the correction which minimizes the quadratic
Newton loss function is given by
\begin{equation}\label{eq:newtoncorrection} \widehat{\Delta
{\mathbf{c}}}^{\text{\tiny Newton}} = -
\mathcal{H}_{\mathbf{c}\mathbf{c}}^{\text{\boldmath \tiny
$\mathbb{C}$}}(\mathbf{\hat c})^{-1} \left( \frac{\partial
\ell(\mathbf{\hat c})}{\partial \mathbf{c}} \right)^H
\end{equation}
provided that the $\mathbf{c}$-complex Hessian
$\mathcal{H}_{\mathbf{c}\mathbf{c}}^{\text{\boldmath \tiny
$\mathbb{C}$}}(\mathbf{\hat c})$ is invertible. Because it defines
the second-order term in the Newton loss function and directly
enters into the Newton correction, we will often refer to
$\mathcal{H}_{\mathbf{c}\mathbf{c}}^{\text{\boldmath \tiny
$\mathbb{C}$}}(\mathbf{\hat c})$ as the \emph{\underline{Newton Hessian}.} If we
block partition the Newton Hessian and solve for the correction
$\widehat{\Delta \mathbf{z}}^{\text{\tiny Newton}}$, we obtain the
solution (\ref{eq:zoptimal2}) which we earlier derived for a more
general (possibly non-quadratic) loss function.

\medskip
We now  determine the form of the cogradient $\frac{\partial
\ell(\mathbf{\hat c})}{\partial \mathbf{c}}$ of the least-squares
loss function (\ref{eq:weightedlsloss}). This is done by utilizing
the $\mathbf{c}$-real perspective which allows us to take (real)
cogradients with respect to $\mathbf{c}$-real. First, however, it is
convenient to \emph{define the compound Jacobian $G(\mathbf{\hat c})$ of
$\mathbf{g}(\mathbf{\hat c})$} as
\begin{equation}\label{eq:compoundjac}
G(\mathbf{\hat c}) \triangleq \frac{\partial \mathbf{g}(\mathbf{\hat
c})}{\partial \mathbf{c}} \triangleq \begin{pmatrix} \frac{\partial
\mathbf{g}(\mathbf{\hat z})}{\partial \mathbf{z}} & \frac{\partial
\mathbf{g}(\mathbf{\hat z})}{\partial \mathbf{\bar z}}\end{pmatrix}
=
\begin{pmatrix} J_\mathbf{g}(\mathbf{c}) & J_\mathbf{g}^c(\mathbf{c}) \end{pmatrix} \in \mathbb{C}^{m \times 2n} \, .
\end{equation}
Setting $\mathbf{e}  = \mathbf{y} - \mathbf{g}(\mathbf{c})$, we
have\footnote{Remember that $\frac{\partial }{\partial \mathbf{c}}$
is only well-defined as a derivative within the $\mathbf{c}$-real
framework.}
\begin{eqnarray*}
\frac{\partial \ell}{\partial \mathbf{c}} & = & \frac{1}{2}
\frac{\partial }{\partial \mathbf{c}} \,  \mathbf{e}^H W
\mathbf{e} \\
& = & \frac{1}{2}\, \mathbf{e}^H W \frac{\partial }{\partial
\mathbf{c}} \, \mathbf{e} + \frac{1}{2} \, \mathbf{e}^T W^T
\frac{\partial }{\partial \mathbf{c}} \, \mathbf{\bar e} \\
& = & - \frac{1}{2}\, \mathbf{e}^H W\,  \frac{\partial
{\mathbf{g}}}{\partial \mathbf{c}} \,  - \frac{1}{2} \, \mathbf{e}^T
W^T \, \frac{\partial \bar{\mathbf{g}}}{\partial
\mathbf{c}} \\
& = & - \frac{1}{2}\, \mathbf{e}^H W\,  G \,  - \frac{1}{2} \,
\mathbf{e}^T W^T \, \overline{\left( \frac{\partial
{\mathbf{g}}}{\partial \mathbf{c}} S \right)} \\
& = & - \frac{1}{2}\, \mathbf{e}^H W\,  G \,  - \frac{1}{2} \,
\mathbf{e}^T W^T \, \overline{G} S
\end{eqnarray*}
or
\begin{equation}\label{eq:ellcogradient}
\frac{\partial \ell}{\partial \mathbf{c}} = - \frac{1}{2}\,
\mathbf{e}^H W\,  G \,  - \frac{1}{2} \, \overline{\mathbf{e}^H W \,
G} S.
\end{equation}
This expression for $\frac{\partial \ell}{\partial \mathbf{c}}$ is
admissible, as required, as it is readily verified that
\[ \overline{\left(\frac{\partial \ell}{\partial \mathbf{c}}\right)^H} = S \left(\frac{\partial \ell}{\partial
\mathbf{c}}\right)^H \] as per the requirement given in
(\ref{eq:vecttest}).

\medskip
 The linear term in
the Newton loss function $\hat{\ell}^{\text{\tiny Newton}}$ is
therefore given by
\begin{eqnarray*}
\frac{\partial \ell}{\partial \mathbf{c}} \Delta \mathbf{c} & = &
 - \frac{1}{2}\,
\mathbf{e}^H W\,  G \, \Delta \mathbf{c}  - \frac{1}{2} \,
\overline{\mathbf{e}^H W \, G} S \, \Delta \mathbf{c} \\
& = &
 - \frac{1}{2}\,
\mathbf{e}^H W\,  G \, \Delta \mathbf{c}  - \frac{1}{2} \,
\overline{\, \mathbf{e}^H W\,  G \, \Delta \mathbf{c}\, } \\ & = &
 - \text{Re} \,\left\{
\mathbf{e}^H W\,  G \, \Delta \mathbf{c} \right\} .
\end{eqnarray*}
Thus
\begin{equation} \label{eq:lininnewton}
\frac{\partial \ell}{\partial \mathbf{c}} \Delta \mathbf{c}=  -
\text{Re} \,\left\{ \mathbf{e}^H W\,  G \, \Delta \mathbf{c}
\right\} =
 - \text{Re} \, \left\{
\left(\mathbf{y} - \mathbf{g}(\mathbf{c})\right)^H W\,  G \, \Delta
\mathbf{c} \right\}.
\end{equation}
If the reader has any doubts as to the validity or correctness of
this derivation, she/he is invited to show that the right-hand side
of (\ref{eq:lininnewton}) is equal to $2\, \text{Re} \left\{
\frac{\partial \ell}{\partial \mathbf{z}} \, \Delta \mathbf{z}
\right\}$ as required from equation (\ref{eq:rcz1terms}).

\medskip
Before continuing on to determine  the functional form of the
$\mathbf{c}$-complex Hessian
$\mathcal{H}_{\mathbf{c}\mathbf{c}}^{\text{\boldmath \tiny
$\mathbb{C}$}}(\mathbf{\hat c})$ needed to form the Newton loss
function and solution, we turn first to a discussion of the
Gauss-Newton algorithm.

\medskip
Whereas the Newton method is based on an iterative quadratic
expansion and minimization of the loss function $\ell(\mathbf{z})$
about a current solution estimation, $\mathbf{\hat z}$, The
Gauss-Newton method is based on iterative ``relinearization'' of the
system equations $\mathbf{y} \approx \mathbf{g}(\mathbf{z})$ about
the current estimate, $\mathbf{\hat z}$ and minimization of the
resulting approximate least-squares problem. \emph{We put
``linearization'' in quotes because (unless the function
$\mathbf{g}$ happens to be holomorphic) generally we are not
linearizing $\mathbf{g}$ with respect to $z$ but, rather, we are
linearizing with respect to $\mathbf{c} =
\text{col}(\mathbf{z},\mathbf{\bar z})$.}

\medskip
Expanding the system equations $\mathbf{y} \approx
\mathbf{g}(\mathbf{z})$ about a current estimate $\mathbf{\hat z}$,
we have
\[  \mathbf{y} - \mathbf{g}(\mathbf{z})= \mathbf{y} - \mathbf{g}(\mathbf{\hat z} + \Delta \mathbf{z})
 \approx \mathbf{y} - \left( \mathbf{g}(\mathbf{\hat
z}) + \frac{\partial \mathbf{g}(\mathbf{\hat z})}{\partial
\mathbf{z}} \Delta \mathbf{z} + \frac{\partial
\mathbf{g}(\mathbf{\hat z})}{\partial \mathbf{\bar z}} \Delta
\mathbf{\bar{z}} \right)
\]
where $\Delta \mathbf{z} = \mathbf{z} - \mathbf{\hat z}$ and $\Delta
\mathbf{\bar z} = \overline{\Delta \mathbf{z}} = \mathbf{\bar z} -
\mathbf{\bar{\hat z}} = \mathbf{\bar z} - \mathbf{\hat{\bar z}}$.
Note that the approximation to $\mathbf{g}$ is \emph{not} a linear function
of $\mathbf{z}$ as complex conjugation is a nonlinear operation on
$\mathbf{z}$. However, if $\mathbf{g}$ is holomorphic, then
$\frac{\partial \mathbf{g} }{\partial \mathbf{\bar z}} \equiv 0$, in
which case the approximation is linear in $\mathbf{z}$. Although the
approximation of $\mathbf{g}$ generally is not linear in
$\mathbf{z}$, it is linear in $\mathbf{c} =
\text{col}(\mathbf{z},\mathbf{\bar z})$, and we rewrite the
approximation as
\begin{equation}\label{eq:gausslin}
 \mathbf{y} - \mathbf{g}(\mathbf{c}) = \mathbf{y} - \mathbf{g}(\mathbf{\hat c} + \Delta \mathbf{c}) \approx
\Delta \mathbf{y} - G(\mathbf{\hat c})\,  \Delta \mathbf{c}
\end{equation}
where $\Delta \mathbf{y} = \mathbf{y} - \mathbf{g}(\mathbf{\hat
z})$, $\mathbf{\hat c} = \text{col}(\mathbf{\hat
z},\mathbf{\hat{\bar z}})$, $\Delta \mathbf{c} = \mathbf{c} -
\mathbf{\hat c}$, and $G(\mathbf{\hat c})$ is the (compound)
Jacobian mapping of $\mathbf{g}$ evaluated at the current estimate
$\mathbf{\hat c}$ given in Equation (\ref{eq:compoundjac}). With
this approximation, the loss function (\ref{eq:weightedlsloss}) is
approximated by the following quadratic loss function (notationally
suppressing the dependence on $\mathbf{\hat c}$),
\[ \ell(\mathbf{c}) = \ell(\mathbf{ \hat c} + \Delta \mathbf{c}) \approx \hat \ell(\Delta
\mathbf{c})^{\text{\tiny Gauss}} \] where
\begin{eqnarray*}
 \hat \ell(\Delta
\mathbf{c})^{\text{\tiny Gauss}} & = & \frac{1}{2} \left\| \Delta
\mathbf{y} -
G \,  \Delta \mathbf{c} \right\|^2_W \\
& = & \frac{1}{2} \left( \Delta \mathbf{y} -  G \, \Delta \mathbf{c}
\right)^H W
\left( \Delta \mathbf{y} - G \,  \Delta \mathbf{c} \right) \\
& = & \frac{1}{2} \| \Delta \mathbf{y}\|^2 - \text{Re} \,\left\{
\Delta \mathbf{y}^H W\,  G \, \Delta \mathbf{c} \right\} +
\frac{1}{2} \Delta \mathbf{c}^H\, G^H W G \, \Delta \mathbf{c} \\
& = & \ell(\mathbf{\hat c}) + \frac{\partial \ell(\mathbf{\hat
c})}{\partial \mathbf{c}} \, \Delta \mathbf{c} + \frac{1}{2} \Delta
\mathbf{c}^H\, G^H W G \, \Delta \mathbf{c}. \qquad \text{(from
(\ref{eq:lininnewton})}
\end{eqnarray*}

Unfortunately, the resulting quadratic form
\begin{equation}\label{eq:gaussloss1}
\hat \ell(\Delta \mathbf{c})^{\text{\tiny Gauss}} =
\ell(\mathbf{\hat c}) + \frac{\partial \ell(\mathbf{\hat
c})}{\partial \mathbf{c}} \, \Delta \mathbf{c} + \frac{1}{2} \Delta
\mathbf{c}^H\, G^H W G \, \Delta \mathbf{c}
\end{equation}
is {\emph{not}} \emph{admissible} as it stands.\footnote{And thus the complex
Gauss-Newton algorithm is generally more complicated in form than the real
Gauss-Newton algorithm for which the quadratic form
(\ref{eq:gaussloss1}) is meaningful.} This is because
the matrix $G^HWG$ is not admissible,
\[ G^H W G = \left(\frac{\partial \mathbf{g}}{\partial \mathbf{c}}
\right)^H W \left(\frac{\partial \mathbf{g}}{\partial \mathbf{c}}
\right) \notin \mathcal{L}(\mathcal{C},\mathcal{C}) . \] This can be
seen by showing that the condition (\ref{eq:ccmappingcond2}) is
violated:
\begin{eqnarray*}
S \, \overline{G^H W G}\,  S & = & S \,
\overline{\left(\frac{\partial \mathbf{g}}{\partial \mathbf{c}}
\right)^H W \left(\frac{\partial \mathbf{g}}{\partial \mathbf{c}}
\right)}\, S \\
& = &  \overline{\left(\frac{\partial \mathbf{g}}{\partial
\mathbf{\bar c}} \right)^H W \left(\frac{\partial
\mathbf{g}}{\partial \mathbf{\bar c}}
\right)} \\
& = & {\left(\frac{\partial \mathbf{\bar g}}{\partial \mathbf{ c}}
\right)^H \bar{W} \left(\frac{\partial \mathbf{\bar g}}{\partial
\mathbf{c}} \right)}
\\
& \ne & {\left(\frac{\partial \mathbf{g}}{\partial \mathbf{ c}}
\right)^H {W} \left(\frac{\partial \mathbf{g}}{\partial \mathbf{c}}
\right)}.
\end{eqnarray*}

\medskip
Fortunately, we can rewrite  the quadratic form
(\ref{eq:gaussloss1}) as an equivalent form which is admissible on
$\mathcal{C}$. To do this note that $G^H W G$ is Hermitian, so that
\[ \Delta \mathbf{c}^H G^HW G \Delta \mathbf{c} = \overline{\Delta \mathbf{c}^H G^HW G \Delta \mathbf{c}}
\in \mathbb{R} \, . \] Also recall from Equation (\ref{eq:Cproj})
that $\mathbf{P}(G^HWG) \in \mathcal{L}(\mathcal{C},\mathcal{C})$
and $\Delta \mathbf{c} \in \mathcal{C} \Rightarrow S \Delta
\mathbf{c} = \Delta \mathbf{\bar c}$. For an admissible variation $\Delta \mathbf{c} \in \mathcal{C}$ we have\footnote{Note that the
ensuing derivation does \emph{not} imply that $G^H W G =
\mathbf{P}({G^H W G})$, a fact which would contradict our claim that
$G^H WG$ is not admissible. This is because in the derivation we are
\emph{not} allowing \emph{arbitrary vectors} in $\mathbb{C}^{2n}$ but are only
admitting vectors $\Delta \mathbf{c}$ constrained to lie in
$\mathcal{C}$, $\Delta \mathbf{c} \in \mathcal{C} \subset
\mathbb{C}^{2n}$.}
\begin{eqnarray*}
\Delta \mathbf{c}^H G^H W G \Delta \mathbf{c} & = & \Delta
\mathbf{c}^H \mathbf{P}(G^H W G) \Delta \mathbf{c} + \Delta
\mathbf{c}^H \left( G^H W G - \mathbf{P}(G^H W G) \right) \Delta
\mathbf{c} \\
& = & \Delta \mathbf{c}^H \mathbf{P}(G^H W G) \Delta \mathbf{c} +
\frac{1}{2} \Delta \mathbf{c}^H \left( G^H W G - S \overline{G^H W
G} S \right) \Delta \mathbf{c} \\& = & \Delta \mathbf{c}^H
\mathbf{P}(G^H W G) \Delta \mathbf{c} + \frac{1}{2}\left( \Delta
\mathbf{c}^H G^H W G \Delta \mathbf{c} - \overline{\Delta
\mathbf{c}^H G^H W G \Delta \mathbf{c}} \right) \\
& = & \Delta \mathbf{c}^H \mathbf{P}(G^H W G) \Delta \mathbf{c} + 0
\\
& = & \Delta \mathbf{c}^H \mathbf{P}(G^H W G) \Delta \mathbf{c} \, .
\end{eqnarray*}

Thus we have shown that on the space of admissible variations,
$\Delta \mathbf{c} \in \mathcal{C}$, the inadmissible quadratic form
(\ref{eq:gaussloss1}) is equivalent to the admissible quadratic form (the \emph{Gauss-Newton approximate loss function})
\begin{equation}\label{eq:gaussloss2}
\hat \ell(\Delta \mathbf{c})^{\text{\tiny Gauss}} =
\ell(\mathbf{\hat c}) + \frac{\partial \ell(\mathbf{\hat
c})}{\partial \mathbf{c}} \, \Delta \mathbf{c} + \frac{1}{2} \Delta
\mathbf{c}^H\, \mathcal{H}_{\mathbf{c}\mathbf{c}}^{\text{\tiny
Gauss}}(\mathbf{\hat c}) \, \Delta \mathbf{c}
\end{equation}
where
\begin{equation}\label{eq:gaussnewtonhessian}
\boxed{\ \textsf{Gauss-Newton Hessian} \qquad \mathcal{H}_{\mathbf{c}\mathbf{c}}^{\text{\tiny Gauss}}(\mathbf{\hat
c}) \triangleq \mathbf{P}\left(G^H(\mathbf{\hat c}) W G
(\mathbf{\hat c})\right) \ }
\end{equation}
denotes the \emph{Gauss-Newton} Hessian. Note that the Gauss-Newton Hessian is the exact Hessian matrix of the Gauss-Newton approximate loss function.

\medskip
Note that the Gauss-Newton Hessian
$\mathcal{H}_{\mathbf{c}\mathbf{c}}^{\text{\tiny
Gauss}}(\mathbf{\hat c})$ is Hermitian and always guaranteed to be
at least positive semi-definite, and guaranteed to be positive
definite if $\mathbf{g}$ is assumed to be one-to-one (and thereby
ensuring that the compound Jacobian matrix $G$ has full column
rank). This is in contrast to the Newton (i.e., the
$\mathbf{c}$-complex) Hessian
$\mathcal{H}_{\mathbf{c}\mathbf{c}}^{\text{\boldmath \tiny
$\mathbb{C}$}}(\mathbf{\hat c})$ which, unfortunately, can be
indefinite or rank deficient even though it is Hermitian and even if
$\mathbf{g}$ is one-to-one.

\medskip
 Assuming that
$\mathcal{H}_{\mathbf{c}\mathbf{c}}^{\text{\tiny
Gauss}}(\mathbf{\hat c})$ is invertible, the correction which
minimizes the Gauss-Newton approximate loss function (\ref{eq:gaussloss2}) is
given by
\begin{equation}\label{eq:gausscorrection} \widehat{\Delta
{\mathbf{c}}}^{\text{\tiny Gauss}} = -
\mathcal{H}_{\mathbf{c}\mathbf{c}}^{\text{\tiny Gauss}}(\mathbf{\hat
c})^{-1} \left( \frac{\partial \ell(\mathbf{\hat c})}{\partial
\mathbf{c}} \right)^H.
\end{equation}
Because of the admissibility of
$\mathcal{H}_{\mathbf{c}\mathbf{c}}^{\text{\tiny Gauss}}$ and
$\left( \frac{\partial \ell(\mathbf{\hat c})}{\partial \mathbf{c}}
\right)^H$, the resulting solution is admissible $\widehat{\Delta
{\mathbf{c}}}^{\text{\tiny Gauss}} \in \mathcal{C}$.

\medskip

Comparing Equations (\ref{eq:newtoncorrection}) and
(\ref{eq:gausscorrection}), it is evident that the difference
between the two algorithms resides in the difference between the
Newton Hessian, $\mathcal{H}_{\mathbf{c}\mathbf{c}}^{\text{\boldmath
\tiny $\mathbb{C}$}}(\mathbf{\hat c})$, which is the actual
$\mathbf{c}$-complex Hessian of the least-squares loss function
$\ell(\mathbf{c})$, and the Gauss-Newton Hessian
$\mathcal{H}_{\mathbf{c}\mathbf{c}}^{\text{\tiny
Gauss}}(\mathbf{\hat c})$ which has an as yet unclear relationship to
$\ell(\mathbf{c})$.\footnote{Note that, \emph{by construction,}
$\mathcal{H}_{\mathbf{c}\mathbf{c}}^{\text{\tiny
Gauss}}(\mathbf{\hat c})$ is the Hessian matrix of the Gauss-Newton
approximate loss function. The question is: what is its relationship to the
least-squares loss function or the Newton approximate loss function?} For this
reason, we now turn to a discussion of the relationship between the Hessians
$\mathcal{H}_{\mathbf{c}\mathbf{c}}^{\text{\boldmath \tiny
$\mathbb{C}$}}(\mathbf{\hat c})$ and
$\mathcal{H}_{\mathbf{c}\mathbf{c}}^{\text{\tiny
Gauss}}(\mathbf{\hat c})$.

\medskip
We can compute the Newton Hessian $\mathcal{H}_{\mathbf{c}
\mathbf{c}}^{\text{\boldmath \tiny $\mathbb{C}$}}$ from the
relationship (see Equation (\ref{eq:hrhc2}))
\[ \mathcal{H}_{\mathbf{c} \mathbf{c}}^{\text{\boldmath \tiny
$\mathbb{C}$}} = S \, \mathcal{H}_{\mathbf{c}
\mathbf{c}}^{\text{\boldmath \tiny $\mathbb{R}$}} = S\,
\frac{\partial}{\partial \mathbf{c}} \left( \frac{\partial
\ell}{\partial \mathbf{c}} \right)^T \] where $\frac{\partial \
}{\partial \mathbf{c}}$ is taken to be a $\mathbf{c}$-real
cogradient operator. Note from (\ref{eq:ellcogradient}) that,
\begin{equation}\label{eq:quadgrad} \left( \frac{\partial \ell}{\partial \mathbf{c}} \right)^H
= - \frac{1}{2} {G^H W \mathbf{e}} - \frac{1}{2} S\overline{G^H W
\mathbf{e}}  = \frac{1}{2} \left( {B} + S \overline{B} \right)
,\end{equation} where
\begin{equation}\label{eq:bdef}
B \triangleq - G^H W \mathbf{e}
\end{equation}
with $\mathbf{e} = \mathbf{y} - \mathbf{g}(\mathbf{c})$. This
results in
\[
\left( \frac{\partial \ell}{\partial \mathbf{c}} \right)^T =
\overline{\left( \frac{\partial \ell}{\partial \mathbf{c}}
\right)^H} =   \frac{1}{2} \left( \bar{B} + S B \right),
\]
Also note that
\[ \frac{\partial \bar{B}}{\partial \mathbf{c}} = \overline{\, \frac{\partial {B}}{\partial
\mathbf{\bar c}} \, }  = \overline{\left( \frac{\partial
{B}}{\partial \mathbf{c}} \, S \right) } = \overline{\,
\frac{\partial {B}}{\partial \mathbf{c}} \, } \, S.\]

\medskip
We have
\[ \mathcal{H}_{\mathbf{c}
\mathbf{c}}^{\text{\boldmath \tiny $\mathbb{R}$}} =
\frac{\partial}{\partial \mathbf{c}} \left( \frac{\partial
\ell}{\partial \mathbf{c}} \right)^T  = \frac{1}{2} \left( S \,
\frac{\partial {B}}{\partial \mathbf{c}} +\frac{\partial
\bar{B}}{\partial \mathbf{c}}   \right) \] or
\begin{equation}\label{eq:hrccls}
\mathcal{H}_{\mathbf{c} \mathbf{c}}^{\text{\boldmath \tiny
$\mathbb{R}$}} = \frac{\partial}{\partial \mathbf{c}} \left(
\frac{\partial \ell}{\partial \mathbf{c}} \right)^T  = \frac{1}{2}
\left( S \, \frac{\partial {B}}{\partial \mathbf{c}} + \overline{\,
\frac{\partial {B}}{\partial \mathbf{c}}\, } \, S   \right).
\end{equation}
This yields
\begin{equation}\label{eq:hrccls1} \mathcal{H}_{\mathbf{c} \mathbf{c}}^{\text{\boldmath \tiny
$\mathbb{C}$}} = S \, \mathcal{H}_{\mathbf{c}
\mathbf{c}}^{\text{\boldmath \tiny $\mathbb{R}$}} = \frac{1}{2}
\left(  \frac{\partial {B}}{\partial \mathbf{c}} + S \, \overline{\,
\frac{\partial {B}}{\partial \mathbf{c}}\, } \, S     \right)
\end{equation}
with $B$ given by (\ref{eq:bdef}), which we can write as
\begin{equation}\label{eq:symmetrizing} \mathcal{H}_{\mathbf{c} \mathbf{c}}^{\text{\boldmath \tiny
$\mathbb{C}$}} = S \, \mathcal{H}_{\mathbf{c}
\mathbf{c}}^{\text{\boldmath \tiny $\mathbb{R}$}} = \mathbf{P}\left(
\frac{\partial {B}}{\partial \mathbf{c}}\right) \, .
\end{equation}

Recall that $\mathcal{H}_{\mathbf{c} \mathbf{c}}^{\text{\boldmath
\tiny $\mathbb{C}$}}$ must be admissible. The function $\text{\bf
P}(\cdot)$ produces admissible matrices which map from $\mathcal{C}$
to $\mathcal{C}$, and thereby ensures that the right-hand side of
equation (\ref{eq:symmetrizing}) is indeed an admissible matrix, as
required for self-consistency of our development. The presence of
the
 operator $\mathbf P$ does not show up in the real case (which is
the standard development given in textbooks) as $\frac{\partial
{B}}{\partial \mathbf{c}}$ is automatically symmetric as required
for admissibility  in the real case.

\medskip
Note that $B$ can be written as
\[ B = - \left(\frac{\partial \mathbf{g}}{\partial \mathbf{c}}\right)^H
W \left(\mathbf{y} - \mathbf{g}\right) = - \sum_{i=1}^m
 \left( \frac{\partial
g_i}{\partial \mathbf{c}}\right)^H \left[ W \left(\mathbf{y}
-\mathbf{g} \right)\, \right]_i \] where $g_i$ and $\left[ W
\left(\mathbf{y} -\mathbf{g} \right)\,  \right]_i$ denote the $i$-th
(scalar) components of the vectors $\mathbf{g}$ and $W\mathbf{e} = W(\mathbf{y}-
\mathbf{g})$ respectively. We can then compute $\frac{\partial
B}{\partial \mathbf{c}}$ as
\begin{eqnarray*}
\frac{\partial B}{\partial \mathbf{c}} & = &  \left(\frac{\partial
\mathbf{g}}{\partial \mathbf{c}}\right)^H W \left(\frac{\partial
\mathbf{g}}{\partial \mathbf{c}}\right) - \sum_{i=1}^m
\frac{\partial}{\partial \mathbf{c}} \left( \frac{\partial
g_i}{\partial \mathbf{c}}\right)^H \left[ W \left(\mathbf{y}
-\mathbf{g} \right) \, \right]_i \\ & = &  G^H W G - \sum_{i=1}^m
\frac{\partial}{\partial \mathbf{c}} \left( \frac{\partial
g_i}{\partial \mathbf{c}}\right)^H \left[ W \left(\mathbf{y}
-\mathbf{g} \right) \, \right]_i
\end{eqnarray*}
or
\begin{equation}\label{eq:bderiv}
\frac{\partial B}{\partial \mathbf{c}} = G^H W G - \sum_{i=1}^m
\frac{\partial}{\partial \mathbf{c}} \left( \frac{\partial
g_i}{\partial \mathbf{c}}\right)^H  \left[ W \mathbf{e} \, \right]_i
\, .
\end{equation}

Equations (\ref{eq:symmetrizing}) and (\ref{eq:bderiv}) result in the following succinct relationship between the complex Newton and Gauss-Newton Hessians,
\begin{equation}\label{eq:syminhh}
\boxed{\ \textsf{Newton Hessian} \qquad \mathcal{H}_{\mathbf{c}\mathbf{c}}^{\text{\tiny Newton}} = \mathcal{H}_{\mathbf{c} \mathbf{c}}^{\text{\boldmath \tiny
$\mathbb{C}$}} = \mathcal{H}_{\mathbf{c}\mathbf{c}}^{\text{\tiny
Gauss}}  - \sum_{i=1}^m
\mathcal{H}_{\mathbf{c}\mathbf{c}}^{\text{\tiny $(i)$}}  \ }
\end{equation}
where the Gauss-Newton Hessian $\mathcal{H}_{\mathbf{c}\mathbf{c}}^{\text{\tiny
Gauss}}$ is given by (\ref{eq:gaussnewtonhessian}) and 
\begin{equation}\label{eq:pmatrices}
\mathcal{H}_{\mathbf{c}\mathbf{c}}^{\text{\tiny $(i)$}} \triangleq
\mathbf{P}\left(\frac{\partial}{\partial \mathbf{c}} \left(
\frac{\partial g_i}{\partial \mathbf{c}}\right)^H \, \left[ W
\mathbf{e} \, \right]_i \right)
 \, , \quad
i = 1, \cdots, m \, .
\end{equation}
Equation (\ref{eq:syminhh}), which is our final result for
the structural form of the Newton Hessian $\mathcal{H}_{\mathbf{c}
\mathbf{c}}^{\text{\boldmath \tiny $\mathbb{C}$}}$, looks very much
like the equivalent result for the real case.\footnote{The
primary difference is due to the presence of the projector
$\mathbf{P}$ in the complex Newton algorithm. Despite the
similarity, note that it takes much more work to rigorously derive
the complex Newton-Algorithm!} The first term on the right-hand-side
of (\ref{eq:syminhh}) is the Gauss-Newton Hessian
$\mathcal{H}_{\mathbf{c}\mathbf{c}}^{\text{\tiny Gauss}}$, which is
admissible, Hermitian and at least positive semidefinite (under the
standard assumption that $W$ is Hermitian positive definite). Below,
we will show that the matrices
$\mathcal{H}_{\mathbf{c}\mathbf{c}}^{\text{\tiny $(i)$}}$, $i = 1,
\cdots, m$, are all \emph{individually} admissible and Hermitian.\footnote{Of course, because $\mathcal{H}_{\mathbf{c} \mathbf{c}}^{\text{\boldmath \tiny
$\mathbb{C}$}}$ and $\mathcal{H}_{\mathbf{c}\mathbf{c}}^{\text{\tiny
Gauss}}$ are both Hermitian and admissible the \emph{total sum}   $\sum_{i=1}^m
\mathcal{H}_{\mathbf{c}\mathbf{c}}^{\text{\tiny $(i)$}}$ must be Hermitian and admissible.}  While the
Gauss-Newton Hessian is always positive semidefinite (and always
positive definite if $\mathbf{g}$ is one-to-one), the presence of
the second term on the right-hand-side of (\ref{eq:syminhh}) can
cause the Newton Hessian to become indefinite, or even negative
definite.

\medskip
We can now understand the relationship between the Gauss-Newton
method and the Newton method when applied to the problem of minimizing
the least-squares loss function. \emph{The Gauss-Newton method is an
approximation to the Newton method which arises from ignoring the
second term on the right-hand-side of (\ref{eq:syminhh}).} This
approximation is not only easier to implement, it will generally
have superior numerical properties as a consequence of the
definiteness of the Gauss-Newton Hessian. Indeed, if the mapping
$\mathbf{g}$ is onto, via the Gauss-Newton algorithm one can produce
a sequence of estimates $\mathbf{\hat c}_k$, $k = 1, 2, 3, \cdots$,
which drives $\mathbf{e(\mathbf{\hat c}_k)} = \mathbf{y} -
\mathbf{g}(\mathbf{\hat c}_k)$, and hence (with some additional smoothness assumptions on $\mathbf{g}$) the second term on the
right-hand-side of (\ref{eq:syminhh}), to zero as $k \rightarrow
\infty$. In which case, asymptotically there will be little
difference in the convergence properties between the Newton and
Gauss-Newton methods.  This property is well known in the classical
optimization literature, which suggests that by working within the
$\mathbf{c}$-real perspective, we may be able to utilize a variety
of insights that have been developed for the Newton and Gauss-Newton
methods when optimizing over real vector spaces.

\medskip
We will now
demonstrate that each \emph{individual} term 
$\mathcal{H}_{\mathbf{c}\mathbf{c}}^{\text{\tiny $(i)$}}$, $i = 1,
\cdots, m$, in (\ref{eq:syminhh}) is admissible and Hermitian. Note that the ``raw''
matrix \[ \left[ W \mathbf{e} \, \right]_i \,
\frac{\partial}{\partial \mathbf{c}} \left( \frac{\partial
g_i}{\partial \mathbf{c}}\right)^H \] is neither Hermitian nor
admissible because of the presence of the complex scalar factor
$\left[ W \mathbf{e} \, \right]_i$. Fortunately, the processing of
the second matrix of partial derivatives by the operator $\mathbf P$
to form the matrix $\mathcal{H}_{\mathbf{c}\mathbf{c}}^{\text{\tiny
$(i)$}}$ via
\[ \mathcal{H}_{\mathbf{c}\mathbf{c}}^{\text{\tiny $(i)$}} =
\mathbf{P}\left( \left[ W \mathbf{e} \, \right]_i \,
\frac{\partial}{\partial \mathbf{c}} \left( \frac{\partial
g_i}{\partial \mathbf{c}}\right)^H  \right) \] creates a matrix which
is both admissible and Hermitian. The fact that
$\mathcal{H}_{\mathbf{c}\mathbf{c}}^{\text{\tiny $(i)$}}$ is
admissible is obvious, as the projector $\mathbf{P}$ is idempotent.
We will now prove that
$\mathcal{H}_{\mathbf{c}\mathbf{c}}^{\text{\tiny $(i)$}}$ is
Hermitian.

\medskip
 Define the matrix
\begin{equation}\label{eq:rawhessian}
A_{\mathbf{c}\mathbf{c}}(g_i) \triangleq  \frac{\partial}{\partial
\mathbf{c}} \left( \frac{\partial g_i}{\partial \mathbf{c}}\right)^H
, \end{equation} and note that
\[ \left[ \frac{\partial}{\partial \mathbf{c}} \left(
\frac{\partial g_i}{\partial \mathbf{c}}\right)^H \right]^H =
\overline{\left[ \frac{\partial}{\partial \mathbf{c}} \left(
\frac{\partial \bar{g}_i}{\partial \mathbf{\bar c}}\right)^T
\right]^T} = \overline{\left[ \frac{\partial}{\partial \mathbf{\bar
c}} \left( \frac{\partial \bar{g}_i}{\partial \mathbf{c}}\right)^T
\right]} = \frac{\partial}{\partial \mathbf{c}} \left(
\frac{\partial \bar{g}_i}{\partial \mathbf{c}}\right)^H , \] which
shows that $A_{\mathbf{c}\mathbf{c}}(g_i)$ has the property that
\begin{equation}\label{eq:sprop1}
A_{\mathbf{c}\mathbf{c}}(g_i)^H =
A_{\mathbf{c}\mathbf{c}}(\bar{g}_i) \, .
\end{equation}
Now note that
\begin{eqnarray*}
S \frac{\partial}{\partial \mathbf{c}} \left( \frac{\partial
g_i}{\partial \mathbf{c}}\right)^H S & = & S
\frac{\partial}{\partial \mathbf{\bar c}} \left( \frac{\partial
g_i}{\partial \mathbf{c}}\right)^H \\
& = & \frac{\partial}{\partial \mathbf{\bar c}} \left[ S \left(
\frac{\partial
g_i}{\partial \mathbf{c}}\right)^H \right]\\
& = & \frac{\partial}{\partial \mathbf{\bar c}} \left(
\frac{\partial g_i}{\partial \mathbf{c}} S \right)^H \\
 & = & \frac{\partial}{\partial \mathbf{\bar c}} \left( \frac{\partial
g_i}{\partial \mathbf{\bar c}}\right)^H \, ,
\end{eqnarray*}
which establishes the second property that
\begin{equation}\label{eq:sprop2}
S A_{\mathbf{c}\mathbf{c}}(g_i) S = A_{\mathbf{\bar c}\mathbf{\bar
c}}(g_i)\, .
\end{equation}
Finally note that properties (\ref{eq:sprop1}) and (\ref{eq:sprop2})
together yield the property
\[ A_{\mathbf{c}\mathbf{c}}(g_i)^H =
A_{\mathbf{c}\mathbf{c}}(\bar{g}_i) = SA_{\mathbf{\bar
c}\mathbf{\bar c}}(\bar g_i) S = S
\overline{A_{\mathbf{c}\mathbf{c}}(g_i)} S \, .
\]

Setting $a_i = \left[ W \mathbf{e} \, \right]_i$, we have
{\footnotesize
\[ \mathcal{H}_{\mathbf{c}\mathbf{c}}^{\text{\tiny
$(i)$}} = \mathbf{P}(a_i \, A_{\mathbf{c}\mathbf{c}}(g_i)) =
\frac{a_i \, A_{\mathbf{c}\mathbf{c}}(g_i) + S \, \overline{a_i \,
A_{\mathbf{c}\mathbf{c}}(g_i)} \, S}{2}  = \frac{a_i \,
A_{\mathbf{c}\mathbf{c}}(g_i) + \bar{a}_i \, S \, \overline{
A_{\mathbf{c}\mathbf{c}}(g_i)} \, S}{2} = \frac{a_i \,
A_{\mathbf{c}\mathbf{c}}(g_i) + \bar{a}_i \,
A_{\mathbf{c}\mathbf{c}}(g_i)^H}{2}
\]} \hspace{-.1in} which is obviously Hermitian. Note that the action of the
projector $\mathbf{P}$ on  ``raw'' matrix \hfill $a_i\, 
A_{\mathbf{c}\mathbf{c}}(g_i)$, Hermitian
symmetrizes the matrix $a_i \, A_{\mathbf{c}\mathbf{c}}(g_i)$.

\medskip
Below, we will examine the least-squares algorithms at the
block-component level, and will show that significant
simplifications occur when $\mathbf{g}(\mathbf{z})$ is holomorphic.

\paragraph{Generalized Gradient Descent Algorithms.}
As in the real case, the Newton and Gauss-Newton
algorithms can be viewed as special instances of a family of
generalized gradient descent algorithms. Given a general real-valued
loss function $\ell(\mathbf{c})$ which we wish to
minimize\footnote{The loss function does \emph{not} have to be
restricted to the least-squares loss considered above.} and a
current estimate, $\mathbf{\hat c}$ of optimal solution, we can
determine an update of our estimate to a new value $\mathbf{\hat
c}_{\text{\tiny new}}$ which will decrease the loss function as
follows.

\medskip
For the loss function $\ell(\mathbf{c})$, with $\mathbf{c} =
\mathbf{\hat c} + d \mathbf{c}$, we have
\[ d \ell(\mathbf{\hat c}) = \ell(\mathbf{\hat c} + d \mathbf{c}) -
\ell(\mathbf{\hat c}) = \frac{\partial \ell(\mathbf{\hat
c})}{\partial \mathbf{c}} d \mathbf{c}  \] which is just the
differential limit of the first order expansion
\[ \Delta \ell(\mathbf{\hat c}; \alpha) = \ell(\mathbf{\hat c} + \alpha \Delta
\mathbf{c}) - \ell(\mathbf{\hat c}) \approx \alpha \frac{\partial
\ell(\mathbf{\hat c})}{\partial \mathbf{c}} \Delta \mathbf{c} \, .
\] The stepsize $\alpha > 0$ is a control parameter which regulates
the accuracy of the first order approximation assuming that \[
\alpha \rightarrow 0 \Rightarrow \alpha \Delta \mathbf{c}
\rightarrow d \mathbf{c} \quad \text{and} \quad \Delta
\ell(\mathbf{\hat c};\alpha) \rightarrow d \ell(\mathbf{\hat c})\, .
\]

\medskip
If we assume that $\mathcal{C}$ is a Cartesian space,\footnote{I.e.,
We assume that $\mathcal{C}$ has identity metric tensor. We call the resulting gradient a Cartesian gradient (if
the metric tensor assumption $\Omega_{\mathbf{c}}  = I$ is true for the space of interest) or
a naive gradient (if the identity metric tensor assumption is false, but made
anyway for convenience).} then the gradient of $\ell(\mathbf{c})$ is
given by\footnote{Note for future reference that the gradient has
been specifically computed in Equation (\ref{eq:quadgrad}) for the
special case when $\ell(\mathbf{c})$ is the least-squares loss
function (\ref{eq:weightedlsloss}).}
\[ \nabla_{\mathbf{c}} \ell(\mathbf{c}) = \left( \frac{\partial \ell(\mathbf{c})}{\partial
\mathbf{c}} \right)^H \, . \]

Take the update to be the \emph{generalized gradient descent correction}
 \begin{equation}\label{eq:ggdupdate}
 \Delta \mathbf{c} =  - Q
(\mathbf{\hat c} )\left( \frac{\partial \ell(\mathbf{\hat
c})}{\partial \mathbf{c}} \right)^H = - Q(\mathbf{\hat c}) \,
\nabla_{\mathbf{c}} \ell(\mathbf{\hat c}) \end{equation} where
$Q(\mathbf{\hat c})$ is a Hermitian matrix function of $\mathbf{c}$
which is assumed to be positive definite when evaluated at the value
$\mathbf{\hat c}$.\footnote{The fact that $Q$ is otherwise \emph{\emph{arbitrary}}
(except for the admissibility criterion discussed below) is what
makes the resulting algorithm a \emph{generalized} gradient descent
algorithm. When $Q = I$, we obtain
the standard (naive) gradient descent algorithm.} This then yields the key
stability condition\footnote{We interpret the stability condition to
mean that for a small enough stepsize $\alpha >0$, we will have
$\Delta \ell(\mathbf{\hat c}; \alpha) \le 0$.}
\begin{equation}\label{eq:stabcon} \Delta \ell(\mathbf{\hat c};
\alpha) \approx - \alpha \| \nabla_{\mathbf{c}} \ell(\mathbf{\hat
c}) \|^2_Q \triangleq - \alpha \, \nabla_{\mathbf{c}}
\ell(\mathbf{\hat c})^H \, Q \, \nabla_{\mathbf{c}}
\ell(\mathbf{\hat c}) \le 0 ,\end{equation} where the
right-hand-side is equal to zero if and only if
\[ \nabla_{\mathbf{c}} \ell(\mathbf{\hat c}) = 0 \, . \]

\medskip
Thus if the stepsize parameter $\alpha$ is chosen small enough,
making the update
\[ \mathbf{\hat c}_{\text{\tiny new}} = \mathbf{\hat c} + \alpha
\Delta \mathbf{c} = \mathbf{\hat c} - Q \, \nabla_{\mathbf{c}}
\ell(\mathbf{\hat c})  \] results in
\[ \ell(\mathbf{\hat c}_{\text{\tiny new}}) = \ell(\mathbf{\hat c} + \alpha
\Delta \mathbf{c}) = \ell(\mathbf{\hat c}) + \Delta
\ell(\mathbf{\hat c}; \alpha) \approx \ell(\mathbf{\hat c}) - \alpha
\| \nabla_{\mathbf{c}} \ell(\mathbf{\hat c}) \|^2_Q \le
\ell(\mathbf{\hat c}) \] showing that we \emph{either} have a
nontrivial update of the value  of $\mathbf{\hat c}$ which results
in a strict decrease in the value of the loss function, \emph{or} we
have no update of $\mathbf{\hat c}$ nor decrease of the loss
function because $\mathbf{\hat c}$ is a stationary point. If the
loss function $\ell(\mathbf{c})$ is bounded from below, iterating on
this procedure starting from a estimate $\mathbf{\hat c}_1$ will
produce a sequence of estimates $\mathbf{\hat c}_i$, $i = 1, 2, 3,
\cdots$, which will converge to a local minimum of the loss
function. This simple procedure is the basis for all generalized
gradient descent algorithms.

\medskip
Assuming that we begin with an admissible estimate, $\mathbf{\hat
c}_1$, for this procedure to be valid, we require that the sequence
of estimates $\mathbf{\hat c}_i$, $i = 1, 2, 3, \cdots$, be
admissible, which is true if the corresponding updates $\Delta
\mathbf{c}$ are admissible,
\[ \Delta \mathbf{c} = - Q(\mathbf{\hat c}_i) \nabla_{\mathbf{\hat c}_i} \ell(\mathbf{\hat c}_i)
= - Q(\mathbf{\hat c}_i) \left( \frac{\partial \ell(\mathbf{\hat
c}_i)}{\partial \mathbf{\hat c}_i} \right)^H  \in \mathcal{C} \, ,
\quad i = 1, 2, \cdots \, .
\] We have established the admissibility of $\nabla_{\mathbf{c}}
\ell(\mathbf{c}) = \left( \frac{\partial \ell(\mathbf{c})}{\partial
\mathbf{c}} \right)^H \in \mathcal{C}$ above. \emph{It is evident
that in order for a generalized gradient descent algorithm (GDA) to
be admissible it must be the case that $Q$ be admissible,}
\[ \boxed{  \ \text{\small $\text{Generalized GDA is Admissible} \Leftrightarrow \text{Generalized Gradient $Q$-Matrix is Admissible,}
\ Q \in \mathcal{L}(\mathcal{C},\mathcal{C})$}  \, . \ } \] Furthermore,
\emph{a sufficient condition that the resulting algorithm be
stablizable\footnote{I.e., that a small enough step size can be chosen to ensure that
the stability condition (\ref{eq:stabcon}) is satisfied.} is that
$Q$ be Hermitian and positive definite.} Note that given a candidate
Hermitian positive definite matrix, $Q'$, which is not admissible,
\[ Q' \notin \mathcal{L}(\mathcal{C},\mathcal{C}) \, , \] we can
transform it into an admissible Hermitian positive definite matrix
via the projection
\[ Q = \mathbf{P}(Q') \in \mathcal{L}(\mathcal{C},\mathcal{C}) \, . \]
It can be much trickier to ensure that $Q$ remains positive definite under the action of $\mathbf{P}$.

\medskip
If we set
\[ Q^{\text{\tiny Newton}}(\mathbf{c}) = \left[  \mathcal{H}_{\mathbf{c}\mathbf{c}}^{\text{\tiny Newton}}(\mathbf{c})\right]^{-1} \] with \[
\mathcal{H}_{\mathbf{c}\mathbf{c}}^{\text{\tiny Newton}} = 
\mathcal{H}_{\mathbf{c}\mathbf{c}}^{\text{\boldmath \tiny
$\mathbb{C}$}} \] then we obtain the Newton algorithm
(\ref{eq:newtoncorrection}). If we take the loss function to be the
least-squares loss function (\ref{eq:weightedlsloss}) and set
\[ Q^{\text{\tiny Gauss}}(\mathbf{c}) = \left[\mathcal{H}_{\mathbf{c}\mathbf{c}}^{\text{\tiny Gauss}}(\mathbf{
c})\right]^{-1} \] we obtain the Gauss-Newton algorithm
(\ref{eq:gausscorrection}). Whereas the Gauss-Newton algorithm
generally has a positive definite $Q$-matrix (assuming that
$g(\mathbf{c})$ is one-to-one), the Newton algorithm can have
convergence problems due to the Newton Hessian
$\mathcal{H}_{\mathbf{c}\mathbf{c}}^{\text{\tiny
Newton}}=\mathcal{H}_{\mathbf{c}\mathbf{c}}^{\text{\boldmath \tiny
$\mathbb{C}$}}$ becoming indefinite. Note that taking
\[ Q = I \, ,\] which we refer to as the ``Cartesian,'' ``standard,'' 
``simple,'' or ``naive'' choice (depending on the context) results in the standard gradient descent
algorithm which
 is stable for a small enough stepsize so that the stability condition (\ref{eq:stabcon}) holds.

\medskip
An important practical issue is the problem of stability
versus speed of convergence. It is well-known that the Newton
algorithm tends to have a very fast rate of convergence, but at the
cost of constructing and inverting the Newton Hessian
$\mathcal{H}_{\mathbf{c}\mathbf{c}}^{\text{\tiny
Newton}}=\mathcal{H}_{\mathbf{c}\mathbf{c}}^{\text{\boldmath \tiny
$\mathbb{C}$}}$ and potentially encountering more difficult
algorithm instability problems. On the other hand, standard gradient
descent ($Q=I$) tends to be very stable and much cheaper to
implement, but can have very long convergence times.

\medskip
The Gauss-Newton algorithm, which is an option available when the
loss function $\ell(\mathbf{c})$ is the least-squares loss function
(\ref{eq:weightedlsloss}), is considered an excellent trade-off
between the Newton algorithm and standard gradient descent. The
Gauss-Newton Hessian
$\mathcal{H}_{\mathbf{c}\mathbf{c}}^{\text{\tiny Gauss}}$ is
generally simpler in form and, if $\mathbf{g}(\mathbf{c})$ is
one-to-one, is always positive definite. Furthermore, if
$\mathbf{g}(\mathbf{c})$ is also onto, assuming the algorithm
converges, the Gauss-Newton and Newton algorithms are asymptotically
equivalent.

\medskip
We can also begin to gain some insight into the proposal by Yan and
Fan \cite{GY00} to ignore the block off-diagonal elements of the
Newton Hessian,\footnote{The values of the block elements of
$\mathcal{H}_{\mathbf{c}\mathbf{c}}^{\text{\tiny Newton}}$ will be
computed for the special case of the least-squares loss function
(\ref{eq:weightedlsloss}) later below.}
\[ \mathcal{H}_{\mathbf{c}\mathbf{c}}^{\text{\tiny Newton}} = \mathcal{H}_{\mathbf{c} \mathbf{c}}^{\text{\boldmath \tiny
$\mathbb{C}$}} = \begin{pmatrix} \mathcal{H}_{\mathbf{z}\mathbf{z}}
& \mathcal{H}_{\mathbf{\bar z}\mathbf{z}} \\
\mathcal{H}_{\mathbf{z}\mathbf{\bar z}} & \mathcal{H}_{\mathbf{\bar
z}\mathbf{\bar z}} \end{pmatrix} . \] As mentioned earlier, Yan and
Fan make the claim in  \cite{GY00}  that the block off-diagonal
elements vanish for a quadratic loss function. As noted above, and
shown in an example below, this is \emph{not} generally
true.\footnote{What is true, as we've noted, is that for a quadratic
loss function, the Gauss-Newton and Newton Hessians asymptotically
become equal.} However, it \emph{is}  reasonable to ask what harm
(if any), or what benefit (if any) can accrue by constructing a
\emph{new}\footnote{I.e., no approximation algorithms are invoked.}
generalized gradient descent algorithm as a modification to the
Newton algorithm created by simply ignoring the block off-diagonal
elements in the Newton Hessian and working instead with the
simplified \emph{quasi-Newton Hessian,}
\[ \mathcal{H}_{\mathbf{c}\mathbf{c}}^{\text{\tiny quasi-Newton}} \triangleq  \widehat{\mathcal{H}}_{\mathbf{c} \mathbf{c}}^{\text{\boldmath \tiny
$\mathbb{C}$}} \triangleq \begin{pmatrix}
\mathcal{H}_{\mathbf{z}\mathbf{z}}
& 0 \\
0 & \mathcal{H}_{\mathbf{\bar z}\mathbf{\bar z}} \end{pmatrix} . \]
This results in  a new generalized gradient descent algorithm, which
we call the \emph{quasi-Newton algorithm}, which is somewhere in
complexity between the Newton algorithm and standard gradient
descent. Note that the hermitian matrix
$\mathcal{H}_{\mathbf{z}\mathbf{z}}$ is positive definite if and
only if $\mathcal{H}_{\mathbf{\bar z}\mathbf{\bar z}}$ is positive
definite. Thus invertibility and positive-definiteness of the
quasi-Newton Hessian
$\mathcal{H}_{\mathbf{c}\mathbf{c}}^{\text{\tiny
quasi-Newton}}=\widehat{\mathcal{H}}_{\mathbf{c}
\mathbf{c}}^{\text{\boldmath \tiny $\mathbb{C}$}}$ is equivalent to
invertibility and positive definiteness of the block element
$\mathcal{H}_{\mathbf{z}\mathbf{z}}$.

\medskip
On the other hand, invertibility and positive definiteness of
$\mathcal{H}_{\mathbf{z}\mathbf{z}}$ is only a necessary condition
for invertibility and positive definiteness of the complete Newton
Hessian $\mathcal{H}_{\mathbf{c}\mathbf{c}}^{\text{\tiny Newton}}=
\mathcal{H}_{\mathbf{c} \mathbf{c}}^{\text{\boldmath \tiny
$\mathbb{C}$}}$. Assuming that $\mathcal{H}_{\mathbf{c}
\mathbf{c}}^{\text{\boldmath \tiny $\mathbb{C}$}}$  is positive
definite, we have the well-known factorization
\begin{equation}\label{eq:wkfactorization}
\begin{pmatrix} I
& 0 \\
- \mathcal{H}_{\mathbf{z}\mathbf{\bar z}}
\mathcal{H}_{\mathbf{z}\mathbf{z}}^{-1} & I
\end{pmatrix}\mathcal{H}_{\mathbf{c} \mathbf{c}}^{\text{\boldmath
\tiny $\mathbb{C}$}} \begin{pmatrix} I
& - \mathcal{H}_{\mathbf{\bar z}\mathbf{z}}\mathcal{H}_{\mathbf{z}\mathbf{z}}^{-1} \\
0 & I \end{pmatrix} =
\begin{pmatrix} \mathcal{H}_{\mathbf{z}\mathbf{z}}
& 0\\
0 & \widetilde{\mathcal{H}}_{\mathbf{z}\mathbf{z}} \end{pmatrix}
\end{equation}
 where \[ \widetilde{\mathcal{H}}_{\mathbf{z}\mathbf{z}} =
\mathcal{H}_{\mathbf{z}\mathbf{z}} - \mathcal{H}_{\mathbf{\bar
z}\mathbf{z}} \mathcal{H}_{\mathbf{\bar z}\mathbf{\bar
z}}^{-1}\mathcal{H}_{\mathbf{z}\mathbf{\bar z}} \] is the Schur
complement of $\mathcal{H}_{\mathbf{z}\mathbf{z}}$ in
$\mathcal{H}_{\mathbf{c} \mathbf{c}}^{\text{\boldmath \tiny
$\mathbb{C}$}}$. From the factorization (\ref{eq:wkfactorization})
we immediately obtain the useful condition
\begin{equation}\label{eq:ranksum}
\text{rank}\left(\mathcal{H}_{\mathbf{c}
\mathbf{c}}^{\text{\boldmath \tiny $\mathbb{C}$}}\right) =
\text{rank}\left(\mathcal{H}_{\mathbf{z}\mathbf{z}}\right) +
\text{rank}\left(\widetilde{\mathcal{H}}_{\mathbf{z}\mathbf{z}}\right)
\, .
\end{equation}

\medskip
Note from condition (\ref{eq:ranksum}) that the Newton Hessian
$\mathcal{H}_{\mathbf{c}\mathbf{c}}^{\text{\tiny
Newton}}=\mathcal{H}_{\mathbf{c} \mathbf{c}}^{\text{\boldmath \tiny
$\mathbb{C}$}}$ is positive definite if and only if
$\mathcal{H}_{\mathbf{z}\mathbf{z}}$ \emph{and} its Schur complement
$\widetilde{\mathcal{H}_{\mathbf{z}\mathbf{z}}}$ are \emph{both}
positive definite. Thus it is obviously a more difficult matter to
ascertain and ensure the stability of the Newton Hessian than to do
the same for the quasi-Newton Hessian.

\medskip

The quasi-Newton algorithm is constructed  by forming the $Q$ matrix
from the quasi-Newton Hessian
$\mathcal{H}_{\mathbf{c}\mathbf{c}}^{\text{\tiny
quasi-Newton}}=\widehat{\mathcal{H}}_{\mathbf{c}
\mathbf{c}}^{\text{\boldmath \tiny $\mathbb{C}$}}$,
\[ Q^{\text{\tiny Pseudo-Newton}}= \left({\mathcal{H}}_{\mathbf{c} \mathbf{c}}^{\text{\tiny
quasi-Newton}}\right)^{-1} = \left(\widehat{\mathcal{H}}_{\mathbf{c}
\mathbf{c}}^{\text{\boldmath \tiny $\mathbb{C}$}}\right)^{-1}
=\begin{pmatrix} \mathcal{H}_{\mathbf{z}\mathbf{z}}^{-1}
& 0 \\
0 & \mathcal{H}_{\mathbf{\bar z}\mathbf{\bar z}}^{-1}
\end{pmatrix}
\] which is admissible and hermitian, and positive definite provided
$\mathcal{H}_{\mathbf{z}\mathbf{z}} =
\overline{\mathcal{H}_{\mathbf{\bar z}\mathbf{\bar z}}}$ is positive
definite. Thus, if $\mathcal{H}_{\mathbf{z}\mathbf{z}} =
\overline{\mathcal{H}_{\mathbf{\bar z}\mathbf{\bar z}}}$ is positive
definite, \emph{the quasi-Newton algorithm is guaranteed to be
stable} (assuming a small enough stepsize $\alpha >0$ so that the
stability condition (\ref{eq:stabcon}) is satisfied). With this
choice of $Q$ in (\ref{eq:ggdupdate}), the quasi-Newton update is
given by\footnote{We can ignore the remaining update equation as it
is just the complex conjugate of the shown update equation.}
\begin{equation}\label{eq:yancorrect} {\Delta \mathbf{ z}}^{\text{\tiny quasi-Newton}} = -
\mathcal{H}_{\mathbf{z}\mathbf{z}}^{-1} \text{\footnotesize
$\left(\frac{\partial f}{\partial \mathbf{z}}\right)^H$}
\end{equation} which is just the simplification shown earlier in
Equation (\ref{eq:fansimplification}) and proposed by Yan and Fan in
\cite{GY00}. However, unlike Yan and Fan, we do not present the
quasi-Newton algorithm as an approximation to the Newton algorithm,
\emph{but rather as \underline{one more algorithm} in the \underline{family} of generalized
Newton algorithms indexed by the choice of the matrix $Q$.}

\medskip
Indeed, recognizing that the Gauss-Newton algorithm potentially has
better stability properties than the Newton algorithm, naturally
leads us to propose a \emph{quasi-Gauss-Newton} algorithm for
minimizing the least-squares lose function (\ref{eq:weightedlsloss})
as follows. Because the hermitian Gauss-Newton Hessian is
admissible, it can be partitioned as
\[ \mathcal{H}_{\mathbf{c}\mathbf{c}}^{\text{\tiny Gauss}} =
\begin{pmatrix} U_{\mathbf{z}\mathbf{z}} & U_{\mathbf{\bar z}\mathbf{z}} \\ \overline{U_{\mathbf{\bar z}\mathbf{z}}} & \overline{U_{\mathbf{z}\mathbf{z}}} \end{pmatrix}
\]
with $U_{\mathbf{\bar z}\mathbf{z}} = U_{\mathbf{\bar
z}\mathbf{z}}^T$.\footnote{The values of these block components will
be computed below.} The Gauss-Newton Hessian is positive-definite if
and only if $U_{\mathbf{z}\mathbf{z}}$ (equivalently
$\overline{U_{\mathbf{z}\mathbf{z}}}$) and its Schur complement
$\widetilde{U_{\mathbf{z}\mathbf{z}}} = U_{\mathbf{z}\mathbf{z}} -
U_{\mathbf{\bar z}\mathbf{z}}
\overline{U_{\mathbf{z}\mathbf{z}}}^{-1} \overline{U_{\mathbf{\bar
z}\mathbf{z}}}$ are invertible.

\medskip
On the other hand the \emph{quasi-Gauss-Newton Hessian,}
\[ \mathcal{H}_{\mathbf{c}\mathbf{c}}^{\text{\tiny quasi-Gauss}}
\triangleq
\begin{pmatrix} U_{\mathbf{z}\mathbf{z}} & 0 \\ 0 & \overline{U_{\mathbf{z}\mathbf{z}}} \end{pmatrix}
\]
is positive definite if and only if $U_{\mathbf{z}\mathbf{z}}$ is
positive definite. Choosing
\[ Q^{\text{\tiny quasi-Gauss}} = \left(\mathcal{H}_{\mathbf{c}\mathbf{c}}^{\text{\tiny
quasi-Gauss}}\right)^{-1}= \begin{pmatrix}
U_{\mathbf{z}\mathbf{z}}^{-1} & 0 \\ 0 &
\overline{U_{\mathbf{z}\mathbf{z}}}^{-1}
\end{pmatrix} \]
results in the \emph{quasi-Gauss-Newton algorithm}
\begin{equation}\label{eq:pgn} {\Delta \mathbf{ z}}^{\text{\tiny quasi-Gauss}} = -
U_{\mathbf{z}\mathbf{z}}^{-1} \text{\footnotesize
$\left(\frac{\partial f}{\partial \mathbf{z}}\right)^H$}
\end{equation}which is
guaranteed to be stable (for a small enough stepsize so that the
stability condition (\ref{eq:stabcon}) is satisfied) if
$U_{\mathbf{z}\mathbf{z}}$ is positive definite.

\medskip
Note that $\mathcal{H}_{\mathbf{z}\mathbf{z}}$ can become indefinite
even while $U_{\mathbf{z}\mathbf{z}}$ remains positive definite.
Thus, the quasi-Gauss-Newton algorithm appears to be generally
easier to stabilize than the quasi-Newton algorithm. Furthermore, if
$\mathbf{g}$ is onto, we expect that asymptotically the
quasi-Gauss-Newton and quasi-Newton algorithm become equivalent.
Thus the quasi-Gauss-Newton algorithm is seen to stand in the same
relationship to the quasi-Newton algorithm as the Gauss-Newton
algorithm does to the Newton algorithm.

\medskip
Without too much effort, we can construct the block matrix
components needed to implement the Newton and Gauss-Newton
algorithms developed above in order to minimize the least-squares
loss function (\ref{eq:weightedlsloss}).\footnote{This, of course,
results in only a special case application of the Newton and
quasi-Newton algorithms, both of which can be applied to more
general loss functions.}

\medskip Let us first look at the elements needed to implement the Gauss-Newton algorithm.
From Equation (\ref{eq:gaussnewtonhessian}) and the derivations
following Equation (\ref{eq:gaussloss1}) one obtains
\begin{equation}\label{eq:huzz} U_{\mathbf{z}\mathbf{z}} =\frac{1}{2}\left(
\left(\frac{\partial \mathbf{g}}{\partial \mathbf{z}}\right)^H \,  W
\left(\frac{\partial \mathbf{g}}{\partial \mathbf{z}}\right) +
\overline{\left(\frac{\partial \mathbf{g}}{\partial \mathbf{\bar
z}}\right)^H \, W \left(\frac{\partial \mathbf{g}}{\partial
\mathbf{\bar z}}\right)}\, \right)
\end{equation}
which is positive definite, assuming that $W$ is positive definite
and that $\mathbf{g}$ is one-to-one. Similarly, one finds that
\begin{equation}\label{eq:hubarzz} U_{\mathbf{\bar z}\mathbf{z}} = \frac{1}{2}\left(
\left(\frac{\partial \mathbf{g}}{\partial \mathbf{z}}\right)^H \,  W
\left(\frac{\partial \mathbf{g}}{\partial \mathbf{\bar z}}\right) +
\overline{\left(\frac{\partial \mathbf{g}}{\partial \mathbf{\bar
z}}\right)^H \, W \left(\frac{\partial \mathbf{g}}{\partial
\mathbf{z}}\right)}\, \right) \, .
\end{equation}
Also $U_{\mathbf{\bar z}\mathbf{\bar z}} =
\overline{U_{\mathbf{z}\mathbf{z}}}$ and $U_{\mathbf{ z}\mathbf{\bar
z}} = \overline{U_{\mathbf{\bar z}\mathbf{z}}}$. We have now
completely specified the Gauss-Newton Hessian
$\mathcal{H}_{\mathbf{c}\mathbf{c}}^{\text{\tiny Gauss}}$ and the
quasi-Gauss-Newton Hessian at the block components level,
\[ \mathcal{H}_{\mathbf{c}\mathbf{c}}^{\text{\tiny Gauss}} =
\begin{pmatrix} U_{\mathbf{z}\mathbf{z}} & U_{\mathbf{\bar z}\mathbf{z}} \\ {U_{\mathbf{z}\mathbf{\bar z}}}
& {U_{\mathbf{\bar z}\mathbf{\bar z}}} \end{pmatrix} \qquad
\mathcal{H}_{\mathbf{c}\mathbf{c}}^{\text{\tiny quasi-Gauss}}
\triangleq
\begin{pmatrix} U_{\mathbf{z}\mathbf{z}} & 0 \\ 0 & {U_{\mathbf{\bar z}\mathbf{\bar
z}}} \end{pmatrix}
\]

Now \emph{note the important fact that $U_{\mathbf{\bar
z}\mathbf{z}} =  {U_{\mathbf{z}\mathbf{\bar z}}} = 0$ when
$\mathbf{g}$ is holomorphic}! Thus, when $\mathbf{g}$ is holomorphic
there is no difference between the Gauss-Newton and
pseudo-Gauss-Newton algorithms.\footnote{Recall that
$\mathbf{g}(\mathbf{z})$ is holomorphic (analytic in $\mathbf{z}$)
if and only if the Cauchy-Riemann condition $\frac{\partial
\mathbf{g}(\mathbf{z})}{\partial \mathbf{\bar z}} = 0$ is
satisfied.} Furthermore, when $\mathbf{g}(\mathbf{z})$ is
holomorphic, $U_{\mathbf{z}\mathbf{z}}$ simplifies to
\begin{equation}\label{eq:simphuzz} U_{\mathbf{z}\mathbf{z}} =\frac{1}{2}
\left(\frac{\partial \mathbf{g}}{\partial \mathbf{z}}\right)^H \,  W
\left(\frac{\partial \mathbf{g}}{\partial \mathbf{z}}\right) =
\frac{1}{2}J_{\mathbf{g}}^H W  J_{\mathbf{g}}\, ,
\end{equation}
where $J_{\mathbf{g}}$ is the Jacobian matrix of $\mathbf{g}$.

\medskip
Now let us turn to the issue of computing the elements need to
implement the Newton Algorithm, recalling that the Newton Hessian is
block partitioned as
\[\mathcal{H}_{\mathbf{c}\mathbf{c}}^{\text{\tiny Newton}} =
\mathcal{H}_{\mathbf{c} \mathbf{c}}^{\text{\boldmath \tiny
$\mathbb{C}$}} = \begin{pmatrix} \mathcal{H}_{\mathbf{z}\mathbf{z}}
& \mathcal{H}_{\mathbf{\bar z}\mathbf{z}} \\
\mathcal{H}_{\mathbf{z}\mathbf{\bar z}} & \mathcal{H}_{\mathbf{\bar
z}\mathbf{\bar z}} \end{pmatrix} \, . \]
 One can readily relate the block components
$\mathcal{H}_{\mathbf{z}\mathbf{z}}$ and $\mathcal{H}_{\mathbf{\bar
z}\mathbf{z}}$ to the matrices $U_{\mathbf{z}\mathbf{z}}$ and
$U_{\mathbf{\bar z}\mathbf{z}}$ used in the Gauss-Newton and
quasi-Gauss-Newton algorithms by use of Equation (\ref{eq:syminhh}).
We find that \[ \mathcal{H}_{\mathbf{z}\mathbf{z}} =
U_{\mathbf{z}\mathbf{z}} - \sum_{i=1}^m
{V}_{\mathbf{z}\mathbf{z}}^{\text{\tiny $(i)$}}
\]  with
\begin{equation}\label{eq:vizz} {V}_{\mathbf{z}\mathbf{z}}^{\text{\tiny $(i)$}} =
\frac{1}{2}\left[ \left(\frac{\partial}{\partial \mathbf{z}}\left(
\frac{\partial g_i(\mathbf{z})}{\partial \mathbf{z}}\right)^H \,
\left[ W \mathbf{e} \, \right]_i \right) +
\overline{\left(\frac{\partial}{\partial \mathbf{\bar z}}\left(
\frac{\partial g_i(\mathbf{z})}{\partial \mathbf{\bar z}}\right)^H
\, \left[ W \mathbf{e} \, \right]_i\right)} \right]
\end{equation}
where $\mathbf{e} = \mathbf{y} - \mathbf{g}(\mathbf{z})$. Similarly,
we find that
\[
\mathcal{H}_{\mathbf{\bar z}\mathbf{z}} = U_{\mathbf{\bar
z}\mathbf{z}} - \sum_{i=1}^m {V}_{\mathbf{\bar
z}\mathbf{z}}^{\text{\tiny $(i)$}}
\]  and
\begin{equation}\label{eq:vibarzz} {V}_{\mathbf{\bar z}\mathbf{z}}^{\text{\tiny $(i)$}} =
\frac{1}{2}\left[ \left(\frac{\partial}{\partial \mathbf{\bar
z}}\left( \frac{\partial g_i(\mathbf{z})}{\partial \mathbf{
z}}\right)^H\, \left[ W \mathbf{e} \, \right]_i\right) +
\overline{\left(\frac{\partial}{\partial \mathbf{ z}}\left(
\frac{\partial g_i(\mathbf{z})}{\partial \mathbf{ \bar
z}}\right)^H\, \left[ W \mathbf{e} \, \right]_i\right)} \right] \,
\end{equation}
Furthermore, ${V}_{\mathbf{\bar z}\mathbf{\bar z}} =
\overline{{V}_{\mathbf{z}\mathbf{z}}}$ and ${V}_{\mathbf{
z}\mathbf{\bar z}} = \overline{{V}_{\mathbf{\bar z}\mathbf{z}}}$.

\medskip
Note that neither ${V}_{\mathbf{z}\mathbf{z}}$ nor
${V}_{\mathbf{\bar z}\mathbf{z}}$  vanish when $\mathbf{g}$ is
holomorphic, but instead simplify to
\begin{equation}\label{eq:simpvizz}
{V}_{\mathbf{z}\mathbf{z}}^{\text{\tiny $(i)$}} = \frac{1}{2}\
\frac{\partial}{\partial \mathbf{z}}\left( \frac{\partial
g_i(\mathbf{z})}{\partial \mathbf{z}}\right)^H \, \left[ W
\mathbf{e} \, \right]_i  \quad \text{and} \quad {V}_{\mathbf{\bar
z}\mathbf{z}}^{\text{\tiny $(i)$}} = \frac{1}{2}\
\frac{\partial}{\partial \mathbf{\bar z}}\left( \frac{\partial
g_i(\mathbf{z})}{\partial \mathbf{z}}\right)^H \, \left[ W
\mathbf{e} \, \right]_i  \, . \end{equation}

\medskip
We have shown that the relationship between the Newton Hessian and
Gauss-Newton Hessian is given by
\[ \underset{\mathcal{H}_{\mathbf{c}\mathbf{c}}^{\text{\tiny Newton}}}{\underbrace{\begin{pmatrix} \mathcal{H}_{\mathbf{z}\mathbf{z}}
& \mathcal{H}_{\mathbf{\bar z}\mathbf{z}} \\
\mathcal{H}_{\mathbf{z}\mathbf{\bar z}} & \mathcal{H}_{\mathbf{\bar
z}\mathbf{\bar z}} \end{pmatrix}}} =
\underset{\mathcal{H}_{\mathbf{c}\mathbf{c}}^{\text{\tiny
Gauss}}}{\underbrace{\begin{pmatrix}
U_{\mathbf{z}\mathbf{z}} & U_{\mathbf{\bar z}\mathbf{z}} \\
{U_{\mathbf{z}\mathbf{\bar z}}} & {U_{\mathbf{\bar z}\mathbf{\bar
z}}} \end{pmatrix}}} - \sum_{i=1}^m\begin{pmatrix}
V_{\mathbf{z}\mathbf{z}}^{\text{\tiny $(i)$}} & V_{\mathbf{\bar z}\mathbf{z}}^{\text{\tiny $(i)$}} \\
{V_{\mathbf{z}\mathbf{\bar z}}^{\text{\tiny $(i)$}}} &
{V_{\mathbf{\bar z}\mathbf{\bar z}}^{\text{\tiny $(i)$}}}
\end{pmatrix}\]
In the special case when $\mathbf{g}(\mathbf{z})$ is holomorphic,
the relationship becomes
\[ \underset{\mathcal{H}_{\mathbf{c}\mathbf{c}}^{\text{\tiny Newton}}}{\underbrace{\begin{pmatrix} \mathcal{H}_{\mathbf{z}\mathbf{z}}
& \mathcal{H}_{\mathbf{\bar z}\mathbf{z}} \\
\mathcal{H}_{\mathbf{z}\mathbf{\bar z}} & \mathcal{H}_{\mathbf{\bar
z}\mathbf{\bar z}} \end{pmatrix}}} =
\underset{\mathcal{H}_{\mathbf{c}\mathbf{c}}^{\text{\tiny
Gauss}}}{\underbrace{\begin{pmatrix}
U_{\mathbf{z}\mathbf{z}} & 0 \\
0 & {U_{\mathbf{\bar z}\mathbf{\bar z}}} \end{pmatrix}}} -
\text{\scriptsize $\frac{1}{2} \sum_{i=1}^m\begin{pmatrix}
\frac{\partial}{\partial \mathbf{z}}\left( \frac{\partial
g_i(\mathbf{z})}{\partial \mathbf{z}}\right)^H \, \left[ W
\mathbf{e} \, \right]_i  &  \frac{\partial}{\partial \mathbf{\bar
z}}\left( \frac{\partial g_i(\mathbf{z})}{\partial
\mathbf{z}}\right)^H \, \left[ W \mathbf{e} \, \right]_i \\
\overline{\frac{\partial}{\partial \mathbf{\bar z}}\left(
\frac{\partial g_i(\mathbf{z})}{\partial \mathbf{z}}\right)^H \,
\left[ W \mathbf{e} \, \right]_i} &
\overline{\frac{\partial}{\partial \mathbf{z}}\left( \frac{\partial
g_i(\mathbf{z})}{\partial \mathbf{z}}\right)^H \, \left[ W
\mathbf{e} \, \right]_i}
\end{pmatrix}$}  \, .\]
This shows that if $\mathbf{g}(\mathbf{z})$ is \emph{holomorphic,}
so that the block off-diagonal elements of the Gauss-Newton Hessian
vanish, and $\mathbf{g}(\mathbf{z})$ is also \emph{onto,} so that
asymptotically we expect that $\mathbf{e} \approx 0$, then the claim
of Yan and Fan in \cite{GY00}  that setting the block off-diagonal
elements of the Hessian matrix can proved a good
\emph{approximation} to the Hessian matrix is reasonable, at least
when optimizing the least-squares loss function. However, when
$\mathbf{e} \approx 0$ the Newton least-squares loss function
(\ref{eq:newtonloss}) reduces to the Gauss-Newton loss function
(\ref{eq:gaussloss2}), so that in the least-squares case one may as
well make the move immediately to the even simpler Gauss-Newton
algorithm (which in this case coincides with the quasi-Gauss-Newton
algorithm).

\medskip
However, the real point to be made is that \emph{any} generalized
gradient descent algorithm is  worthy of
consideration,\footnote{I.e., we don't have to necessarily invoke an
approximation argument.} provided that it is admissible, provably
stable, and (at least locally) convergent to the desired optimal
solution. After all the standard gradient descent algorithm
corresponds to the cheapest ``approximation'' of all, namely that
\[ \mathcal{H}_{\mathbf{c}\mathbf{c}}^{\text{\tiny Newton}} \approx   I \,
\]
and very few will deny the utility of this algorithm, even though as
an ``approximation'' to the Newton algorithm it might be far from
correct. The resulting algorithm has intrinsic merit as an algorithm
in its own right, namely as the member of the family of gradient
descent algorithms corresponding to the simplest choice of the
$Q$-matrix,
\[ Q = I \, . \] In the end, if the
algorithm \emph{works}, it's ok. As it is said, ``the proof is in
the pudding.''\footnote{Of course, we are allowed to ask what the
performance of the $Q = I$ standard gradient-descent algorithm is relative
to the $Q^{\text{\tiny Newton}}$ algorithm.}

\medskip
We see, then, that we have a variety of algorithms at hand which fit
within the framework of generalized gradient descent algorithms.
These algorithms are characterized by the specific choice of the
$Q$-matrix in the gradient descent algorithm, and include (roughly
in the expected order of decreasing complexity, decreasing ideal
performance, and increasing stability when applied to the
least-squares loss function): 1) the Newton algorithm, 2) the
quasi-Newton algorithm, 3) the Gauss-Newton algorithm,  4) the
quasi-Gauss-Newton algorithm, and 5) standard gradient descent. Note
that the Newton, quasi-Newton, and standard gradient descent
algorithms are algorithms for minimizing a \emph{general}  loss function, while the Gauss-Newton and
quasi-Gauss-Newton algorithms are methods for minimizing the
\emph{least-squares} loss function (\ref{eq:weightedlsloss}).

\medskip
For convenience, we will now summarize the generalized gradient
descent algorithms that we have developed in this note.  In all of
the algorithms, the update step is given by
\[ \mathbf{\hat c} \leftarrow \mathbf{\hat c} + \alpha \Delta
\mathbf{c} \] or, equivalently,
\[ \mathbf{\hat z} \leftarrow \mathbf{\hat z} + \alpha \Delta
\mathbf{z} \] for a specific choice of the stepsize $\alpha
> 0$. The stability claims made are based on the assumption that
$\alpha$ has been chosen small enough to ensure that the stability
condition (\ref{eq:stabcon}) is valid. Furthermore, we use the
shorthand notation
\[ G(\mathbf{c}) = \frac{\partial \mathbf{g}(\mathbf{c})}{\partial
\mathbf{c}} \] and \[ \mathbf{e}(\mathbf{c}) = \mathbf{y} -
\mathbf{g}(\mathbf{c}) \, . \]
Note that in the taxonomy given below only the Newton Hessian $\mathcal{H}_{\mathbf{c}\mathbf{c}}^{\text{\tiny
Newton}}$ is generally the true Hessian of the loss function.
\newpage
\begin{enumerate}
\item {\bf Standard (a.k.a.\ Simple, Cartesian, or Naive) Gradient Descent.}

Applies to any smooth loss function which is bounded from below.



$Q^{\text{\tiny standard}}(\mathbf{\hat c}) = I$

$\Delta \mathbf{c}^{\text{\tiny standard}} = - \nabla_{\mathbf{z}}
\ell(\mathbf{\hat c}) = - \left(\frac{\partial \ell(\mathbf{\hat
c})}{\partial \mathbf{c}}\right)^H$

$\Delta \mathbf{z}^{\text{\tiny standard}} = - \nabla_{\mathbf{ z}}
\ell(\mathbf{\hat z}) = - \left(\frac{\partial \ell(\mathbf{\hat
z})}{\partial \mathbf{z}}\right)^H$

Application to Least-Squares Loss Function
(\ref{eq:weightedlsloss}):
\begin{itemize}
\item[]
$\left( \frac{\partial \ell}{\partial \mathbf{c}} \right)^H = -
\frac{1}{2} {G^H W \mathbf{e}} - \frac{1}{2} S\overline{G^H W
\mathbf{e}}  = \frac{1}{2} \left( {B}(\mathbf{\hat c}) + S
\overline{B}(\mathbf{\hat c}) \right)$

where \ $B(\mathbf{\hat c}) = - G(\mathbf{\hat c})^H W
e(\mathbf{\hat c})$

$\Delta \mathbf{c}^{\text{\tiny standard}} = - \frac{1}{2} \left[
B(\mathbf{\hat c}) + S \overline{B(\mathbf{\hat c})} \right]$ \quad

$\text{$\left(\frac{\partial \ell}{\partial \mathbf{z}}\right)^H$} =
- \frac{1}{2} \left[ \left( \frac{\mathbf g(\mathbf{\hat
z})}{\partial \mathbf{z}} \right)^H W \mathbf{e}(\mathbf{\hat z}) +
\overline{\left( \frac{\mathbf g(\mathbf{\hat z})}{\partial
\mathbf{\bar z}} \right)^H W \mathbf{e}(\mathbf{\hat z})} \, \,
\right]$

$\Delta \mathbf{z}^{\text{\tiny standard}} =  \frac{1}{2} \left[
\left( \frac{\mathbf g(\mathbf{\hat z})}{\partial \mathbf{z}}
\right)^H W \mathbf{e}(\mathbf{\hat z}) + \overline{\left(
\frac{\mathbf g(\mathbf{\hat z})}{\partial \mathbf{\bar z}}
\right)^H W \mathbf{e}(\mathbf{\hat z})}\,  \, \right]$

$\mathbf{g}(\mathbf{z})$ holomorphic:

$\text{$\left(\frac{\partial \ell}{\partial \mathbf{z}}\right)^H$} =
- \frac{1}{2} \left( \frac{\mathbf g(\mathbf{\hat z})}{\partial
\mathbf{z}} \right)^H W \mathbf{e}(\mathbf{\hat z})$

$\Delta \mathbf{z}^{\text{\tiny standard}} = \frac{1}{2}  \left(
\frac{\mathbf g(\mathbf{\hat z})}{\partial \mathbf{z}} \right)^H W
\mathbf{e}(\mathbf{\hat z})$
\end{itemize}
Generally stable but slow.

\item {\bf Gauss-Newton Algorithm.}

Applies to the least-squares loss function
(\ref{eq:weightedlsloss}).

$\mathcal{H}_{\mathbf{c}\mathbf{c}}^{\text{\tiny
Gauss}}(\mathbf{\hat c}) = \begin{pmatrix}
U_{\mathbf{z}\mathbf{z}} & U_{\mathbf{\bar z}\mathbf{z}} \\
{U_{\mathbf{z}\mathbf{\bar z}}} & {U_{\mathbf{\bar z}\mathbf{\bar
z}}} \end{pmatrix}$

where $U_{\mathbf{z}\mathbf{z}}$ is given by (\ref{eq:huzz}), \,
$U_{\mathbf{\bar z}\mathbf{\bar z}} =
\overline{U_{\mathbf{z}\mathbf{z}}}$,  $U_{\mathbf{\bar
z}\mathbf{z}}$ is given by (\ref{eq:hubarzz}), and
$U_{\mathbf{z}\mathbf{\bar z}} = \overline{U_{\mathbf{\bar
z}\mathbf{z}}}$.

$Q^{\text{\tiny Gauss}}(\mathbf{\hat c}) =
\mathcal{H}_{\mathbf{c}\mathbf{c}}^{\text{\tiny Gauss}}(\mathbf{
\hat c})^{-1}$

$\Delta \mathbf{c}^{\text{\tiny Gauss}} =  - Q^{\text{\tiny
Gauss}}(\mathbf{\hat c}) \left(\frac{\partial \ell(\mathbf{\hat
c})}{\partial \mathbf{c}}\right)^H$ \ where

$\left( \frac{\partial \ell}{\partial \mathbf{c}} \right)^H = -
\frac{1}{2} {G^H W \mathbf{e}} - \frac{1}{2} S\overline{G^H W
\mathbf{e}}  = \frac{1}{2} \left( {B}(\mathbf{\hat c}) + S
\overline{B}(\mathbf{\hat c}) \right)$

with \ $B(\mathbf{\hat c}) = - G(\mathbf{\hat c})^H W e(\mathbf{\hat
c})$

${\Delta \mathbf{z}}^{\text{\tiny Gauss}} = \left(
U_{\mathbf{z}\mathbf{z}} - U_{\mathbf{\bar z}\mathbf{z}}
U_{\mathbf{\bar z}\mathbf{\bar z}}^{-1}U_{\mathbf{z}\mathbf{\bar
z}}\right)^{-1} \left\{ U_{\mathbf{\bar z}\mathbf{z}}
U_{\mathbf{\bar z}\mathbf{\bar z}}^{-1} \text{ $\left(\frac{\partial
\ell}{\partial \mathbf{\bar z}}\right)^H$} -
\text{$\left(\frac{\partial \ell}{\partial \mathbf{z}}\right)^H$}
\right\}$ \quad where

$\text{$\left(\frac{\partial \ell}{\partial \mathbf{z}}\right)^H$} =
- \frac{1}{2} \left[ \left( \frac{\mathbf g(\mathbf{\hat
z})}{\partial \mathbf{z}} \right)^H W \mathbf{e}(\mathbf{\hat z}) +
\overline{\left( \frac{\mathbf g(\mathbf{\hat z})}{\partial
\mathbf{\bar z}} \right)^H W \mathbf{e}(\mathbf{\hat z})} \, \,
\right]$; \quad $ \text{$\left(\frac{\partial \ell}{\partial
\mathbf{\bar z}}\right)^H$} = \overline{\text{$\left(\frac{\partial
\ell}{\partial \mathbf{z}}\right)^H$}}$

$\mathbf{g}(\mathbf{z})$ holomorphic:

$U_{\mathbf{z}\mathbf{z}}$ takes the simpler form
(\ref{eq:simphuzz}), \, $U_{\mathbf{\bar z}\mathbf{\bar z}} =
\overline{U_{\mathbf{z}\mathbf{z}}}$, and $U_{\mathbf{z}\mathbf{\bar
z}} = \overline{U_{\mathbf{\bar z}\mathbf{z}}} = 0$.

$\mathcal{H}_{\mathbf{c}\mathbf{c}}^{\text{\tiny
Gauss}}(\mathbf{\hat c}) = \begin{pmatrix}
U_{\mathbf{z}\mathbf{z}} & 0 \\
0 & {U_{\mathbf{\bar z}\mathbf{\bar z}}} \end{pmatrix}= \frac{1}{2}
\begin{pmatrix} \left(\frac{\partial \mathbf{g}}{\partial
\mathbf{z}}\right)^H \, W \left(\frac{\partial \mathbf{g}}{\partial
\mathbf{z}}\right)
 &  0 \\  0 & \overline{\left(\frac{\partial
\mathbf{g}}{\partial \mathbf{z}}\right)^H \,  W \left(\frac{\partial
\mathbf{g}}{\partial \mathbf{z}}\right)}
\end{pmatrix}$

$\text{$\left(\frac{\partial \ell}{\partial \mathbf{z}}\right)^H$} =
- \frac{1}{2} \left( \frac{\mathbf g(\mathbf{\hat z})}{\partial
\mathbf{z}} \right)^H W \mathbf{e}(\mathbf{\hat z})$

${\Delta \mathbf{z}}^{\text{\tiny Gauss}} =
U_{\mathbf{z}\mathbf{z}}^{-1} \text{$\left(\frac{\partial
\ell}{\partial \mathbf{z}}\right)^H$}=  \left[ \left(\frac{\partial
\mathbf{g(\mathbf{\hat z})}}{\partial \mathbf{z}}\right)^H \, W
\left(\frac{\partial \mathbf{g(\mathbf{\hat z})}}{\partial
\mathbf{z}}\right)\right]^{-1} \left( \frac{\mathbf g(\mathbf{\hat
z})}{\partial \mathbf{z}} \right)^H W \mathbf{e}(\mathbf{\hat z})$

Stability generally requires positive definiteness of both
$U_{\mathbf{z}\mathbf{z}}$ and its Schur complement:
$\widetilde{U}_{\mathbf{z}\mathbf{z}} = U_{\mathbf{z}\mathbf{z}} -
U_{\mathbf{\bar z}\mathbf{z}} U_{\mathbf{\bar z}\mathbf{\bar
z}}^{-1}U_{\mathbf{z}\mathbf{\bar z}}$.  The need to step for
positive-definiteness of the Schur complement can significantly
increase the complexity of an on-line adaptive filtering algorithm.

If $\mathbf{g}(\mathbf{z})$ is holomorphic, then stability only
requires positive definiteness of the matrix
$U_{\mathbf{z}\mathbf{z}} = \left(\frac{\partial
\mathbf{g(\mathbf{\hat z})}}{\partial \mathbf{z}}\right)^H \, W
\left(\frac{\partial \mathbf{g(\mathbf{\hat z})}}{\partial
\mathbf{z}}\right)$, which will be the case if
$\mathbf{g}(\mathbf{z})$ is one-to-one. Thus, the algorithm may be
easier to stabilize when $\mathbf{g}(\mathbf{z})$ is holomorphic.

Convergence tends to be fast.

\item {\bf Pseudo-Gauss-Newton Algorithm.}

Applies to the least-squares loss function
(\ref{eq:weightedlsloss}).

$\mathcal{H}_{\mathbf{c}\mathbf{c}}^{\text{\tiny
Gauss}}(\mathbf{\hat c}) = \begin{pmatrix}
U_{\mathbf{z}\mathbf{z}} & 0 \\
0 & {U_{\mathbf{\bar z}\mathbf{\bar z}}} \end{pmatrix}$

where $U_{\mathbf{z}\mathbf{z}}$ is given by (\ref{eq:huzz}) \, and
$U_{\mathbf{\bar z}\mathbf{\bar z}} =
\overline{U_{\mathbf{z}\mathbf{z}}}$.

$Q^{\text{\tiny pseudo-Gauss}}(\mathbf{\hat c}) =
\left[\mathcal{H}_{\mathbf{c}\mathbf{c}}^{\text{\tiny
pseudo-Gauss}}(\mathbf{ \hat c})\right]^{-1} = \begin{pmatrix}
U_{\mathbf{z}\mathbf{z}}^{-1} & 0 \\
0 & {U_{\mathbf{\bar z}\mathbf{\bar z}}}^{-1} \end{pmatrix}$

$\Delta \mathbf{c}^{\text{\tiny pseudo-Gauss}} =  - Q^{\text{\tiny
pseudo-Gauss}}(\mathbf{\hat c}) \left(\frac{\partial
\ell(\mathbf{\hat c})}{\partial \mathbf{c}}\right)^H$ \ where

$\left( \frac{\partial \ell}{\partial \mathbf{c}} \right)^H = -
\frac{1}{2} {G^H W \mathbf{e}} - \frac{1}{2} S\overline{G^H W
\mathbf{e}}  = \frac{1}{2} \left( {B}(\mathbf{\hat c}) + S
\overline{B}(\mathbf{\hat c}) \right)$

with \ $B(\mathbf{\hat c}) = - G(\mathbf{\hat c})^H W e(\mathbf{\hat
c})$

${\Delta \mathbf{z}}^{\text{\tiny pseudo-Gauss}} = - \left[
U_{\mathbf{z}\mathbf{z}}(\mathbf{\hat z})\right]^{-1}
\text{$\left(\frac{\partial \ell(\mathbf{\hat z})}{\partial
\mathbf{z}}\right)^H$}=  \left[ \left(\frac{\partial
\mathbf{g}}{\partial \mathbf{z}}\right)^H \,  W \left(\frac{\partial
\mathbf{g}}{\partial \mathbf{z}}\right) +
\overline{\left(\frac{\partial \mathbf{g}}{\partial \mathbf{\bar
z}}\right)^H \, W \left(\frac{\partial \mathbf{g}}{\partial
\mathbf{\bar z}}\right)} \right]^{-1} \text{$\left(\frac{\partial
\ell(\mathbf{\hat z})}{\partial \mathbf{z}}\right)^H$}$ where

$\text{$\left(\frac{\partial \ell}{\partial \mathbf{z}}\right)^H$} =
- \frac{1}{2} \left[ \left( \frac{\mathbf g(\mathbf{\hat
z})}{\partial \mathbf{z}} \right)^H W \mathbf{e}(\mathbf{\hat z}) +
\overline{\left( \frac{\mathbf g(\mathbf{\hat z})}{\partial
\mathbf{\bar z}} \right)^H W \mathbf{e}(\mathbf{\hat z})} \, \,
\right]$

$\mathbf{g}(\mathbf{z})$ holomorphic:

$U_{\mathbf{z}\mathbf{z}}$ takes the simpler form of
(\ref{eq:simphuzz}) , and $U_{\mathbf{\bar z}\mathbf{\bar z}} =
\overline{U_{\mathbf{z}\mathbf{z}}}$.

 $\mathcal{H}_{\mathbf{c}\mathbf{c}}^{\text{\tiny
pseudo-Gauss}}(\mathbf{\hat c}) = \begin{pmatrix}
U_{\mathbf{z}\mathbf{z}} & 0 \\
0 & {U_{\mathbf{\bar z}\mathbf{\bar z}}} \end{pmatrix}= \frac{1}{2}
\begin{pmatrix} \left(\frac{\partial \mathbf{g}}{\partial
\mathbf{z}}\right)^H \, W \left(\frac{\partial \mathbf{g}}{\partial
\mathbf{z}}\right)
 &  0 \\  0 & \overline{\left(\frac{\partial
\mathbf{g}}{\partial \mathbf{z}}\right)^H \,  W \left(\frac{\partial
\mathbf{g}}{\partial \mathbf{z}}\right)}
\end{pmatrix}$

$\text{$\left(\frac{\partial \ell}{\partial \mathbf{z}}\right)^H$} =
- \frac{1}{2}  \left( \frac{\mathbf g(\mathbf{\hat z})}{\partial
\mathbf{z}} \right)^H W \mathbf{e}(\mathbf{\hat z})$

${\Delta \mathbf{z}}^{\text{\tiny pseudo-Gauss}} = \left[
\left(\frac{\partial \mathbf{g(\mathbf{\hat z})}}{\partial
\mathbf{z}}\right)^H \, W \left(\frac{\partial
\mathbf{g(\mathbf{\hat z})}}{\partial \mathbf{z}}\right)\right]^{-1}
\left( \frac{\mathbf g(\mathbf{\hat z})}{\partial \mathbf{z}}
\right)^H W \mathbf{e}(\mathbf{\hat z})$

Stability requires positive definiteness of
$U_{\mathbf{z}\mathbf{z}}(\mathbf{\hat z}) = \left(\frac{\partial
\mathbf{g(\mathbf{\hat z})}}{\partial \mathbf{z}}\right)^H \, W
\left(\frac{\partial \mathbf{g(\mathbf{\hat z})}}{\partial
\mathbf{z}}\right)$ which will be the case if
$\mathbf{g}(\mathbf{z})$ is one-to-one.

Convergence is expected to be quick but generally slower than for
Gauss-Newton due to loss of efficiency due to neglecting the block
off-diagonal terms in the Gauss-Newton Hessian (off-set, however, by
reduced complexity and possible gains in stability), except for the
case when $\mathbf{g}(\mathbf{z})$ is holomorphic, in which case the
two algorithms coincide.

\item {\bf Newton-Algorithm.}

Applies to any smooth loss function which is bounded from below.

$\mathcal{H}_{\mathbf{c}\mathbf{c}}^{\text{\tiny
Newton}}(\mathbf{\hat c}) =
\begin{pmatrix} \mathcal{H}_{\mathbf{z}\mathbf{z}}(\mathbf{\hat c})
& \mathcal{H}_{\mathbf{\bar z}\mathbf{z}}(\mathbf{\hat c}) \\
\mathcal{H}_{\mathbf{z}\mathbf{\bar z}}(\mathbf{\hat c}) &
\mathcal{H}_{\mathbf{\bar z}\mathbf{\bar z}}(\mathbf{\hat c})
\end{pmatrix}$

$Q^{\text{\tiny Newton}}(\mathbf{\hat c}) =
\left[\mathcal{H}_{\mathbf{c}\mathbf{c}}^{\text{\tiny
Newton}}(\mathbf{ \hat c})\right]^{-1}$

$\Delta \mathbf{c}^{\text{\tiny Newton}} =  - Q^{\text{\tiny
Newton}}(\mathbf{\hat c}) \left(\frac{\partial \ell(\mathbf{\hat
c})}{\partial \mathbf{c}}\right)^H$

$\Delta \mathbf{z}^{\text{\tiny Newton}} = \left(
\mathcal{H}_{\mathbf{z}\mathbf{z}} - \mathcal{H}_{\mathbf{\bar
z}\mathbf{z}} \mathcal{H}_{\mathbf{\bar z}\mathbf{\bar
z}}^{-1}\mathcal{H}_{\mathbf{z}\mathbf{\bar z}}\right)^{-1} \left\{
\mathcal{H}_{\mathbf{\bar z}\mathbf{z}} \mathcal{H}_{\mathbf{\bar
z}\mathbf{\bar z}}^{-1} \text{\footnotesize $\left(\frac{\partial
\ell}{\partial \mathbf{\bar z}}\right)^H$} - \text{\footnotesize
$\left(\frac{\partial \ell}{\partial \mathbf{z}}\right)^H$}
\right\}$

Application to the Least-Squares Loss Function
(\ref{eq:weightedlsloss}):
\begin{itemize}
\item[] {\small $\mathcal{H}_{\mathbf{c}\mathbf{c}}^{\text{\tiny Newton}}
= \begin{pmatrix} \mathcal{H}_{\mathbf{z}\mathbf{z}}
& \mathcal{H}_{\mathbf{\bar z}\mathbf{z}} \\
\mathcal{H}_{\mathbf{z}\mathbf{\bar z}} & \mathcal{H}_{\mathbf{\bar
z}\mathbf{\bar z}} \end{pmatrix} =
\begin{pmatrix}
U_{\mathbf{z}\mathbf{z}} & U_{\mathbf{\bar z}\mathbf{z}} \\
{U_{\mathbf{z}\mathbf{\bar z}}} & {U_{\mathbf{\bar z}\mathbf{\bar
z}}} \end{pmatrix} - \sum_{i=1}^m\begin{pmatrix}
V_{\mathbf{z}\mathbf{z}}^{\text{\tiny $(i)$}} & V_{\mathbf{\bar z}\mathbf{z}}^{\text{\tiny $(i)$}} \\
{V_{\mathbf{z}\mathbf{\bar z}}^{\text{\tiny $(i)$}}} &
{V_{\mathbf{\bar z}\mathbf{\bar z}}^{\text{\tiny $(i)$}}}
\end{pmatrix} \\ \text{\hskip1.47in} = \mathcal{H}_{\mathbf{c}\mathbf{c}}^{\text{\tiny
Gauss}}(\mathbf{\hat c}) - \sum_{i=1}^m\begin{pmatrix}
V_{\mathbf{z}\mathbf{z}}^{\text{\tiny $(i)$}} & V_{\mathbf{\bar z}\mathbf{z}}^{\text{\tiny $(i)$}} \\
{V_{\mathbf{z}\mathbf{\bar z}}^{\text{\tiny $(i)$}}} &
{V_{\mathbf{\bar z}\mathbf{\bar z}}^{\text{\tiny $(i)$}}}
\end{pmatrix}$}

$U_{\mathbf{z}\mathbf{z}}$ is given by (\ref{eq:huzz}), \,
$U_{\mathbf{\bar z}\mathbf{\bar z}} =
\overline{U_{\mathbf{z}\mathbf{z}}}$, \, $U_{\mathbf{\bar
z}\mathbf{z}}$ is given by (\ref{eq:hubarzz}),
$U_{\mathbf{z}\mathbf{\bar z}} = \overline{U_{\mathbf{\bar
z}\mathbf{z}}}$

$V_{\mathbf{z}\mathbf{z}}^{\text{\tiny $(i)$}}$ is given by
(\ref{eq:vizz}), \, $V_{\mathbf{\bar z}\mathbf{\bar z}}^{\text{\tiny
$(i)$}} = \overline{V_{\mathbf{z}\mathbf{z}}^{\text{\tiny $(i)$}}}$,
\, $V_{\mathbf{\bar z}\mathbf{z}}^{\text{\tiny $(i)$}}$ is given by
(\ref{eq:vibarzz}), $V_{\mathbf{z}\mathbf{\bar z}}^{\text{\tiny
$(i)$}} = \overline{V_{\mathbf{\bar z}\mathbf{z}}^{\text{\tiny
$(i)$}}}$.

$\Delta \mathbf{c}^{\text{\tiny Newton}} =  - Q^{\text{\tiny
Newton}}(\mathbf{\hat c}) \left(\frac{\partial \ell(\mathbf{\hat
c})}{\partial \mathbf{c}}\right)^H$ \ where

$\left( \frac{\partial \ell}{\partial \mathbf{c}} \right)^H = -
\frac{1}{2} {G^H W \mathbf{e}} - \frac{1}{2} S\overline{G^H W
\mathbf{e}}  = \frac{1}{2} \left( {B}(\mathbf{\hat c}) + S
\overline{B}(\mathbf{\hat c}) \right)$

with \ $B(\mathbf{\hat c}) = - G(\mathbf{\hat c})^H W e(\mathbf{\hat
c})$

$\Delta \mathbf{z}^{\text{\tiny Newton}} = \left(
\mathcal{H}_{\mathbf{z}\mathbf{z}} - \mathcal{H}_{\mathbf{\bar
z}\mathbf{z}} \mathcal{H}_{\mathbf{\bar z}\mathbf{\bar
z}}^{-1}\mathcal{H}_{\mathbf{z}\mathbf{\bar z}}\right)^{-1} \left\{
\mathcal{H}_{\mathbf{\bar z}\mathbf{z}} \mathcal{H}_{\mathbf{\bar
z}\mathbf{\bar z}}^{-1} \text{\footnotesize $\left(\frac{\partial
\ell}{\partial \mathbf{\bar z}}\right)^H$} - \text{\footnotesize
$\left(\frac{\partial \ell}{\partial \mathbf{z}}\right)^H$}
\right\}$ \ where

$\text{$\left(\frac{\partial \ell}{\partial \mathbf{z}}\right)^H$} =
- \frac{1}{2} \left[ \left( \frac{\mathbf g(\mathbf{\hat
z})}{\partial \mathbf{z}} \right)^H W \mathbf{e}(\mathbf{\hat z}) +
\overline{\left( \frac{\mathbf g(\mathbf{\hat z})}{\partial
\mathbf{\bar z}} \right)^H W \mathbf{e}(\mathbf{\hat z})} \, \,
\right]$; \quad $ \text{$\left(\frac{\partial \ell}{\partial
\mathbf{\bar z}}\right)^H$} = \overline{\text{$\left(\frac{\partial
\ell}{\partial \mathbf{z}}\right)^H$}}$

$\mathbf{g}(\mathbf{z})$ holomorphic:

$\mathcal{H}_{\mathbf{c}\mathbf{c}}^{\text{\tiny Newton}} =
\begin{pmatrix}
U_{\mathbf{z}\mathbf{z}} & 0 \\
0 & {U_{\mathbf{\bar z}\mathbf{\bar z}}} \end{pmatrix} -
\sum_{i=1}^m\begin{pmatrix}
V_{\mathbf{z}\mathbf{z}}^{\text{\tiny $(i)$}} & V_{\mathbf{\bar z}\mathbf{z}}^{\text{\tiny $(i)$}} \\
{V_{\mathbf{z}\mathbf{\bar z}}^{\text{\tiny $(i)$}}} &
{V_{\mathbf{\bar z}\mathbf{\bar z}}^{\text{\tiny $(i)$}}}
\end{pmatrix} = \mathcal{H}_{\mathbf{c}\mathbf{c}}^{\text{\tiny
pseudo-Gauss}}(\mathbf{\hat c}) - \sum_{i=1}^m\begin{pmatrix}
V_{\mathbf{z}\mathbf{z}}^{\text{\tiny $(i)$}} & V_{\mathbf{\bar z}\mathbf{z}}^{\text{\tiny $(i)$}} \\
{V_{\mathbf{z}\mathbf{\bar z}}^{\text{\tiny $(i)$}}} &
{V_{\mathbf{\bar z}\mathbf{\bar z}}^{\text{\tiny $(i)$}}}
\end{pmatrix}$

$V_{\mathbf{z}\mathbf{z}}^{\text{\tiny $(i)$}}$ and $V_{\mathbf{\bar
z}\mathbf{z}}^{\text{\tiny $(i)$}}$ take the simpler forms of
(\ref{eq:simpvizz}), \, $V_{\mathbf{\bar z}\mathbf{\bar
z}}^{\text{\tiny $(i)$}} =
\overline{V_{\mathbf{z}\mathbf{z}}^{\text{\tiny $(i)$}}}$, \,
$V_{\mathbf{z}\mathbf{\bar z}}^{\text{\tiny $(i)$}} =
\overline{V_{\mathbf{\bar z}\mathbf{z}}^{\text{\tiny $(i)$}}}$

$U_{\mathbf{z}\mathbf{z}}$ takes the simpler form of
(\ref{eq:simphuzz}), \, $U_{\mathbf{\bar z}\mathbf{\bar z}} =
\overline{U_{\mathbf{z}\mathbf{z}}}$

$\Delta \mathbf{z}^{\text{\tiny Newton}} = \left(
\mathcal{H}_{\mathbf{z}\mathbf{z}} - \mathcal{H}_{\mathbf{\bar
z}\mathbf{z}} \mathcal{H}_{\mathbf{\bar z}\mathbf{\bar
z}}^{-1}\mathcal{H}_{\mathbf{z}\mathbf{\bar z}}\right)^{-1} \left\{
\mathcal{H}_{\mathbf{\bar z}\mathbf{z}} \mathcal{H}_{\mathbf{\bar
z}\mathbf{\bar z}}^{-1} \text{\footnotesize $\left(\frac{\partial
\ell}{\partial \mathbf{\bar z}}\right)^H$} - \text{\footnotesize
$\left(\frac{\partial \ell}{\partial \mathbf{z}}\right)^H$}
\right\}$ \ where

$\text{$\left(\frac{\partial \ell}{\partial \mathbf{z}}\right)^H$} =
- \frac{1}{2} \left( \frac{\mathbf g(\mathbf{\hat z})}{\partial
\mathbf{z}} \right)^H W \mathbf{e}(\mathbf{\hat z})$; \quad $
\text{$\left(\frac{\partial \ell}{\partial \mathbf{\bar
z}}\right)^H$} = \overline{\text{$\left(\frac{\partial
\ell}{\partial \mathbf{z}}\right)^H$}}$
\end{itemize}

Stability generally requires positive definiteness of both
$\mathcal{H}_{\mathbf{z}\mathbf{z}}$ and its Schur complement
$\widetilde{\mathcal{H}_{\mathbf{\bar z}\mathbf{\bar z}}} = \left(
\mathcal{H}_{\mathbf{z}\mathbf{z}} - \mathcal{H}_{\mathbf{\bar
z}\mathbf{z}} \mathcal{H}_{\mathbf{\bar z}\mathbf{\bar
z}}^{-1}\mathcal{H}_{\mathbf{z}\mathbf{\bar z}}\right)$. The need to
step for positive-definiteness of the Schur complement can
significantly increase the complexity of an on-line adaptive
filtering algorithm.

When minimizing the least-squares loss function, we expect stability
to be greater when $\mathbf{g}(\mathbf{c})$ is holomorphic. This is
particularly true if $\mathbf{g}(\mathbf{c})$ is also onto and the
algorithm is convergent, as we then expect the difference between
the Newton and Gauss-Newton Hessians (and hence the difference
between the Newton and Gauss-Newton algorithms) to become negligible
asymptotically.

The Newton algorithm is known to have very fast convergence
properties, provided it can be stabilized.

\item {\bf Pseudo-Newton Algorithm.}

Applies to any smooth loss function which is bounded from below.

$\mathcal{H}_{\mathbf{c}\mathbf{c}}^{\text{\tiny
pseudo-Newton}}(\mathbf{\hat c}) =
\begin{pmatrix} \mathcal{H}_{\mathbf{z}\mathbf{z}}(\mathbf{\hat c})
& 0 \\
0 & \mathcal{H}_{\mathbf{\bar z}\mathbf{\bar z}}(\mathbf{\hat c})
\end{pmatrix}$

$Q^{\text{\tiny pseudo-Newton}}(\mathbf{\hat c}) =
\left[\mathcal{H}_{\mathbf{c}\mathbf{c}}^{\text{\tiny
pseudo-Newton}}(\mathbf{ \hat c})\right]^{-1}$

$\Delta \mathbf{c}^{\text{\tiny psedudo-Newton}} =  - Q^{\text{\tiny
pseudo-Newton}}(\mathbf{\hat c}) \left(\frac{\partial
\ell(\mathbf{\hat c})}{\partial \mathbf{c}}\right)^H$

$\Delta \mathbf{z}^{\text{\tiny pseudo-Newton}} = - \left[
\mathcal{H}_{\mathbf{z}\mathbf{z}}(\mathbf{\hat z})\right]^{-1}
\left(\frac{\partial \ell(\mathbf{\hat z})}{\partial
\mathbf{z}}\right)^H$

Application to the Least-Squares Loss Function
(\ref{eq:weightedlsloss}):
\begin{itemize}
\item[] {\small $\mathcal{H}_{\mathbf{c}\mathbf{c}}^{\text{\tiny pseudo-Newton}}
=
\begin{pmatrix} \mathcal{H}_{\mathbf{z}\mathbf{z}}(\mathbf{\hat c})
& 0 \\
0 & \mathcal{H}_{\mathbf{\bar z}\mathbf{\bar z}}(\mathbf{\hat c})
\end{pmatrix}=
\begin{pmatrix}
U_{\mathbf{z}\mathbf{z}} - \sum\limits_{i=1}^m V_{\mathbf{z}\mathbf{z}}^{\text{\tiny $(i)$}}& 0 \\
0 & {U_{\mathbf{\bar z}\mathbf{\bar z}}} - \sum\limits_{i=1}^m
{V_{\mathbf{\bar z}\mathbf{\bar z}}^{\text{\tiny $(i)$}}}
\end{pmatrix} \\ \text{\hskip0.65in} = \mathcal{H}_{\mathbf{c}\mathbf{c}}^{\text{\tiny
pseudo-Gauss}}(\mathbf{\hat c}) -
\begin{pmatrix}\sum\limits_{i=1}^m
V_{\mathbf{z}\mathbf{z}}^{\text{\tiny $(i)$}} & 0 \\
0 & \sum\limits_{i=1}^m {V_{\mathbf{\bar z}\mathbf{\bar
z}}^{\text{\tiny $(i)$}}}
\end{pmatrix}$}

$V_{\mathbf{z}\mathbf{z}}^{\text{\tiny $(i)$}}$ is given by
(\ref{eq:vizz}) and $V_{\mathbf{\bar z}\mathbf{\bar z}}^{\text{\tiny
$(i)$}} = \overline{V_{\mathbf{z}\mathbf{z}}^{\text{\tiny $(i)$}}}$.
\ $U_{\mathbf{z}\mathbf{z}}$ is given by (\ref{eq:huzz}) and
$U_{\mathbf{\bar z}\mathbf{\bar z}} =
\overline{U_{\mathbf{z}\mathbf{z}}}$

$\Delta \mathbf{c}^{\text{\tiny pseudo-Newton}} =  - Q^{\text{\tiny
pseudo-Newton}}(\mathbf{\hat c}) \left(\frac{\partial
\ell(\mathbf{\hat c})}{\partial \mathbf{c}}\right)^H$ \ where

$\left( \frac{\partial \ell}{\partial \mathbf{c}} \right)^H = -
\frac{1}{2} {G^H W \mathbf{e}} - \frac{1}{2} S\overline{G^H W
\mathbf{e}}  = \frac{1}{2} \left( {B}(\mathbf{\hat c}) + S
\overline{B}(\mathbf{\hat c}) \right)$

with \ $B(\mathbf{\hat c}) = - G(\mathbf{\hat c})^H W e(\mathbf{\hat
c})$

$\Delta \mathbf{z}^{\text{\tiny pseudo-Newton}} = - \left[
\mathcal{H}_{\mathbf{z}\mathbf{z}}(\mathbf{\hat z})\right]^{-1}
\left(\frac{\partial \ell(\mathbf{\hat z})}{\partial
\mathbf{z}}\right)^H = - \left[ U_{\mathbf{z}\mathbf{z}} -
\sum\limits_{i=1}^m V_{\mathbf{z}\mathbf{z}}^{\text{\tiny
$(i)$}}\right]^{-1} \left(\frac{\partial \ell(\mathbf{\hat
z})}{\partial \mathbf{z}}\right)^H$ \ where

$\text{$\left(\frac{\partial \ell}{\partial \mathbf{z}}\right)^H$} =
- \frac{1}{2} \left[ \left( \frac{\mathbf g(\mathbf{\hat
z})}{\partial \mathbf{z}} \right)^H W \mathbf{e}(\mathbf{\hat z}) +
\overline{\left( \frac{\mathbf g(\mathbf{\hat z})}{\partial
\mathbf{\bar z}} \right)^H W \mathbf{e}(\mathbf{\hat z})} \, \,
\right]$

$\mathbf{g}(\mathbf{z})$ holomorphic \quad {\boldmath $\Rightarrow$}

$U_{\mathbf{z}\mathbf{z}}$ takes the simpler form of
(\ref{eq:simphuzz}), \, $U_{\mathbf{\bar z}\mathbf{\bar z}} =
\overline{U_{\mathbf{z}\mathbf{z}}}$.

$V_{\mathbf{z}\mathbf{z}}^{\text{\tiny $(i)$}}$ takes the simpler
form (\ref{eq:simpvizz}), \, $V_{\mathbf{\bar z}\mathbf{\bar
z}}^{\text{\tiny $(i)$}} =
\overline{V_{\mathbf{z}\mathbf{z}}^{\text{\tiny $(i)$}}}$

$\text{$\left(\frac{\partial \ell(\mathbf{\hat z})}{\partial
\mathbf{z}}\right)^H$} = - \frac{1}{2} \left( \frac{\mathbf
g(\mathbf{\hat z})}{\partial \mathbf{z}} \right)^H W
\mathbf{e}(\mathbf{\hat z}) $

$\Delta \mathbf{z}^{\text{\tiny pseudo-Newton}} = \frac{1}{2} \left[
U_{\mathbf{z}\mathbf{z}} - \sum\limits_{i=1}^m
V_{\mathbf{z}\mathbf{z}}^{\text{\tiny $(i)$}}\right]^{-1} \left(
\frac{\mathbf g(\mathbf{\hat z})}{\partial \mathbf{z}} \right)^H W
\mathbf{e}(\mathbf{\hat z}) \\ \text{\hskip0.75in} = \left[
\left(\frac{\partial \mathbf{g}}{\partial \mathbf{z}}\right)^H \,  W
\left(\frac{\partial \mathbf{g}}{\partial \mathbf{z}}\right) -
\sum\limits_{i=1}^m \frac{\partial}{\partial \mathbf{z}}\left(
\frac{\partial g_i(\mathbf{z})}{\partial \mathbf{z}}\right)^H \,
\left[ W \mathbf{e} \, \right]_i \right]^{-1} \left( \frac{\mathbf
g(\mathbf{\hat z})}{\partial \mathbf{z}} \right)^H W
\mathbf{e}(\mathbf{\hat z})$
\end{itemize}
Stability generally requires positive definiteness of
$\mathcal{H}_{\mathbf{z}\mathbf{z}}$.

The pseudo-Newton is expected to be fast, but have a loss of
efficiency relative to the Newton algorithm. When
$\mathbf{g}(\mathbf{z})$ is holomorphic and onto, we expect good
performance as asymptotically a stabilized pseudo-Newton algorithm
will coincide with the Newton algorithm. If $\mathbf{g}(\mathbf{z})$
is nonholomorphic, the pseudo-Newton and Newton algorithms will not
coincide asymptotically, so the speed of the pseudo-Newton algorithm
is expected to always lag the Newton algorithm.
\end{enumerate}

The algorithm suggested by Yan and Fan in \cite{GY00} corresponds in
the above taxonomy to the pseudo-Newton algorithm. We see that for
obtaining a least-squares solution to the nonlinear inverse problem
$\mathbf{y} = \mathbf{g}(\mathbf{z})$, if $\mathbf{g}$ is
holomorphic, then the Yan and Fan suggestion can result in a good
approximation to the Newton algorithm. However, for nonholomorphic
least-squares inverse problems and for other types of optimization
problems (including the problem considered by Yan and Fan in
\cite{GY00}), the approximation suggested by Yan and Fan is
\emph{not} guaranteed to provide a good \emph{approximation} to the
Newton algorithm.\footnote{Such a claim \emph{might} be true.
However, it would have to be justified.} However, as we have
discussed, it \emph{does} result in an admissible generalized gradient
descent method \emph{in its own right,} and, as such, one can judge the
resulting algorithm on its own merits and in comparison with other
competitor algorithms.

\paragraph{Equality Constraints.} The classical approach to incorporating equality constraints into the problem of
optimizing a scalar cost function is via the method of Lagrange
multipliers. The theory of Lagrange multipliers is well-posed when
the objective function and constraints are \emph{real-valued}
functions of real unknown variables. Note that a  vector of $p$
\emph{complex} equality constraint conditions,
\[ \mathbf{g}(\mathbf{z}) = 0 \in \mathbb{C}^p \] is equivalent to $2p$
\emph{real} equality constraints corresponding to the conditions
\[ \text{Re} \, \mathbf{g}(\mathbf{z}) = 0 \in \mathbb{R}^p \quad \text{and}
\quad \text{Im} \, \mathbf{g}(\mathbf{z}) = 0 \in \mathbb{R}^p \, .
\]

Thus, given the problem of optimizing a real scalar-valued loss
function $\ell(\mathbf{z})$ subject to a vector of $p$ complex
equality constraints $\mathbf{h}(\mathbf{z}) = 0$, one
can construct a well-defined lagrangian as
\begin{equation}\label{eq:lagrangian1} \mathfrak{L} =
\ell(\mathbf{z}) + \lambda_R^T \, \text{Re} \,
\mathbf{g}(\mathbf{z}) + \lambda_I^T \, \text{Im} \,
\mathbf{g}(\mathbf{z}) \, , \end{equation} for \emph{real-valued}
$p$-dimensional lagrange multiplier vectors $\lambda_R$ and
$\lambda_I$.

\medskip
If we define the \emph{complex lagrange multiplier vector} $\lambda$
by \[ \lambda = \lambda_R + j\,  \lambda_I \in \mathbb{C}^p \] it is
straightforward to show that  the lagrangian (\ref{eq:lagrangian1})
can be equivalently written as
\begin{equation}\label{eq:lagrangian2}
\mathfrak{L} = \ell(\mathbf{z}) + \text{Re} \, \lambda^H
\mathbf{g}(\mathbf{z}) \, .
\end{equation}

\medskip
One can now apply the multivariate {$\mathbb{CR}$}-Calculus
developed in this note to find a stationary solution to the
Lagrangian (\ref{eq:lagrangian2}). Of course, subtle issues
involving the application of the $\mathbf{z}$, $\mathbf{c}$-complex,
and $\mathbf{c}$-real perspectives to the problem will likely arise
on a case-by-case basis.

\paragraph{Final Comments on the 2nd Order Analysis.} It is evident
 that the analysis of second-order properties of a real-valued
 function on $\mathbb{C}^n$ is much more complicated than in the
 purely real case, perhaps even dauntingly so. Thus, it is perhaps not surprising that
 very little analysis of second properties can be found in any single location in the
 literature.\footnote{That I could find. Please alert me to any
 relevant survey references that I am ignorant of.} By far, the most
 illuminating is the paper by van den Bos \cite{VDB94a}, which,
 unfortunately, is very sparse in its explanation.\footnote{Likely a
 result of page limitations imposed by the publisher.} A careful reading of van den Bos
 indicates that he is fully aware that there are two interpretations
 of $\mathbf{c}$, \emph{viz}~the real interpretation and the complex
 interpretation. \emph{This is a key insight.} As we have seen
 above, it provides a very powerful analysis and algorithm development
 tool which allows us to switch between the $\mathbf{c}$-real
 interpretation (which enables us to use the tools and insights of real
 analysis) and the $\mathbf{c}$-complex perspective (which is
 shorthand for working at the algorithm implementation level of $\mathbf{z}$ and
 $\mathbf{\bar z}$). The now-classic paper by Brandwood \cite{Br83} presents
 a development of the complex vector calculus using
 the $\mathbf{c}$-complex perspective which, although adequate for
 the development of first-order algorithms, presents greater difficulties when used as a tool
 for second order algorithm development. In this note,
 we've exploited the insights provided by van den Bos \cite{VDB94a}
 to perform a more careful analysis of second-order Newton and
 Gauss-Newton algorithms. Of course, much work remains to explore the
 analytical, structural, numerical, and implementation properties of
 these, and other second order, algorithms.

\section{Applications}

\paragraph{1. A Simple ``Nonlinear'' Least Squares Problem - I.} This is a
simple, but interesting, problem which is nonlinear in $z \in
\mathbb{C}$ yet linear in $\mathbf{c} \in \mathcal{C} \subset
\mathbb{C}^2$.

\medskip
Let $z \in \mathbb{C}$ be an unknown scalar complex quantity we wish
to estimate from multiple iid noisy measurements,
\[ y_k = s + n_k\, ,
\] $k = 1, \cdots, n$, of a scalar signal $s \in \mathbb{C}$ which
is related to $z$ via
\[ s = g(z), \quad g(z) = \alpha z + \beta \bar{z}.  \]
where $\alpha \in \mathbb{C}$ and $\beta \in \mathbb{C}$ are known
complex numbers. It is assumed that the measurement noise $n_k$ is
iid and (complex) Gaussian, $n_k \sim N(0,\sigma^2 I)$, with
$\sigma^2$ known. Note that the function $g(z)$ is both nonlinear in
$z$ (because complex conjugation is a nonlinear operation on $z$)
and nonholomorphic (nonanalytic in $z$). However, because the
problem must be linear in the underlying real space $\mathcal{R} =
\mathbf{R}^2$ (a fact which shows up in the obvious fact that the
function $g$ is linear in $\mathbf{c}$), we expect that this problem
should be exactly solvable, as will be shown to indeed be the case.

\medskip
Under the above assumptions the maximum likelihood estimate (MLE) is
found by minimizing the loss function \cite{Ka93}\footnote{The
additional overall factor of $\frac{1}{n}$ has been added for
convenience.}
\begin{eqnarray*}
\ell(z) & = & \frac{1}{2n} \sum_{k=1}^n \| y_k - g(z) \|^2 \\  & = &
\frac{1}{n} \sum_{k=1}^n
\| y_k - \alpha z- \beta \bar{z} \|^2 \\
& = & \frac{1}{2n} \sum_{k=1}^n \overline{( y_k - \alpha z- \beta
\bar{z})} (y_k
- \alpha z- \beta \bar{z} )\\
& = & \frac{1}{2n} \sum_{k=1}^n ( \bar{y}_k - \bar{\alpha} \bar{z}-
\bar{\beta} {z}) (y_k - \alpha z- \beta \bar{z} ).
\end{eqnarray*}
Note that this is a \emph{nonlinear least-squares problem} as the
function $g(z)$ is \emph{nonlinear} in $z$.\footnote{Recall that
complex conjugation is a nonlinear operation.} Furthermore, $g(z)$
is \emph{nonholomorphic} (nonanalytic in $z$). Note, however, that
although $g(z)$ is \emph{nonlinear in $z$,} it \emph{is linear} in
$\mathbf{c} = (z,\bar{z})^T$, and that as a consequence the loss
function $\ell(z) = \ell(\mathbf{c})$ has an \emph{exact} second
order expansion in $\mathbf{c}$ of the form
(\ref{eq:creal2ndorder}), which can be verified by a simple
expansion of $\ell(z)$ in terms of $z$ and $\bar z$ (see below). The
corresponding $\mathbf{c}$-complex Hessian matrix (to be computed
below) \emph{does not} have zero off-diagonal entries, which shows
that a loss function being quadratic does not \emph{alone} ensure
that $\mathcal{H}_{\mathbf{\bar z}\mathbf{z}} = 0$, a fact which
contradicts the claim made in \cite{GY00}.

\medskip
Defining the sample average of n samples $\{ \xi_1, \cdots, \xi_k
\}$ by
\[ \left< \xi \right> \triangleq \frac{1}{n}  \sum_{k=1}^n \xi_k \]
the loss function $\ell(z)$ can be expanded and rewritten as {\small
\begin{equation}\label{eq:loss1} 2\, \ell(z) = \left< \left| y \right|^2 \right> + \alpha
\bar{\beta} z^2 - \left( \alpha \left< \bar{y} \right> + \bar{\beta}
\left< y \right> \right) z + \left( \left| \alpha \right|^2 + \left|
\beta \right|^2 \right) z\bar{z} - \left( \bar{\alpha} \left< {y}
\right> + {\beta} \left< \bar{y} \right> \right) \bar{z} +
\bar{\alpha} \beta \bar{z}^2 \end{equation}} or {\small
\[
\ell(z) = \frac{1}{2} \left< \left| y \right|^2 \right>  -
\frac{1}{2} \begin{pmatrix} \alpha \left< \bar{y} \right> +
\bar{\beta} \left< y \right> & \bar{\alpha} \left< y \right> + \beta
\left< \bar{y} \right> \end{pmatrix}\binom{z}{\bar z} + \frac{1}{4}
\binom{z}{\bar z}^H \begin{pmatrix} \left| \alpha \right|^2 + \left|
\beta \right|^2 & 2 \bar{\alpha}\beta \\ 2 {\alpha}\bar{\beta} &
\left| \alpha \right|^2 + \left| \beta \right|^2\end{pmatrix}
\binom{z}{\bar z} \, . \]} Since this expansion is done using the
$\mathbf{z}$-perspective, we expect that it corresponds to a second
order expansion about the value $\mathbf{\hat z} = 0$, {\small
\begin{equation}\label{eq:loss2}
\ell(z) = \ell(0) +\frac{\partial \ell(0)}{\partial \mathbf{c}}
\mathbf{c} + \frac{1}{2} \mathbf{c}^H \mathcal{H}_{\mathbf{c}
\mathbf{c}}^{\text{\boldmath \tiny $\mathbb{C}$}}(0) \mathbf{c}
\end{equation}} with
\[ \frac{\partial \ell(0)}{\partial \mathbf{c}} = \begin{pmatrix} \frac{\partial \ell(0)}{\partial \mathbf{z}} &
\frac{\partial \ell(0)}{\partial \mathbf{\bar z}} \end{pmatrix} = -
\frac{1}{2}
\begin{pmatrix} \alpha \left< \bar{y} \right> + \bar{\beta} \left< y
\right> & \bar{\alpha} \left< y \right> + \beta \left< \bar{y}
\right> \end{pmatrix} \] and
\[ \mathcal{H}_{\mathbf{c} \mathbf{c}}^{\text{\boldmath \tiny
$\mathbb{C}$}}(0) =\frac{1}{2} \begin{pmatrix} \left| \alpha
\right|^2 +
\left| \beta \right|^2 & 2 \, \bar{\alpha}\beta \\
2\, {\alpha}\bar{\beta} & \left| \alpha \right|^2 + \left| \beta
\right|^2 \end{pmatrix} . \] And indeed this turns out to be the
case. Simple differentiation of (\ref{eq:loss1}) yields,
\[ \frac{\partial \ell(z)}{\partial z} = \alpha \bar{\beta} z + \frac{1}{2} \left( \left| \alpha \right|^2
+ \left| \beta \right|^2 \right) \bar{z} - \frac{1}{2} \left( \alpha
\left< \bar{y} \right> + \bar{\beta} \left< {y} \right> \right) \]
\[ \frac{\partial \ell(z)}{\partial \bar{z}} =  \bar \alpha {\beta}
\bar z + \frac{1}{2} \left( \left| \alpha \right|^2 + \left| \beta
\right|^2 \right) {z} - \frac{1}{2} \left( \bar{\alpha} \left< {y}
\right> + {\beta} \left< \bar{y} \right> \right) \] which evaluated
at zero give the linear term in the quadratic loss function, and
further differentiations yield,
\[ \mathcal{H}_{\mathbf{c} \mathbf{c}}^{\text{\boldmath \tiny
$\mathbb{C}$}}(z) = \begin{pmatrix} \mathcal{H}_{{z}{z}}
& \mathcal{H}_{{\bar z}{z}} \\
\mathcal{H}_{{z}{\bar z}} & \mathcal{H}_{{\bar z}{\bar z}}
\end{pmatrix} = \frac{1}{2} \begin{pmatrix} \left| \alpha \right|^2
+
\left| \beta \right|^2 & 2 \, \bar{\alpha}\beta \\
2\, {\alpha}\bar{\beta} & \left| \alpha \right|^2 + \left| \beta
\right|^2 \end{pmatrix}  \] which is independent of $z$. Note that,
as expected,
\[ \frac{\partial \ell(z)}{\partial \bar z}= \overline{\frac{\partial \ell(z)}{\partial
z}} \,. \]

If we set the two partial derivatives to zero, we obtain two
stationarity equations for the two stationary quantities $z$ and
$\bar z$. Solving for $z$ then yields the least-squares estimate of
$z$,\footnote{Note that this answer reduces to the obvious solutions
for the two special cases $\alpha = 0$ and $\beta = 0$.}
\[ \hat{z}_{\text{\tiny opt}} = \frac{1}{\left| \alpha \right|^2 - \left| \beta \right|^2}
\left( \bar{\alpha} \left< y \right> - \beta \left< \bar{y} \right>
\right)\ . \] This solution can also be obtained by completing the
square on (\ref{eq:loss2}) to obtain
\[ \mathbf{\hat c}_{\text{\tiny opt}} = - \left(\mathcal{H}_{\mathbf{c} \mathbf{c}}^{\text{\boldmath \tiny
$\mathbb{C}$}}\right)^{-1} \left(\frac{\partial \ell(0)}{\partial
\mathbf{c}}\right)^H\]

An obvious necessary condition for the least-squares solution to
exist is that
\[ \left| \alpha \right|^2 \ne \left| \beta \right|^2 . \]
The solution will be a global\footnote{Because the Hessian is
independent of $z$.} minimum if the Hessian matrix is positive
definite. This will be true if the two leading principal minors are
strictly positive, which is true if and only if, again, $\left|
\alpha \right|^2 \ne \left| \beta \right|^2$. Thus, if $\left|
\alpha \right|^2 \ne \left| \beta \right|^2$ the solution given
above is a global minimum to the least squares problem.

\medskip
The condition $\left| \alpha \right|^2 = \left| \beta \right|^2$
corresponds to \emph{loss of identifiability} of the model \[ g(z) =
\alpha z + \beta \bar{z}\, .
\] To see this, first note that to identify a complex number is
equivalent to identifying both the real and imaginary parts of the
number. If either of them is unidentifiable, then so is the number.

\medskip
Now note that the condition $\left| \alpha \right|^2 = \left| \beta
\right|^2$ says that $\alpha$ and $\beta$ have the same magnitude,
but, in general, a different phase. If we call the phase difference
$\phi$, then the condition $\left| \alpha \right|^2 = \left| \beta
\right|^2$ is equivalent to the condition
\[ \alpha = e^{j \phi} \beta \, , \]
which yields
\[ g(z) = e^{j \phi} \beta z + \beta \bar{z} = e^{j \frac{\phi}{2}} \beta
\left( e^{j \frac{\phi}{2}} z  + e^{- j \frac{\phi}{2}} \bar{z}
\right) = e^{j \frac{\phi}{2}} \beta \left( e^{j \frac{\phi}{2}} z +
\overline{e^{j \frac{\phi}{2}}{z}} \right) = e^{j \frac{\phi}{2}}
\beta \,  \text{Re} \left\{ e^{j \frac{\phi}{2}} z \right\}.
\]
Thus, it is evident that the imaginary part of $e^{j \frac{\phi}{2}}
z$ is unidentifiable, and thus the complex number $e^{j
\frac{\phi}{2}} z$ itself is unidentifiable. And, since
\[ z = e^{-j \frac{\phi}{2}} \left(
 e^{j \frac{\phi}{2}} z  \right)= e^{-j \frac{\phi}{2}}
\left( \text{Re} \left\{ e^{j \frac{\phi}{2}} z \right\} + j \,
\text{Im} \left\{ e^{j \frac{\phi}{2}} z\right\} \right) , \] it is
obvious that $z$ is unidentifiable.

\medskip
Note for the simplest case of $\alpha = \beta$ ($\phi = 0$), we have
\[ g(z) = \alpha z + \alpha \bar{z} = \alpha \, \text{Re} \left\{ z \right\} \]
in which case $\text{Im} \left\{ z \right\}$, and hence $z$, is
unidentifiable.

\paragraph{2. A Simple ``Nonlinear'' Least Squares Problem - II.}
The ``nonlinearity'' encountered in the previous example, is in a
sense ``bogus'' and is not a nonlinearity at all, at least when
viewed from the $\mathbf{c}$-real perspective. Not surprisingly
then, we were able to compute an exact solution. Here, we will
briefly look at the Newton and Gauss-Newton algorithms applied to
the simple problem of Example 1.

\medskip
In the previous example, we computed the Newton Hessian of the
least-squares loss function (\ref{eq:loss1}). The difference between
the Newton and Gauss-Newton algorithm resides in the difference
between the Newton Hessian and the Gauss-Newton Hessian. To compute
the Gauss-Newton Hessian, note that
\[ y = g(\mathbf{c}) = (\alpha \ \beta)\binom{z}{\bar z} = G
\mathbf{c} \] and therefore (since the problem is linear in
$\mathbf{c}$) we have the not surprising result that
\[ G \Delta \mathbf{c} = \frac{\partial g(\mathbf{c})}{\partial \mathbf{c}}
\Delta \mathbf{c} \] with
\[ G = (\alpha \ \beta) \, . \]
In this example, the least-squares weighting matrix is $W = I$ and
we have
\[ G^H W G = G^H G = \binom{\bar{\alpha}}{\bar{\beta}} (\alpha \
\beta) = \begin{pmatrix} \left| \alpha \right|^2 & \bar{\alpha}
\beta \\ \bar{\beta} \alpha & \left| \beta \right|^2 \end{pmatrix}
\]
which is seen to be independent of $\mathbf{c}$.
 From
(\ref{eq:gaussnewtonhessian}), we construct the Gauss-Newton Hessian
as  \[ \small \mathcal{H}_{\mathbf{c}\mathbf{c}}^{\text{\tiny
Gauss}} = \mathbf{P}\left(G^H G \right) = \frac{\begin{pmatrix}
\left| \alpha \right|^2 & \bar{\alpha} \beta \\ \bar{\beta} \alpha &
\left| \beta \right|^2 \end{pmatrix} + S \,
\overline{\begin{pmatrix} \left| \alpha \right|^2 & \bar{\alpha}
\beta \\ \bar{\beta} \alpha & \left| \beta \right|^2
\end{pmatrix}}\, S }{2} = \frac{1}{2}
\begin{pmatrix} \left| \alpha \right|^2 +
\left| \beta \right|^2 & 2 \, \bar{\alpha}\beta \\
2\, {\alpha}\bar{\beta} & \left| \alpha \right|^2 + \left| \beta
\right|^2 \end{pmatrix} = \mathcal{H}_{\mathbf{c}
\mathbf{c}}^{\text{\boldmath \tiny $\mathbb{C}$}}
\]
showing that for this simple example the Newton and Gauss-Newton
Hessians are the same, and therefore \emph{the Newton and
Gauss-Newton algorithms are identical.} As seen from Equations
(\ref{eq:syminhh}) and (\ref{eq:rawhessian}), this is a consequence
of the fact that $g(\mathbf{c})$ is linear in $\mathbf{c}$ as then
the matrix of second partial derivatives of $g$ required to compute
the difference between the Newton and Gauss-Newton algorithms
vanishes
\[
A_{\mathbf{c}\mathbf{c}}(g) \triangleq \frac{\partial}{\partial
\mathbf{c}} \left( \frac{\partial g}{\partial \mathbf{c}}\right)^H =
0 . \]

From the derivatives computed in the previous example, we can
compute $\left(\frac{\partial \ell(\mathbf{\hat c})}{\partial
\mathbf{c}}\right)^H$ as
\[ \left(\frac{\partial \ell(\mathbf{\hat c})}{\partial \mathbf{c}}\right)^H
= \begin{pmatrix}{\left(\frac{\partial
\ell(\mathbf{\hat c})}{\partial \mathbf{z}}\right)^H} \\
{\left(\frac{\partial \ell(\mathbf{\hat c})}{\partial \mathbf{\bar
z}}\right)^H}
\end{pmatrix} = \begin{pmatrix}\left(\frac{\partial \ell(\mathbf{0})}{\partial
\mathbf{z}}\right)^H \\ \left(\frac{\partial
\ell(\mathbf{0})}{\partial \mathbf{\bar z}}\right)^H \end{pmatrix} +
\frac{1}{2}
\begin{pmatrix} \left| \alpha \right|^2 +
\left| \beta \right|^2 & 2 \, \bar{\alpha}\beta \\
2\, {\alpha}\bar{\beta} & \left| \alpha \right|^2 + \left| \beta
\right|^2 \end{pmatrix}\binom{\hat z}{\hat{\bar z}}
\]
or
\[ \left(\frac{\partial \ell(\mathbf{\hat c})}{\partial \mathbf{c}}\right)^H
= \left(\frac{\partial \ell(0)}{\partial \mathbf{c}}\right)^H +
\mathcal{H}_{\mathbf{c} \mathbf{c}}^{\text{\boldmath \tiny
$\mathbb{C}$}}\mathbf{\hat c} . \] The optimal update in the Newton
algorithm is therefore given by \[ \widehat{\Delta \mathbf{c}} = -
\left(\mathcal{H}_{\mathbf{c} \mathbf{c}}^{\text{\boldmath \tiny
$\mathbb{C}$}}\right)^{-1}\left(\frac{\partial \ell(\mathbf{\hat
c})}{\partial \mathbf{c}}\right)^H = - \left(\mathcal{H}_{\mathbf{c}
\mathbf{c}}^{\text{\boldmath \tiny
$\mathbb{C}$}}\right)^{-1}\left(\frac{\partial \ell(0)}{\partial
\mathbf{c}}\right)^H  - \mathbf{\hat c} = \mathbf{\hat
c}_{\text{\tiny opt}} - \mathbf{\hat c} \, . \] The update step in
the Newton algorithm is given by
\[ \mathbf{\hat c}_{\text{\tiny new}} = \mathbf{\hat c} + \alpha \widehat{\Delta
\mathbf{c}} \, . \] If we take the ``Newton stepsize'' $\alpha = 1$,
we obtain
\[ \mathbf{\hat c}_{\text{\tiny new}} = \mathbf{\hat c} +  \widehat{\Delta
\mathbf{c}} = \mathbf{\hat c} + \mathbf{\hat c}_{\text{\tiny opt}} -
\mathbf{\hat c} = \mathbf{\hat c}_{\text{\tiny opt}} \] showing that
we can attain the optimal solution in only one update step. For the
real case, it is well-known that the Newton algorithm attains the
optimum in one step for a quadratic loss function. Thus our result
is not surprising given that the problem is a linear least-squares
problem in $\mathbf{c}$.

\medskip
Note that the off-diagonal elements of the constant-valued  Hessian
$\mathcal{H}_{\mathbf{c} \mathbf{c}}^{\text{\boldmath \tiny
$\mathbb{C}$}}$ \emph{are never zero} and generally are not small
relative to the size of the diagonal elements of
$\mathcal{H}_{\mathbf{c} \mathbf{c}}^{\text{\boldmath \tiny
$\mathbb{C}$}}$. This contradicts the statement made in \cite{GY00}
that for a quadratic loss function, the diagonal elements must be
zero.\footnote{It is true, as we noted above, that for the quadratic
loss function associated with a holomorphic nonlinear inverse
problem the off-diagonal elements of the Hessian are zero. However,
the statement is not true in general.} However, the pseudo-Newton
algorithm proposed in \cite{GY00} will converge to the correct
solution when applied to our problem, but at a slower convergent
rate than the full Newton algorithm, which is seen to be capable of
providing one-step convergence. We have a trade off between
complexity (the less complex pseudo-Newton algorithm versus the more
complex Newton algorithm) versus speed of convergence (the slower
converging pseudo-Newton algorithm versus the fast Newton
algorithm).

\paragraph{3. The Complex LMS Algorithm.} Consider the problem of determining the complex \emph{vector}
parameter $a \in \mathbb{C}^n$ which minimizes the following
generalization of the loss function (\ref{eq:e0}) to the vector
parameter case,
\begin{equation}\label{eq:vectorcase}
    \ell(a) = \text{E}\left\{ \left| e_k \right|^2 \right\} , \qquad
    e_k = \eta_k - a^H \xi_k ,
\end{equation}
for $\eta_k \in \mathbb{C}$ and $\xi_k \in \mathbb{C}^n$. We will
assume throughout that the parameter space is Euclidean so that
$\Omega_a = I$. The cogradient of $\ell(a)$ with respect to the
unknown parameter vector $a$ is given by
\[ \frac{\partial \ }{\partial a}\,  \ell(a)  = \text{E}\left\{ \frac{\partial \ }{\partial a}\,  \left| e \right|^2
\right\} \, . \]

To determine the cogradient of
\[ \left| e_k \right|^2 = \bar{e}_k e_k =
e_k \bar{e}_k = (\eta_k - a^H \xi_k)\overline{(\eta_k - a^H \xi_k)}
\] note  that
\[ \bar{e}_k = \overline{(\eta_k - a^H \xi_k)} = (\bar{\eta}_k -  \xi_k^H a) \]
and that $e_k = (\eta_k - a^H \xi_k)$ is independent of $a$. Then we
have
\begin{eqnarray*}
\frac{\partial \ }{\partial a}\, e_k \bar{e}_k & = & e_k \, \frac{\partial \ }{\partial a}\, (\bar{\eta}_k -  \xi_k^H a) \\
& = &   - e_k \frac{\partial \ }{\partial a}\, \xi_k^H a \\
& = & - \, e_k \, \xi_k^H \, .
\end{eqnarray*}
The gradient of $| e_k |^2 = e_k \bar{e}_k$ is given by
\[ \nabla_a e_k \bar{e}_k = \left(\frac{\partial \ }{\partial a}\, e_k
\bar{e}_k\right)^H = - \left( e_k \, \xi_k^H \right)^H = - \xi_k
\bar{e}_k \, . \]

\medskip
Thus, we readily have that the gradient (direction of steepest
ascent) of the loss function $\ell(a) = \text{E}\left\{ \left| e_k
\right|^2 \right\}$ is
\[  \nabla_{a} \, \ell(a) = - \text{E} \left\{ \xi_k \bar{e}_k
\right\} = - \text{E} \left\{ \xi_k \, (\bar{\eta}_k - \xi_k^H a)
\right\} \, .
\]
If we set this (or the cogradient) equal to zero to determine a
stationary point of the loss function we obtain the standard
Wiener-Hopf equations for the MMSE estimate of $a$.\footnote{Which,
as mentioned earlier, can also be obtained from the orthogonality
principle or completing the square. Thus, if the Wiener-Hopf
equations are our only goal there is no need to discuss complex
derivatives at all.  It is only when a direction of steepest descent
is needed in order to implement an on-line adaptive descent-like
algorithm that the need for the extended or conjugate derivative
arises.}

\medskip
Alternatively, if we make the \emph{instantaneous
stochastic-gradient approximation,}
\[ \nabla_{a} \ell(a) \approx \widehat \nabla_{a} \ell(\widehat a_k)
\triangleq \nabla_{a} | e_k |^2 = - \xi_k \bar{e}_k = \xi_k \left(
\bar{\eta}_k - \xi_k^H \widehat a_k\right) \, ,
\] where $\widehat a_k$ is a current estimate of the MMSE value of
$a$ and $-\nabla_a \ell(a)$ gives the direction of steepest descent
of $\ell(a)$, we obtain the standard LMS on-line stochastic
gradient-descent algorithm for learning an estimate of the complex
vector $a$,
\begin{eqnarray*}
\widehat a_{k+1} & = & \widehat a_k - \alpha_k \widehat
\nabla_{a} \ell(\widehat a_k) \\
& = & \widehat a_k + \alpha_k \xi_k \bar{e}_k \\
& = & \widehat a_k + \alpha_k \xi_k \left( \bar{\eta}_k - \xi^H_k
\widehat a_k \right)
\\ & = & \left( I - \alpha_k \xi_k \xi^H_k \right) \widehat a_k +
\alpha_k \xi_k \bar{\eta}_k \, .
\end{eqnarray*}

Thus, we have easily derived the complex LMS algorithm,
\begin{equation}
\text{\sf Complex LMS Algorithm:} \qquad \widehat a_{k+1} = \left( I
- \alpha_k \xi_k \xi^H_k \right) \widehat a_k + \alpha_k \xi_k
\bar{\eta}_k \, .
\end{equation}


\end{document}